\documentclass[preprint,12pt,authoryear,a4paper]{elsarticle}

\usepackage[utf8]{inputenc}
\usepackage[T1]{fontenc}
\usepackage{lmodern}
\usepackage{amsmath,amssymb,amsfonts}
\usepackage{booktabs}
\usepackage{array}
\usepackage{graphicx}
\usepackage{rotating}
\usepackage{pdflscape}
\usepackage{subcaption}
\usepackage{xcolor}
\usepackage{enumitem}
\usepackage{hyperref}
\hypersetup{hidelinks}
\usepackage{algorithm}
\usepackage{algpseudocode}
\usepackage{longtable}
\usepackage{multirow}
\usepackage{siunitx}
\usepackage{upgreek}
\usepackage{bm}
\usepackage{lineno}
\usepackage{geometry}
\usepackage{setspace}
\usepackage[section]{placeins}
\let\oldsubsection\subsection
\renewcommand{\subsection}{\FloatBarrier\oldsubsection}

\journal{Energy Conversion and Management}

\begin{document}

\begin{frontmatter}

\title{Exploring Low-Regret Energy System Transition Pathways with Modeling to Generate Alternatives}

\author[inst1]{Hendrik Schricker}
\author[inst1]{Maria Louise Schmitt}
\author[inst1]{Philipp Kleusberg}
\author[inst1]{Valentin Dinges}
\author[inst1]{Benedikt Schuler}
\author[inst1]{Frederike Kuperjans}
\author[inst1]{Boyung Jürgens}
\author[inst1]{Niklas von der Assen\corref{cor1}}
\cortext[cor1]{Corresponding author. E-mail: niklas.vonderassen@ltt.rwth-aachen.de}

\affiliation[inst1]{organization={Institute of Technical Thermodynamics, RWTH Aachen University},
addressline={Schinkelstr. 8},
postcode={52062},
city={Aachen},
country={Germany}}

\begin{abstract}
Long-term planning of energy system transitions is complex due to interdependent investment decisions and uncertain market and policy conditions. Stochastic optimization models account for this uncertainty by identifying a single optimal hedging strategy. However, the flexibility in technology choice and investment timing, i.e., the transition flexibility, usually remains unexplored. Modeling to generate alternatives (MGA) explores near-optimal solution spaces but has not been coupled with stochastic optimization to quantify the transition flexibility.
We explore the low-regret pathway space, i.e., all transition pathways that stay within a chosen regret tolerance in every scenario. We introduce a termination criterion for MGA that bounds the exploration error in high-dimensional solution spaces and measure how far each pathway decision, once taken, narrows the capacity ranges of the others. We apply the method to an industrial energy system in Germany. The pathway that minimizes the worst-case regret keeps it at 2.2\%, compared with up to 30.9\% for perfect-foresight pathways. The transition flexibility increases most at low regret tolerances. At 3.5\%, the capacity ranges are on average about 60\% as wide as at the highest regret tolerance of 10\%. A large share of these capacity ranges is surplus capacity, i.e., capacity installed in addition to, not instead of, other capacity. Surplus capacity is largely not operated but, once installed, narrows the other capacity ranges more strongly than operated capacity. Exploring the low-regret pathway space, rather than reporting a single optimal pathway, therefore is the more informative basis for planning energy system transitions under uncertainty.
\end{abstract}


\end{frontmatter}

\section{Introduction}
\label{sec:introduction}

Energy system transitions require long-lived investments in renewable generation, storage, and sector-coupling technologies. These investments are capital-intensive and exposed to uncertain future energy price and policy conditions. Transition planning is therefore not only a problem of minimizing cost, but of identifying pathways that remain robust and flexible under uncertain futures. We refer to the flexibility in technology choice and investment timing across such pathways as their \emph{transition flexibility}. Transition flexibility describes the agency that remains for decision makers despite the uncertainty in future conditions. Quantifying transition flexibility requires methods that account for uncertainty in the transition planning.

In the literature, a large variety of methods supports energy system planning under uncertainty \citep{Paredes-Vergara2024,Yue2018}. \citet{Yue2018} distinguish parametric uncertainty, i.e., uncertainty in input parameters such as future energy prices, technology costs, or weather conditions, from structural uncertainty in the model formulation itself, which modeling to generate alternatives (MGA) addresses. Methods for parametric uncertainty either consider the uncertainty within the optimization model or, like Monte Carlo analysis, sample parameter sets and solve the model for each sample. Stochastic programming models the uncertainty endogenously in the optimization and identifies solutions that hedge best against uncertainty depending on the chosen decision criterion. For example \citet{Zhou2013} formulate a two-stage stochastic program for the design of distributed energy systems under uncertain energy demands and renewable energy availability. \citet{Mavromatidis2018} study the influence of decision criteria in two-stage stochastic programming for the design of distributed energy systems. Their results show that risk-aware criteria such as minimax regret, which minimizes the largest cost increase relative to the perfect-foresight optimum of each scenario, can lead to more robust solutions compared to risk-neutral criteria such as expected cost without being overly conservative. Consequently, approaches to identify low-regret solutions in energy planning problems are increasingly researched, e.g., by \citet{Wiest2026}.

Several studies consider transition pathway optimization under uncertainty. \citet{Bohlayer2021}, \citet{Frischmuth2024}, and \citet{McDonald2024} formulate and compare two-stage and multi-stage stochastic programs for transition pathways of distributed and national energy systems under parametric uncertainty. Beyond parametric uncertainty, \citet{Coppitters2026} study the robustness of national transition pathways to unexpected events by solving an ensemble of scenarios with perfect and myopic foresight.

Independent of how uncertainty is treated, energy system optimization studies often report only a single cost-optimal solution. \citet{Keepin1984}, \citet{Trutnevyte2016}, and \citet{FossasTenas2026} show that cost-optimal solutions can be highly sensitive to hidden assumptions and differ substantially from real-world transitions. \citet{DeCarolis2016} argue that energy system optimization models should identify feasible decision alternatives because not all preferences, e.g., social and political preferences, can be fully represented in a cost function. Consequently, the cost-optimal solution is not necessarily the most useful decision alternative in practice.

MGA overcomes the focus on single optimal solutions by relaxing the optimality criterion within a user-defined tolerance and searching for solutions that remain near-optimal while differing in decision-relevant variables. \citet{Brill1979} and \citet{Brill1982} formalize this idea by identifying maximally different alternatives in optimization problems. \citet{Lombardi2025} review recent methodological advancements in MGA and categorize them by their focus: efficient computation of near-optimal solution spaces \citep[e.g.,][]{Lau2024b}, spatial diversity for large-scale energy systems \citep[e.g.,][]{Lombardi2020}, intuitive results communication \citep[e.g.,][]{Baader.11.112023}, and the combination with methods for parametric uncertainty \citep[e.g.,][]{Li2017,Neumann2023,Grochowicz2023}.

To incorporate parametric uncertainty in MGA, \citet{Li2017} combine Monte Carlo analysis with MGA for national electricity system pathways to generate near-optimal alternatives for each sampled parameter set. However, each alternative is near-optimal only for its own parameter sample. Whether an alternative remains near-optimal in the other samples is not assessed. \citet{Neumann2023} extend this idea with surrogate models to quantify the probability that a given design remains near-optimal across parameter samples. \citet{Grochowicz2023}, \citet{vanGreevenbroek2025b}, and \citet{vanGreevenbroek2025a} intersect the near-optimal solution spaces of multiple scenarios to identify robust solutions, i.e., solutions that remain near-optimal in every scenario. Exploring the near-optimal solution space of each scenario independently reduces the computational effort for the MGA optimization runs, as each run considers only a single scenario instead of all scenarios simultaneously. However, the exploration is not focused on those regions that lie in the intersection of all spaces. As a result, many runs may be spent on solutions that are not robust in all scenarios, and many solutions that ultimately define the robust, near-optimal solution space may be missed. In the worst case, the intersection might not even be constructible \citep{vanGreevenbroek2024}. Furthermore, the original method of \citet{Grochowicz2023} requires convex hull computations that become intractable in the high-dimensional near-optimal solution spaces typical for transition pathway problems. All abovementioned approaches share one limitation: robustness is not embedded in the construction of the near-optimal solution space itself, which leaves a gap in MGA research.

Applications of MGA to transition pathways remain scarce. \citet{Price2017} compare near-optimal transition pathways for the global energy system and \citet{Sinha2024} generate diverse near-optimal US decarbonization pathways. However, both neglect robustness to parametric uncertainty. Furthermore, both studies limit the exploration to pathway decisions cumulated over the planning horizon, i.e., they do not explicitly explore the flexibility in the timing of these decisions. The study of \citet{vanGreevenbroek2025a} extends the consideration of near-optimal transition pathways by incorporating parametric uncertainty. They intersect near-optimal solution spaces of individual scenarios to identify robust transition pathways for European green hydrogen deployment. The flexibility they explore is the range of a single pathway variable, the annual green hydrogen production, that remains near-optimal in all scenarios for each modeled year. Due to the computational limitation of the intersection-based method described above, their approach would not be directly applicable to the joint exploration of many pathway variables in high-dimensional near-optimal solution spaces. Recently, \citet{Kchaou2026} combine MGA with a multi-armed bandit approach to identify robust early-stage energy transition decisions under unexpected events. However, their work assesses robustness only after generating near-optimal alternatives. Thus, robustness is again not embedded in the construction of the near-optimal space.

\citet{Lombardi2025} identify the application of MGA to transition pathways and its coupling with methods for parametric uncertainty as key open challenges in MGA research. The literature lacks a method that combines both dimensions: a near-optimal transition pathway space that accounts for uncertainty by construction and that is explicit in both the choice of technologies and the timing of the pathway decisions.

Our own prior work explored near-optimal and near-Pareto-optimal design spaces of industrial energy systems, but only for single-year design decisions without parametric uncertainty \citep{Schricker2025b,Schricker2026}. This paper characterizes the near-optimal space of transition pathways as a low-regret pathway space around an optimal pathway that minimizes the worst-case regret across all considered scenarios. Bounding the scenario-wise regret instead of the cost transfers the near-optimality tolerance of MGA to the scenario set of the stochastic program: a transition pathway is low-regret only if it is near-optimal in every scenario. The paper makes three contributions. First, we define the low-regret pathway space by scenario-wise regret constraints in a two-stage stochastic program. Second, in contrast to approaches that intersect near-optimal solution spaces after exploration, we explore the intersection directly. For this purpose, we introduce an exploration method that directs each search toward the largest gaps between the inner and outer approximation of the pathway space, together with a termination criterion that bounds the exploration error. Both, the exploration method and the termination criterion remain tractable in high-dimensional near-optimal solution spaces. Third, we quantify the transition flexibility of low-regret pathway spaces and assess how the regret tolerance and single pathway decisions influence it. We demonstrate the method on the transition planning of an industrial energy system in Germany under uncertain future energy and CO$_2$ certificate prices, as industrial energy systems combine a high need for decarbonization with high exposure to future uncertainty \citep{Bataille2018,Brodnicke2025}. With the explored pathway space, we answer three research questions:
\begin{enumerate}[label=\textbf{RQ\arabic*:}, ref=RQ\arabic*, leftmargin=*, itemsep=2pt]
    \item\label{rq:pathway-space} How much transition flexibility exists in low-regret pathway spaces?
    \item\label{rq:cost-of-flexibility} How much transition flexibility is gained with a higher regret tolerance?
    \item\label{rq:residual-flexibility} How strongly do single pathway decisions constrain the transition flexibility?
\end{enumerate}
\section{Methods}
\label{sec:methods}

\begin{figure}[htbp]
\centering
\includegraphics[width=\textwidth]{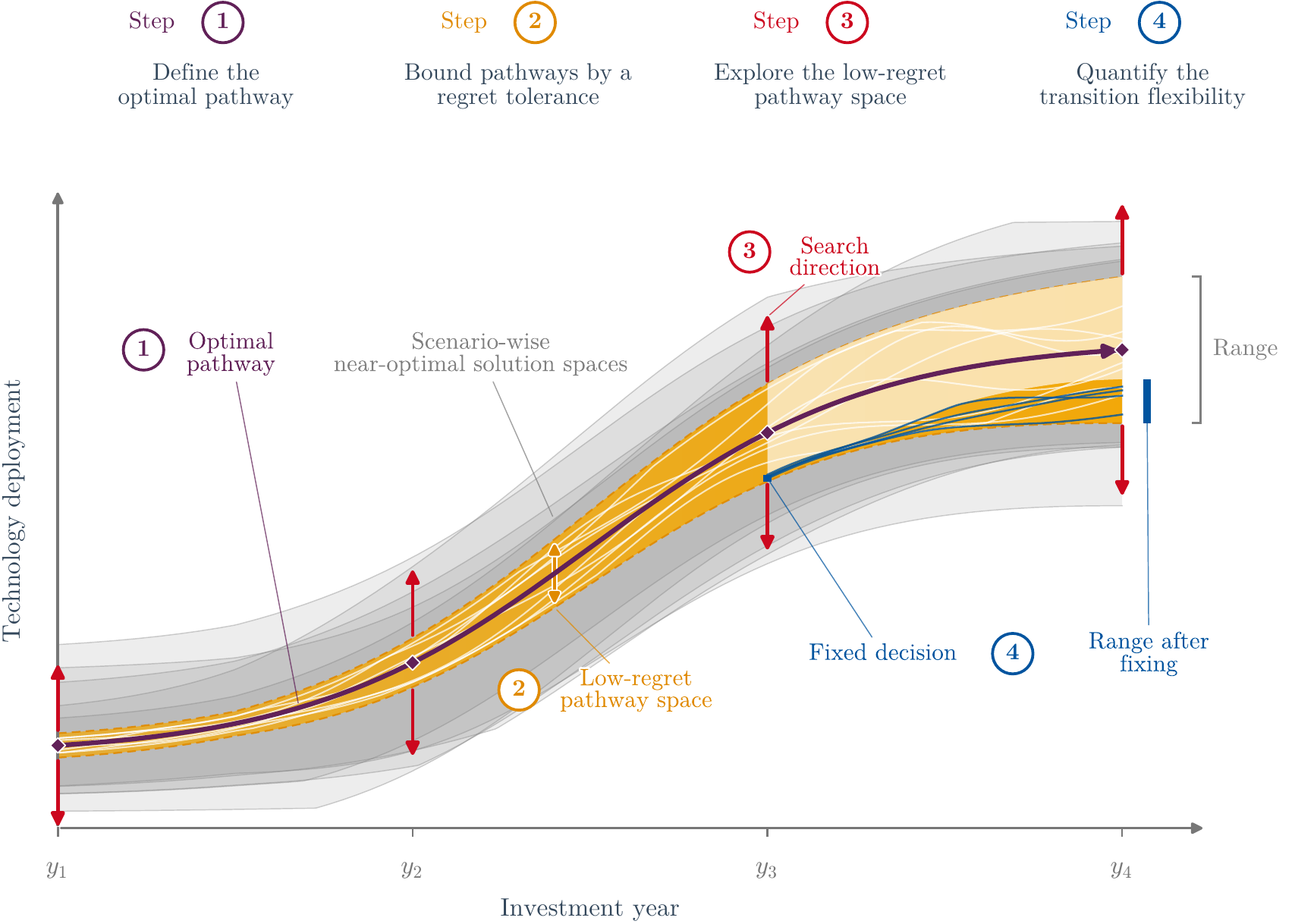}
\caption{Schematic overview of the proposed method for exploring low-regret transition pathways under uncertainty. In \textbf{Step~1}, we solve a two-stage stochastic energy system optimization model to identify the optimal pathway, which minimizes the worst-case regret across all scenarios (purple line). In \textbf{Step~2}, we bound the worst-case regret by a regret tolerance, which defines the intersection of the scenario-wise spaces (grey bands), i.e., the low-regret pathway space (orange band). In \textbf{Step~3}, we explore this space by optimizations along search directions (red arrows), which yield explored pathways (white lines). In \textbf{Step~4}, we fix one pathway decision (blue bar) and quantify the transition flexibility as the range of the pathways consistent with this decision (blue lines), compared with the full range (grey bracket).}
\label{fig:method-overview}
\end{figure}

Our method is designed to support investment decisions for energy system transitions under uncertainty. Figure~\ref{fig:method-overview} gives an overview of the method in four steps, which correspond to the four subsections of this section. In Step~1, we formulate a two-stage stochastic energy system optimization model and define the optimal pathway that minimizes the worst-case regret across all scenarios (Section~\ref{subsec:model}). In Step~2, we define the low-regret pathway space (Section~\ref{subsec:lowregret}), and in Step~3, we explore this space by optimizations along search directions (Section~\ref{subsec:exploration}). In Step~4, we quantify the transition flexibility within the low-regret pathway space (Section~\ref{subsec:decision-support}).

\subsection{Two-stage stochastic transition model and optimal pathway}
\label{subsec:model}

We first formulate the two-stage stochastic energy system optimization model that underlies all optimization problems in this study (Section~\ref{subsubsec:model-formulation}). Based on this model, we then define the optimal pathway as the transition pathway that minimizes the worst-case regret across all scenarios (Section~\ref{subsubsec:minimax}), corresponding to Step~1 in Figure~\ref{fig:method-overview}.

\subsubsection{Model formulation}
\label{subsubsec:model-formulation}

Let $\mathcal{Y}$ denote the set of investment years and $\Omega$ the set of scenarios representing the uncertain input parameters such as energy prices or CO$_2$ certificate prices. The modeled energy system consists of conversion and storage components $c\in \mathcal{C}$, which convert and store energy carriers.

We formulate the transition optimization problems as a two-stage stochastic linear program with decision vector
\begin{equation}
\boldsymbol{x}=\left[\boldsymbol{d},\{\boldsymbol{o}_\omega\}_{\omega\in\Omega}\right]\in \mathcal{X} .
\label{eq:decision-vector}
\end{equation}
The first-stage investment vector $\boldsymbol{d}=\{d_{c,y}\}_{c\in \mathcal{C},\,y\in \mathcal{Y}}$ collects the investment decision $d_{c,y}$ for each component $c$ and year $y$. We assume the investment vector to be non-anticipative, i.e., the same investments apply across all scenarios. Furthermore, we assume that uncertainty enters the model only in the second stage, i.e., the scenarios affect operational parameters and costs but not the investment stage. The second-stage decision vector $\boldsymbol{o}_\omega$ collects the scenario-dependent operational decisions for all years and all time steps within a year: energy conversion flows, storage charging and discharging, and external purchases and sales. The second-stage decision vector may adapt to the realized scenario $\omega$. For brevity, we refer to the full decision vector $\boldsymbol{x}$ as a transition pathway, comprising the investments together with their scenario-wise operation.

The feasible set $\mathcal{X}$ collects all decision vectors that satisfy the constraints of the stochastic model. The underlying transition model follows the formulation of the SecMOD framework \citep{Reinert2022}, which we cast as a two-stage stochastic program for design under uncertainty following \citet{Schricker2026b}.

The total annualized cost of a transition pathway in scenario $\omega$ is
\begin{equation}
\mathit{TAC}_\omega(\boldsymbol{x})=\sum_{y\in \mathcal{Y}}\alpha_y\, \mathit{COST}_{y,\omega}(\boldsymbol{x}),
\qquad
\alpha_y=(1+i)^{-(y-y_{\mathrm{ref}})},
\label{eq:tac}
\end{equation}
where $\mathit{COST}_{y,\omega}(\boldsymbol{x})$ collects the investment, maintenance, and operation costs, the procurement costs of energy carriers, and the revenues from sales of energy carriers in year $y$ and scenario~$\omega$. The discount factor $\alpha_y$ is computed with the interest rate $i$ and reference year $y_{\mathrm{ref}}$, which we set to the first year of the optimization horizon. We model all constraints and cost terms as linear functions of continuous decision variables. Consequently, the feasible set $\mathcal{X}$ is a polyhedron and all considered optimization problems are linear programs.

\subsubsection{Optimal pathway}
\label{subsubsec:minimax}

For each scenario $\omega$, we first compute the optimal investment and operational decisions as if the scenario were known in advance, i.e., the perfect-foresight pathway. The perfect-foresight optimization problem is:
\begin{equation}
\mathit{TAC}_\omega^*=
\min_{\boldsymbol{x}_\omega} \mathit{TAC}_\omega(\boldsymbol{x}_\omega)
\quad
\mathrm{s.t.}\quad \boldsymbol{x}_\omega\in \mathcal{X}_\omega .
\label{eq:scenario-optimum}
\end{equation}
Here, $\boldsymbol{x}_\omega$ denotes the perfect-foresight pathway of scenario $\omega$ and $\mathcal{X}_\omega$ the scenario-specific feasible set, i.e., the set of transition pathways whose investments $\boldsymbol{d}$ and operational decisions $\boldsymbol{o}_\omega$ satisfy the constraints of scenario $\omega$. Since $\mathit{TAC}_\omega$ accounts for the operation in scenario $\omega$ only, the perfect-foresight problem effectively optimizes the investments $\boldsymbol{d}$ together with the operational decisions $\boldsymbol{o}_\omega$, while the operational decisions of all other scenarios do not enter Eq.~\eqref{eq:scenario-optimum}.

For any transition pathway $\boldsymbol{x}$, the scenario-wise relative regret is
\begin{equation}
r_\omega(\boldsymbol{x})=\frac{\mathit{TAC}_\omega(\boldsymbol{x})-\mathit{TAC}_\omega^*}{\mathit{TAC}_\omega^*}.
\label{eq:relative-regret}
\end{equation}
The relative regret $r_\omega(\boldsymbol{x})$ measures the relative cost increase of transition pathway $\boldsymbol{x}$ compared with the perfect-foresight optimum in scenario $\omega$. We define the optimal pathway as the transition pathway that minimizes the worst-case regret across all scenarios, which is known as the minimax-regret criterion \citep{Savage1951}. With auxiliary variable $R$, the minimax-regret problem is
\begin{equation}
r^*=
\min_{\boldsymbol{x},R} R
\quad
\mathrm{s.t.}\quad
\boldsymbol{x}\in \mathcal{X},\quad
r_\omega(\boldsymbol{x})\le R \quad \forall \omega\in\Omega .
\label{eq:minimax-regret}
\end{equation}
The optimal value $r^*$ is the minimum achievable worst-case regret.\footnote{Since in Eq.~\eqref{eq:minimax-regret} only the constraint of the scenario with the largest regret is active, we add the mean regret over all scenarios, weighted with a factor of $10^{-3}$, to its objective, so that the operation is cost-minimal in the remaining scenarios.} Computing the optimal pathway requires $|\Omega|+1$ optimization runs of the transition model: one perfect-foresight problem (Eq.~\eqref{eq:scenario-optimum}) per scenario, and the minimax-regret problem (Eq.~\eqref{eq:minimax-regret}), which contains all scenarios simultaneously.

\subsection{Low-regret pathway space}
\label{subsec:lowregret}

Building on the optimal pathway defined in Section~\ref{subsec:model}, we constrain the permitted worst-case regret to define the low-regret pathway space, corresponding to Step~2 in Figure~\ref{fig:method-overview}. Let $\bar r\ge r^*$ denote the regret tolerance, i.e., the largest worst-case regret allowed across all scenarios. For a given tolerance $\bar r$, the low-regret pathway space is
\begin{equation}
\mathcal{X}_{\bar r}=\{\boldsymbol{x}\in \mathcal{X}\mid r_\omega(\boldsymbol{x})\le \bar r,\;\forall\omega\in\Omega\}.
\label{eq:low-regret-decision-space}
\end{equation}
The lowest meaningful tolerance is $\bar r=r^*$, at which the space reduces to the set of minimax-regret optimal pathways. Because $\mathit{TAC}_\omega^*$ is constant, the regret tolerance is equivalent to the linear constraint $\mathit{TAC}_\omega(\boldsymbol{x})\le (1+\bar r)\,\mathit{TAC}_\omega^*$ for all $\omega\in\Omega$. Since the operational decisions of each scenario enter only that scenario's constraint, the first-stage decisions of $\mathcal{X}_{\bar r}$ form the intersection of the near-optimal design sets of the single scenarios at the same tolerance.

Constraining regret rather than absolute cost generalizes the near-optimal cost slack of MGA for stochastic optimization problems. Each scenario's cost is constrained relative to that scenario's optimum, so a single tolerance $\bar r$ is comparable across scenarios regardless of their absolute cost.

\subsection{Directional exploration of the low-regret pathway space}
\label{subsec:exploration}

This section details Step~3 in Figure~\ref{fig:method-overview} and proceeds in three parts. We first project the low-regret pathway space onto a lower-dimensional subspace of the decision variables (Section~\ref{subsubsec:projection}). We then explore the boundary of the projected space by search directions and derive an inner and an outer approximation (Section~\ref{subsubsec:inner-outer}). Finally, we define how search directions are selected and how we measure the exploration progress (Section~\ref{subsubsec:directions-termination}).

\subsubsection{Projection onto decision-relevant variables}
\label{subsubsec:projection}

The low-regret pathway space $\mathcal{X}_{\bar r}$ is not explored in its full dimension. The decision vector $\boldsymbol{x}$ contains the first-stage investment and the second-stage operational decisions of every scenario, year, and time step, whose dimension can reach millions. Individual low-regret pathways can still be computed in this space, but its geometry cannot be characterized. We therefore analyze a lower-dimensional subspace of decision-relevant variables, obtained by projecting the full pathway space onto that subspace.

The exploration variables, i.e., the pathway variables we explore, are selected by a linear projection $\hat{\boldsymbol{x}}=\mathbf{S}\boldsymbol{x}\in\mathbb{R}^{m}$ of the full pathway space, where $\mathbf{S}$ denotes the projection matrix and $m$ the number of exploration variables. In our case study, the projection aggregates the first-stage investment decisions into technology families per investment year (Section~\ref{subsec:exploration-variables}). The projected low-regret pathway space is the image of $\mathcal{X}_{\bar r}$ under $\mathbf{S}$,
\begin{equation}
\mathcal{P}_{\bar r}=\{\hat{\boldsymbol{x}}\in\mathbb{R}^{m}\mid \exists \boldsymbol{x}\in \mathcal{X}_{\bar r}\;\mathrm{with}\;\hat{\boldsymbol{x}}=\mathbf{S}\boldsymbol{x}\}.
\label{eq:projected-space}
\end{equation}
Because $\mathcal{X}$ is a convex polyhedron and both the regret constraints and $\mathbf{S}$ are linear, $\mathcal{X}_{\bar r}$ and $\mathcal{P}_{\bar r}$ are convex polyhedra. Thus, convex combinations of low-regret pathways, i.e., of their decision vectors $\boldsymbol{x}$ or projections $\hat{\boldsymbol{x}}$, are again low-regret pathways.

\subsubsection{Inner and outer approximation}
\label{subsubsec:inner-outer}

Following common MGA approaches, we explore the projected pathway space $\mathcal{P}_{\bar r}$ by optimization runs along search directions. For a given search direction $\boldsymbol{n}\in\mathbb{R}^{m}$, we explore the boundary of $\mathcal{P}_{\bar r}$ by solving
\begin{equation}
\beta(\boldsymbol{n}) =
\max_{\boldsymbol{x},\hat{\boldsymbol{x}}} \boldsymbol{n}^\top \hat{\boldsymbol{x}}
\quad
\mathrm{s.t.}\quad
\boldsymbol{x}\in \mathcal{X}_{\bar r},\quad
\hat{\boldsymbol{x}}=\mathbf{S}\boldsymbol{x} .
\label{eq:directional-lp}
\end{equation}
Solving Eq.~\eqref{eq:directional-lp} returns a boundary point of $\mathcal{P}_{\bar r}$ together with a supporting halfspace $\boldsymbol{n}^\top\hat{\boldsymbol{x}}\le\beta(\boldsymbol{n})$ with support value $\beta(\boldsymbol{n})$. We collect the solved search directions as the searched directions $\boldsymbol{n}_h$, $h=1,\dots,H$, where $H$ is the number of solutions of Eq.~\eqref{eq:directional-lp} so far.\footnote{Since the objective of Eq.~\eqref{eq:directional-lp} determines only the boundary point, we subtract the mean regret over all scenarios, weighted with a factor of $10^{-3}$, from this objective, so that the operation is cost-minimal in the scenarios that do not reach the regret tolerance.}

\begin{figure}[p]
\centering
\includegraphics[width=\textwidth,height=0.80\textheight,keepaspectratio]{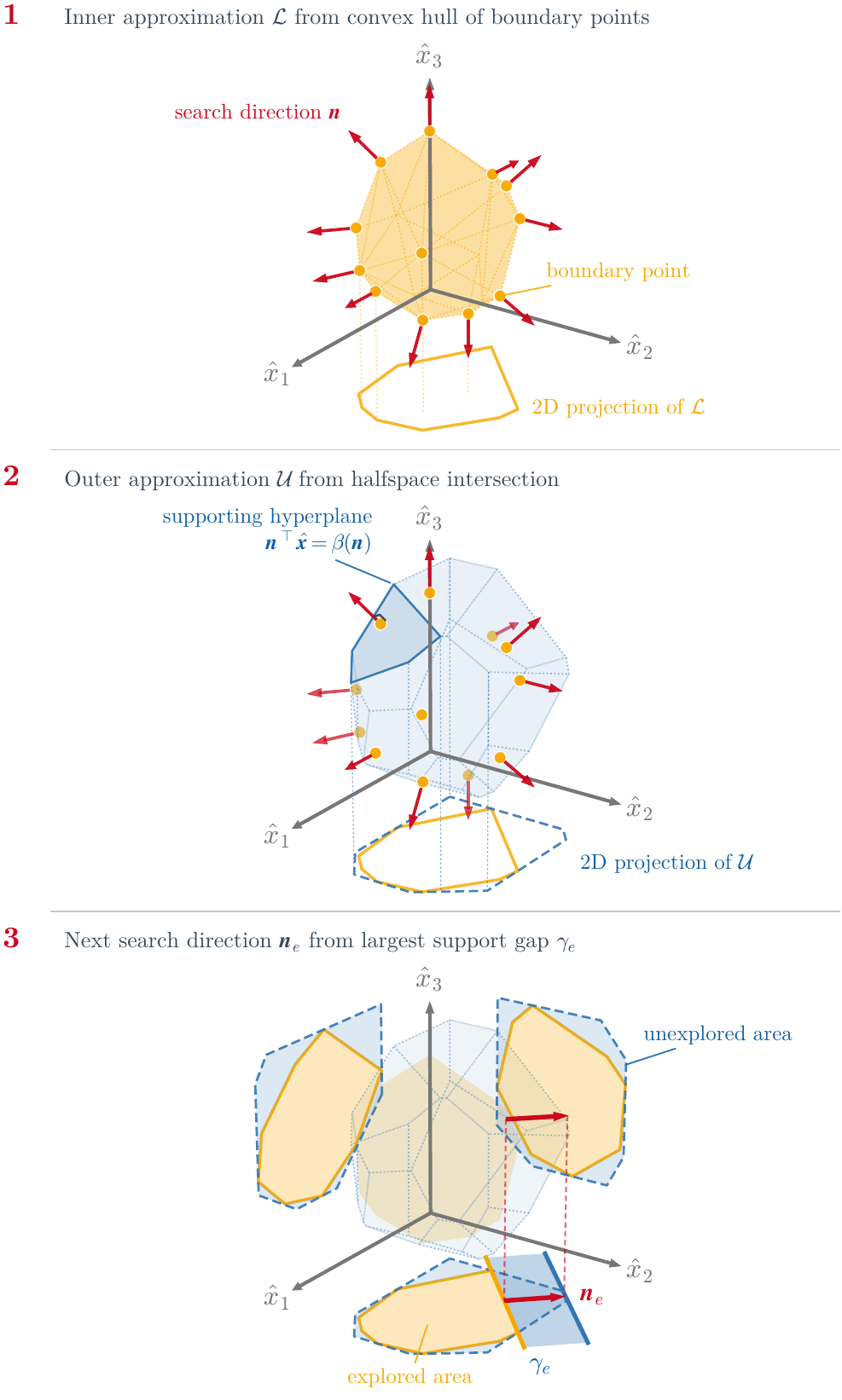}
\caption{Schematic overview of the exploration of the projected low-regret pathway space, illustrated for three exploration variables. \textbf{(1)}~Each solution of Eq.~\eqref{eq:directional-lp} along a search direction $\boldsymbol{n}$ returns a boundary point whose convex hull forms the inner approximation $\mathcal{L}$, and \textbf{(2)}~a supporting halfspace $\boldsymbol{n}^\top\hat{\boldsymbol{x}}\le\beta(\boldsymbol{n})$ whose intersection forms the outer approximation $\mathcal{U}$. \textbf{(3)}~In each two-dimensional projection, the edges $e$ with the largest support gap $\gamma_e$ between the inner and outer approximation are searched next. The search direction $\boldsymbol{n}_e$ (arrow) is the outward edge normal.}
\label{fig:approx-sampling}
\end{figure}

Because $\mathcal{P}_{\bar r}$ is convex, the convex hull of the boundary points discovered so far forms an inner approximation $\mathcal{L}$ \citep{Pedersen2021}. The intersection of the discovered supporting halfspaces yields an outer approximation $\mathcal{U}=\{\hat{\boldsymbol{x}} : \boldsymbol{n}_h^\top\hat{\boldsymbol{x}}\le\beta(\boldsymbol{n}_h)\ \text{for all}\ h\}$ \citep{Schricker2025b}, so that
\begin{equation}
\mathcal{L}\subseteq \mathcal{P}_{\bar r}\subseteq \mathcal{U}.
\label{eq:inner-outer}
\end{equation}
Both approximations are monotone in the number $H$ of solutions of Eq.~\eqref{eq:directional-lp}, since an additional boundary point can only enlarge $\mathcal{L}$ and an additional halfspace can only shrink $\mathcal{U}$. Figure~\ref{fig:approx-sampling} schematically illustrates the inner and outer approximation.

\subsubsection{Search direction selection and exploration progress}
\label{subsubsec:directions-termination}

The exploration serves the transition flexibility assessment of Section~\ref{subsec:decision-support}, which reads the range of each exploration variable from its coordinate-wise bounds and the combinability of the pathway decisions from the two-dimensional coordinate projections of the pathway space. Our exploration method is therefore designed to provide exact results for precisely these bounds and projections rather than to cover the full $m$-dimensional volume. Thus, full-dimensional coverage is neither required nor intended.

Designing an exploration method requires two methodological choices: a search strategy \citep[for a review, see][]{Lau2024b} and a criterion for terminating the search. Existing termination criteria fix the number of searches in advance, focus on convergence of the inner approximation \citep{Pedersen2021}, or measure the explored share of the near-optimal space, e.g., as the volume ratio of the inner and outer approximation \citep{Schricker2025b} or as the largest gap between the inner and outer approximation \citep{Turan2026}. Neither choice is trivial in high dimensions, since, e.g., exact volumes of high-dimensional polyhedra are intractable to compute \citep{Barber1996,Turan2026}. We therefore design a tailored method that advances the state of the art in two ways. First, our method bounds coverage in the two-dimensional coordinate projections used by the volume estimation by shadow addition (VESA) of \citet{Lau2024b}. These projections remain tractable in high dimensions. To this end, we transfer the degree of exploration of \citet{Schricker2025b}, i.e., the ratio of the volumes of inner and outer approximation, from the intractable full-dimensional volumes to the projection areas. Second, we direct new searches toward large gaps between the inner and outer approximation as in ORACLE \citep{Turan2026}. Since the bilevel problem of ORACLE dominates the computational effort, G-ORACLE replaces it by an enumeration of the vertices of the outer approximation \citep{Hajikazemi2026}. The outer approximation starts from the box of the coordinate-wise bounds, which already has $2^m$ vertices. We instead evaluate the gaps at the edges of the $m(m-1)/2$ two-dimensional coordinate projections with one small linear program per edge. Each projection has at most as many edges as explored boundary points.

We initialize the searched directions $\boldsymbol{n}_h$ with the positive and negative coordinate unit directions, i.e., $H=2m$ solutions of Eq.~\eqref{eq:directional-lp}, which yield exact coordinate-wise bounds. Each exploration variable is then rescaled by its full range to $[0,1]$. From the explored boundary points, we construct the two-dimensional coordinate projections of the inner approximation $\mathcal{L}$, one for every pair of exploration variables. After the rescaling, each of these projections lies in the unit square. Every subsequent search direction is derived from these projections of $\mathcal{L}$.

New candidate search directions are the outward edge normals of these projections (Figure~\ref{fig:approx-sampling}, panel~3). An edge normal is a two-dimensional vector, whereas a search direction for Eq.~\eqref{eq:directional-lp} is a vector $\boldsymbol{n}\in\mathbb{R}^{m}$ with one entry per exploration variable. For an edge~$e$ of the projection onto the variables $i$ and $j$, the search direction $\boldsymbol{n}_e$ therefore takes the two entries of the edge normal for the variables $i$ and $j$ and zero for all other exploration variables.

For each candidate search direction, we compute the gap between inner and outer approximation along that search direction. The inner support value $\beta_{\mathcal{L}}(\boldsymbol{n}_e)$, i.e., the support value of $\mathcal{L}$ along $\boldsymbol{n}_e$, is the largest value of $\boldsymbol{n}_e^\top\hat{\boldsymbol{x}}$ over the explored boundary points, i.e., the position of the edge itself. The outer support value $\beta_{\mathcal{U}}(\boldsymbol{n}_e)$ follows from the auxiliary linear program
\begin{equation}
\beta_{\mathcal{U}}(\boldsymbol{n}_e) =
\max_{\hat{\boldsymbol{x}}} \boldsymbol{n}_e^\top \hat{\boldsymbol{x}}
\quad
\mathrm{s.t.}\quad
\hat{\boldsymbol{x}}\in \mathcal{U} .
\label{eq:outer-support-lp}
\end{equation}
In contrast to Eq.~\eqref{eq:directional-lp}, Eq.~\eqref{eq:outer-support-lp} does not contain the energy system model, since $\mathcal{U}$ consists of the halfspaces of the searched directions $\boldsymbol{n}_h$ alone. The difference of the outer and the inner support value is the support gap
\begin{equation}
\gamma_e=\beta_{\mathcal{U}}(\boldsymbol{n}_e)-\beta_{\mathcal{L}}(\boldsymbol{n}_e)\ge 0.
\label{eq:support-gap}
\end{equation}

With the support gap, we decide which edges are selected for further exploration. All further optimization runs are allocated in batches to the edges in the order of decreasing support gap and solved in parallel. Each solved run adds its direction $\boldsymbol{n}_e$ to the searched directions $\boldsymbol{n}_h$ and its halfspace to $\mathcal{U}$.

We track the exploration progress by comparing the inner and outer approximation in the projections. The outer support values of the edges bound each projection by a polygon. The area of this polygon is the projection area of $\mathcal{U}$ in that coordinate pair. Summing the projection areas of $\mathcal{L}$ and $\mathcal{U}$ over all coordinate pairs yields the inner and outer projection areas $A_L$ and $A_U$. Their ratio is the degree of exploration
\begin{equation}
\mathit{DoE}_{\mathrm{VESA}}=\frac{A_L}{A_U}.
\label{eq:doe-vesa}
\end{equation}
The subscript distinguishes this ratio of summed two-dimensional projection areas as in VESA \citep{Lau2024b} from the volume-based degree of exploration of \citet{Schricker2025b}. The exploration terminates once the degree of exploration reaches a prespecified tolerance or a compute budget is exhausted. In both cases, the final degree of exploration quantifies the coverage achieved.

\subsection{Quantitative assessment of transition flexibility}
\label{subsec:decision-support}

\begin{figure}[!htbp]
\centering
\includegraphics[width=0.66\textwidth]{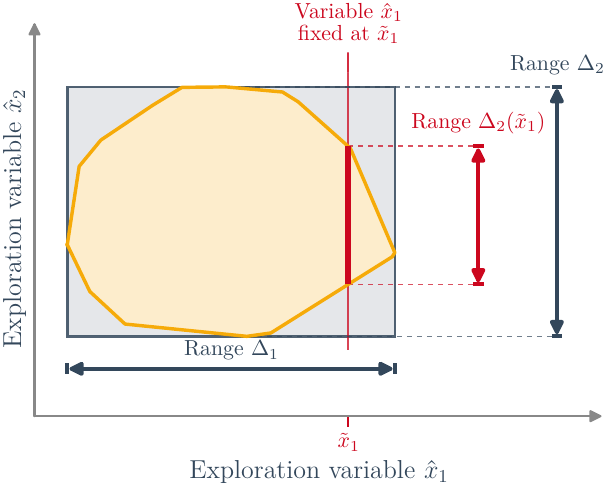}
\caption{Schematic illustration of the transition flexibility indicators in the two-dimensional projection of the pathway space onto the exploration variables $\hat{x}_1$ and $\hat{x}_2$. The ranges $\Delta_1$ and $\Delta_2$ of the two exploration variables span the gray box. The first indicator is the range $\Delta_1$. Fixing variable 1 at the single value $\tilde{x}_1$ (red line) leaves a narrower range $\Delta_2(\tilde{x}_1)$ of variable 2. The second indicator, the conditional range $R_1$, is the mean ratio of the ranges after fixing $\Delta_j(\tilde{x}_1)$ to the full ranges $\Delta_j$ over all other exploration variables $j\neq 1$.}
\label{fig:conditional-flexibility}
\end{figure}

Finally, we quantify the transition flexibility of the explored pathway space, i.e., the freedom of decision makers to choose among low-regret pathways (Step~4 in Figure~\ref{fig:method-overview}). In the work of \citet{vanGreevenbroek2025b}, design flexibility is measured as the width of the near-optimal space that remains once a single investment is fixed. We transfer this idea from a single investment year to transition pathways. Let $\mathcal{E}\subseteq \mathcal{P}_{\bar r}$ denote the explored pathway set, $\hat{x}_i$ the $i$-th exploration variable, i.e., the $i$-th component of $\hat{\boldsymbol{x}}$, and $\tilde{x}_i$ a fixed value of it. We define two indicators of transition flexibility, illustrated in Figure~\ref{fig:conditional-flexibility}.

The first indicator is the range $\Delta_i$ of exploration variable $i$, i.e., the range of values that this variable takes across the low-regret pathway space. It is defined as $\Delta_i=\hat{x}_i^{\max}-\hat{x}_i^{\min}$ with $\hat{x}_i^{\max}=\max_{\hat{\boldsymbol{x}}\in\mathcal{E}} \hat{x}_i$ and $\hat{x}_i^{\min}$ defined analogously. To compare pathway spaces at different regret tolerances, we use the mean range $\bar{\Delta}$, i.e., the ranges $\Delta_i(\bar r)$ at tolerance $\bar r$ averaged over all exploration variables after normalizing each to its range $\Delta_i(\bar r_{\max})$ at the highest explored tolerance $\bar r_{\max}$,
\begin{equation}
\bar{\Delta}(\bar r)=\frac{1}{|\mathcal{I}|}\sum_{i\in\mathcal{I}}\frac{\Delta_i(\bar r)}{\Delta_i(\bar r_{\max})},
\label{eq:mean-range}
\end{equation}
where $\mathcal{I}$ contains all exploration variables with a nonzero range $\Delta_i(\bar r_{\max})$.

The second indicator is the conditional range $R_i$ of exploration variable $i$, i.e., the mean share of the ranges $\Delta_j$ that the other exploration variables $j\neq i$ keep once variable $i$ is fixed at the value $\tilde{x}_i$. The conditional range $R_i$ follows the conditional-width indicator of \citet{vanGreevenbroek2025b}. Fixing variable $i$ leaves every other variable $j$ a narrower range $\Delta_j(\tilde{x}_i)$, which we read from the two-dimensional projection onto the plane of the variable pair $(i,j)$ (Figure~\ref{fig:conditional-flexibility}). The conditional range $R_i$ is defined as
\begin{equation}
R_i(\tilde{x}_i)
=\frac{1}{|\mathcal{J}_i|}\sum_{j\in\mathcal{J}_i}\frac{\Delta_j(\tilde{x}_i)}{\Delta_j},
\label{eq:conditional-range}
\end{equation}
where $\mathcal{J}_i$ contains all exploration variables $j\neq i$ with a nonzero range $\Delta_j$. We refer to this share as the combinability of the pathway decision. A conditional range $R_i$ of one means that the fixed pathway decision leaves the ranges $\Delta_j$ of all other exploration variables unchanged, whereas a conditional range close to zero means that it determines almost all other pathway decisions as well.

If the exploration variables are installed capacities, as in our case study, the upper part of a range $\Delta_i$ can be \emph{surplus capacity}, i.e., capacity installed in addition to, not instead of, other capacity. We identify it as the capacity above the largest value that the capacity-minimizing searches reach, i.e., the searches whose direction coefficients are all non-positive.
\clearpage
\section{Case study}
\label{sec:case-study}

We model the transition pathway of an industrial energy system from 2025 to 2045 in five year steps. This section introduces the energy system superstructure (Section~\ref{subsec:structure}), the scenarios for the considered uncertainties (Section~\ref{subsec:scenarios}), and the exploration variables (Section~\ref{subsec:exploration-variables}).

\subsection{Energy system structure}
\label{subsec:structure}

The modeled industrial multi-energy system supplies steam, heating, cooling, and electricity to the site. Excess heat from processes can either be reintegrated or released unused to the environment. In the existing system, heat and steam are supplied by natural-gas-fired boilers and combined heat and power (CHP) units, cooling by absorption and compression chillers, and electricity and natural gas are purchased at the spot market.

\begin{landscape}
\begin{figure}[!ht]
\centering
\includegraphics[width=0.93\linewidth,keepaspectratio]{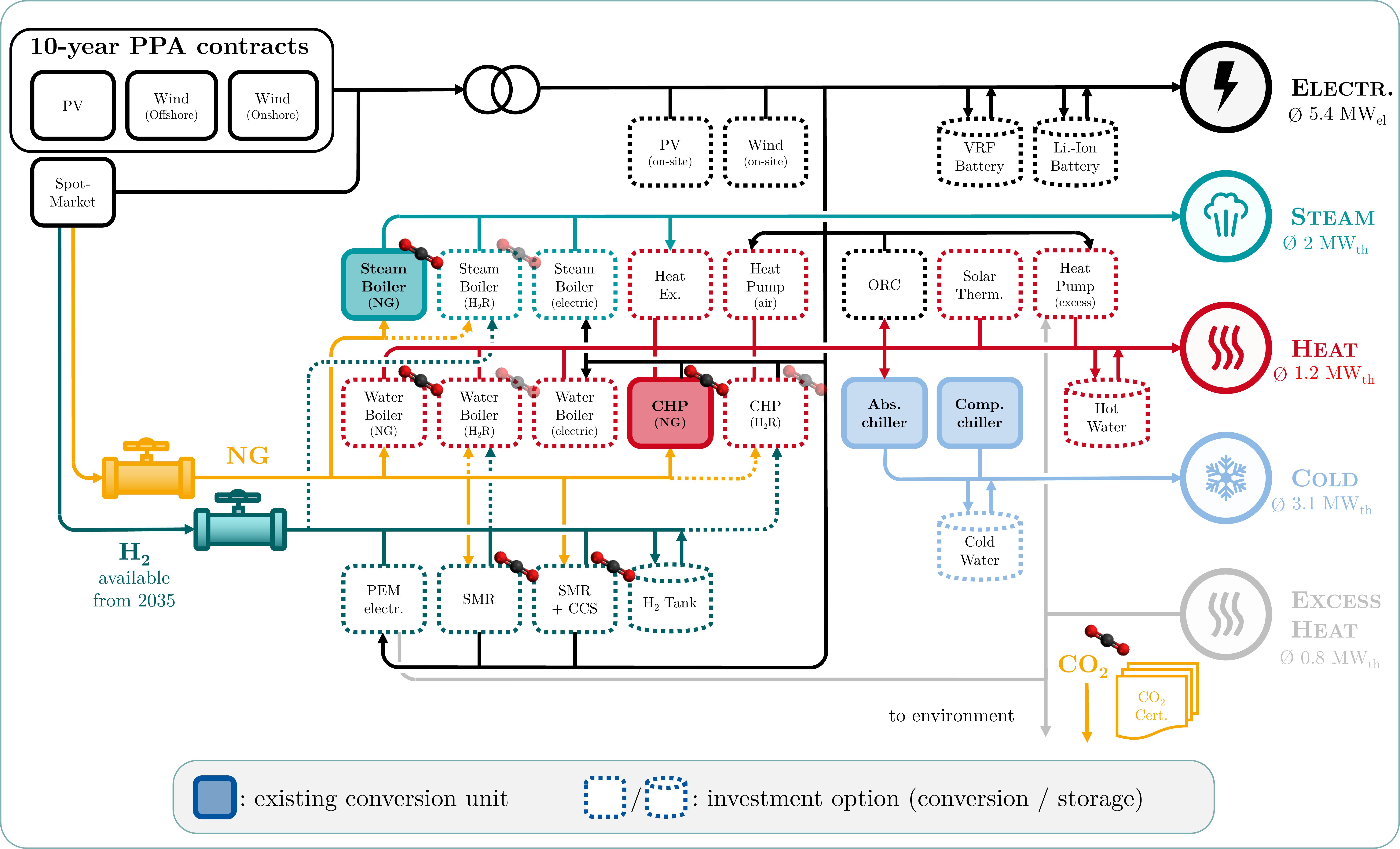}
\caption{Superstructure of the industrial energy system. The system supplies electricity, steam, heating, and cooling. Connecting flows denote the energy carriers exchanged between components and demands. CCS: carbon capture and storage, CHP: combined heat and power, H$_2$: hydrogen, H$_2$R: hydrogen-ready, Li-ion battery: lithium-ion battery, NG: natural gas, ORC: organic Rankine cycle, PEM electr.: proton exchange membrane electrolyzer, PPA: power purchase agreement, PV: photovoltaic, SMR: steam methane reformer, VRF battery: vanadium redox flow battery.}
\label{fig:superstructure}
\end{figure}
\end{landscape}

Figure~\ref{fig:superstructure} shows the superstructure of the energy system with the existing components and all investment options. The future system considers a broad set of investment options: renewable electricity and heat supply (on-site or contracted via power purchase agreements (PPAs)), sector-coupling technologies, hydrogen technologies, and energy storage. We assume that the hydrogen grid is available from 2035 onwards. In 2025, only existing components and PPAs are available, while all further investment options are available from 2030 onwards. Hydrogen-ready technologies can be operated flexibly with either natural gas or hydrogen. All relevant techno-economic parameters and assumptions of the case study are reported in Section~S1 of the Supplementary Information.

\newgeometry{margin=2cm}
\restoregeometry

\subsection{Scenario and time series data}
\label{subsec:scenarios}

\begin{figure}[htbp]
\centering
\includegraphics[width=\textwidth]{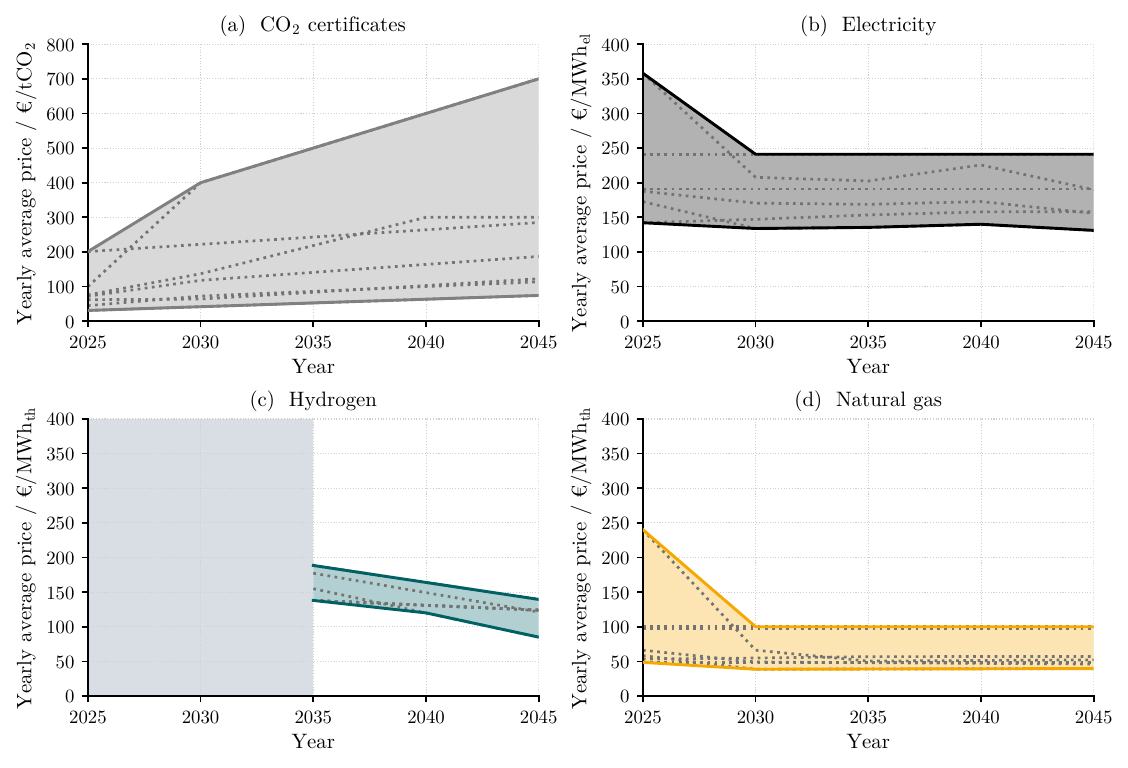}
\caption{Price trajectories from 2025 to 2045 for the four uncertain input parameters. Gray dotted lines show the seven scenarios. The shaded band spans the scenario range (minimum to maximum across scenarios). \textbf{(a)}~CO$_2$ certificate prices. \textbf{(b)}~ Spot market electricity purchase prices. \textbf{(c)}~Hydrogen prices, shown from the availability of the hydrogen grid connection in 2035, with the gray region marking the years without hydrogen supply. \textbf{(d)}~Natural gas prices.}
\label{fig:scenarios}
\end{figure}

We model future energy carrier and CO$_2$ certificate price uncertainty by a set of seven scenarios. Each scenario consists of one long-term price trajectory per uncertain price (spot market electricity, natural gas, hydrogen, and CO$_2$ certificates) from 2025 to 2045. Figure~\ref{fig:scenarios} shows the resulting price trajectories of the seven scenarios, and Table~\ref{tab:scenario-prices} reports the CO$_2$ certificate price and the annual average purchase prices of each scenario in the first and last year of its trajectory.

Each scenario uses hourly resolved demand, price, and renewable availability time series for every investment year. The intra-annual profiles are identical across scenarios. The price time series are scaled by each scenario's long-term trajectory. The electricity price profile is year-specific: its intra-annual profile becomes more volatile in later investment years to reflect the rising share of renewable electricity generation \citep{Prognos2024}. To reduce computational complexity, the hourly time series of each investment year and scenario are aggregated into four typical days. With this aggregation, each two-stage stochastic linear program of the transition model comprises $550{,}751$ variables and $464{,}011$ constraints. Details on the price trajectories, the time series data, and the time series aggregation are provided in Sections~S1 and~S2.1 of the Supplementary Information.

\begin{table}[htbp]
\centering
\caption{The seven price scenarios. Reported are the CO$_2$ certificate price and the annual average electricity, hydrogen, and natural gas purchase prices. The full price trajectories are reported in Section~S1.4.3 of the Supplementary Information.}
\label{tab:scenario-prices}
\footnotesize
\setlength{\tabcolsep}{3pt}
\begin{tabular}{c rr rr rr rr c}
\toprule
& \multicolumn{2}{c}{CO$_2$ certificate} & \multicolumn{2}{c}{Electricity} & \multicolumn{2}{c}{Hydrogen$^{\mathrm{*}}$} & \multicolumn{2}{c}{Natural gas} & \\
& \multicolumn{2}{c}{/ EUR/tCO$_2$} & \multicolumn{2}{c}{/ EUR/MWh} & \multicolumn{2}{c}{/ EUR/MWh} & \multicolumn{2}{c}{/ EUR/MWh} & \\
\cmidrule(lr){2-3}\cmidrule(lr){4-5}\cmidrule(lr){6-7}\cmidrule(lr){8-9}
Scenario & 2025 & 2045 & 2025 & 2045 & 2035 & 2045 & 2025 & 2045 & Reference \\
\midrule
1 & 76 & 300 & 358 & 189.8 & 138 & 124 & 240.1 & 52.2 & $^{\mathrm{a,b,c}}$ \\
2 & 72 & 187 & 187.6 & 155.7 & 138 & 124 & 66.1 & 46.2 & $^{\mathrm{a,b,c}}$ \\
3 & 62 & 122.5 & 172.5 & 131.1 & 138 & 124 & 58.3 & 39.6 & $^{\mathrm{a,b,c}}$ \\
4 & 31 & 74.4 & 142.2 & 158.1 & 138 & 124 & 53.3 & 57 & $^{\mathrm{c,d}}$ \\
5 & 45 & 113.8 & 191 & 191 & 188.5 & 139.5 & 48.6 & 48.6 & $^{\mathrm{e}}$ \\
6 & 201 & 285 & 241 & 241 & 155 & 85 & 100 & 100 & $^{\mathrm{e}}$ \\
7 & 100 & 700 & 191 & 191 & 177.3 & 121.3 & 97.2 & 97.2 & $^{\mathrm{e}}$ \\
\bottomrule
\end{tabular}
\par\smallskip
{\footnotesize\raggedright $^{\mathrm{*}}$~The hydrogen grid is assumed to be available from 2035 onwards, so hydrogen prices start in 2035 instead of 2025.\\
$^{\mathrm{a}}$~\citet{Prognos2024}, $^{\mathrm{b}}$~\citet{Bundesnetzagentur2024}, $^{\mathrm{c}}$~\citet{Agora2021}, $^{\mathrm{d}}$~\citet{TUM_thesis2020}, $^{\mathrm{e}}$~\citet{Gorbach2022}.\par}
\end{table}

\subsection{Exploration variables}
\label{subsec:exploration-variables}

We project each transition pathway onto the site's first-stage investment decisions aggregated into 14 technology families, summarized in Table~\ref{tab:exploration-variables}. The exploration variable of a technology family in an investment year is its summed installed capacity in that investment year. Each technology family aggregates the components with similar function. ``Thermal integration'' aggregates the technologies that convert thermal streams into other energy carriers (heat from steam via heat exchangers, heat from excess heat via heat pumps, electricity via organic Rankine cycles, and cooling via absorption chillers).

Resolving each technology family's installed capacity in the five investment years gives a projected pathway space of $m=14\times 5=70$ dimensions. Because the 2025 capacities of all technology families except the PPAs and the grid connections are fixed (Section~\ref{subsec:structure}) and the hydrogen grid connection is fixed at zero until 2035, twelve of the $m=70$ dimensions are constant. Thus, we explore $58$ dimensions of the low-regret pathway space.

\begin{table}[htbp]
\centering
\caption{Technology families whose installed capacities define the low-regret pathway space. Each technology family aggregates the installed capacity of the listed technologies, resolved in every investment year. The capacities of renewable supply, conversion, hydrogen production, and grid connections are measured in MW, those of storage in MWh. CCS: carbon capture and storage, CHP: combined heat and power, H$_2$R: hydrogen-ready, ORC: organic Rankine cycle, PEM: proton exchange membrane, PPA: power purchase agreement, PV: photovoltaic, SMR: steam methane reformer.}
\label{tab:exploration-variables}
\small
\begin{tabular}{c >{\raggedright\arraybackslash}p{3.4cm} >{\raggedright\arraybackslash}p{7.2cm}}
\toprule
No. & Technology family & Aggregated technologies \\
\midrule
\multicolumn{3}{l}{\emph{Renewable supply}} \\
1 & PPA & Photovoltaic and wind (on-/offshore) PPAs \\
2 & Renewables on-site & On-site photovoltaics, wind, and solar-thermal collectors \\
\addlinespace
\multicolumn{3}{l}{\emph{Conversion}} \\
3 & NG fired & Natural-gas-fired CHP units and boilers \\
4 & H$_2$ ready & Hydrogen-ready CHP units and boilers \\
5 & Electrified heating and cooling & Electric boilers, air-source heat pumps, and compression chillers \\
6 & Thermal integration & Heat exchangers, excess-heat heat pumps, ORCs, and absorption chillers \\
\addlinespace
\multicolumn{3}{l}{\emph{Hydrogen production}} \\
7 & Electrolysis & PEM electrolyzers \\
8 & SMR & Steam methane reformers, with and without CCS \\
\addlinespace
\multicolumn{3}{l}{\emph{Storage}} \\
9 & Battery storage & Lithium-ion and vanadium-redox-flow batteries \\
10 & Thermal storage & Hot-water and cold-water tanks \\
11 & Hydrogen storage & Hydrogen tanks \\
\addlinespace
\multicolumn{3}{l}{\emph{Grid connections}} \\
12 & Grid electricity & Electricity grid connection \\
13 & Grid natural gas & Natural gas grid connection \\
14 & Grid hydrogen & Hydrogen grid connection \\
\bottomrule
\end{tabular}
\end{table}

\clearpage

\section{Results}
\label{sec:results}

This section presents the results of applying the method of Section~\ref{sec:methods} to the industrial energy system of Section~\ref{sec:case-study}. We explore low-regret pathway spaces at regret tolerances $\bar r$ of 2.5\%, 3.5\%, 5\%, and 10\%. Section~\ref{sec:results-casestudy} introduces the minimax-regret optimal pathway and compares it with the perfect-foresight pathways. Section~\ref{sec:results-space} explores the low-regret pathway spaces and analyzes the resulting capacity ranges (\ref{rq:pathway-space}). Section~\ref{sec:results-cost-of-flexibility} analyzes the trade-off between regret tolerance and transition flexibility (\ref{rq:cost-of-flexibility}), and Section~\ref{sec:results-conditional} quantifies the combinability of the pathway decisions (\ref{rq:residual-flexibility}).

\subsection{Minimax-regret pathway versus perfect-foresight pathways}
\label{sec:results-casestudy}

\begin{figure}[htbp]
\centering
\includegraphics[width=\textwidth]{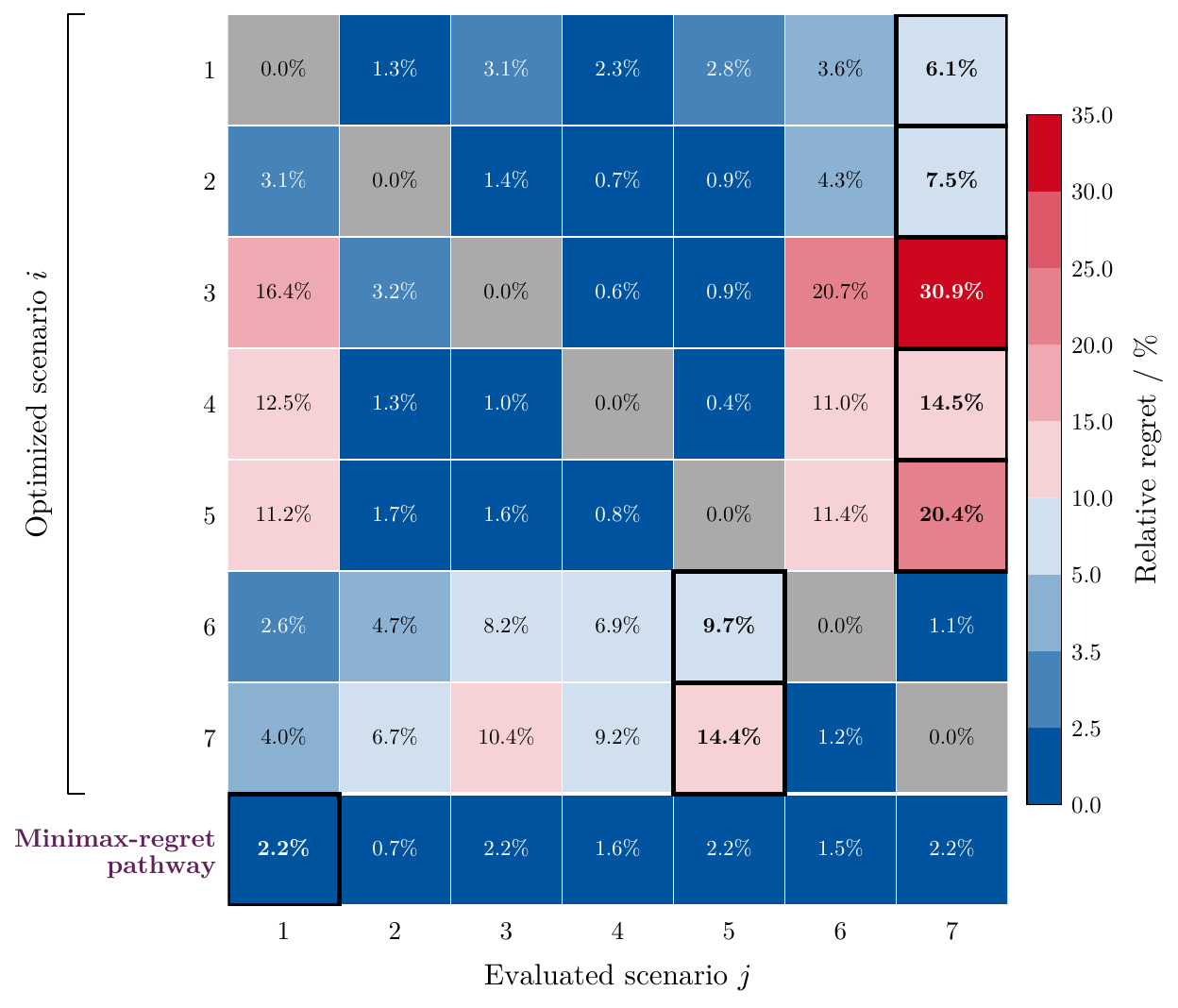}
\caption{Relative regret (\%) of the seven perfect-foresight pathways, each optimized for one scenario $i$ of Table~\ref{tab:scenario-prices}, and of the stochastic minimax-regret pathway. All pathways are evaluated under all seven scenarios $j$. The gray diagonal entries are zero by construction. Bold borders mark each pathway's worst-case
scenario.}
\label{fig:regret-ranges}
\end{figure}

Figure~\ref{fig:regret-ranges} compares the stochastic minimax-regret pathway with the seven perfect-foresight pathways. Each transition pathway is evaluated against all seven price scenarios. The perfect-foresight pathways lead to worst-case regrets of 6.1--30.9\%, whereas the minimax-regret pathway achieves the lowest worst-case regret of 2.2\%.

\begin{figure}[htbp]
\centering
\includegraphics[width=\textwidth]{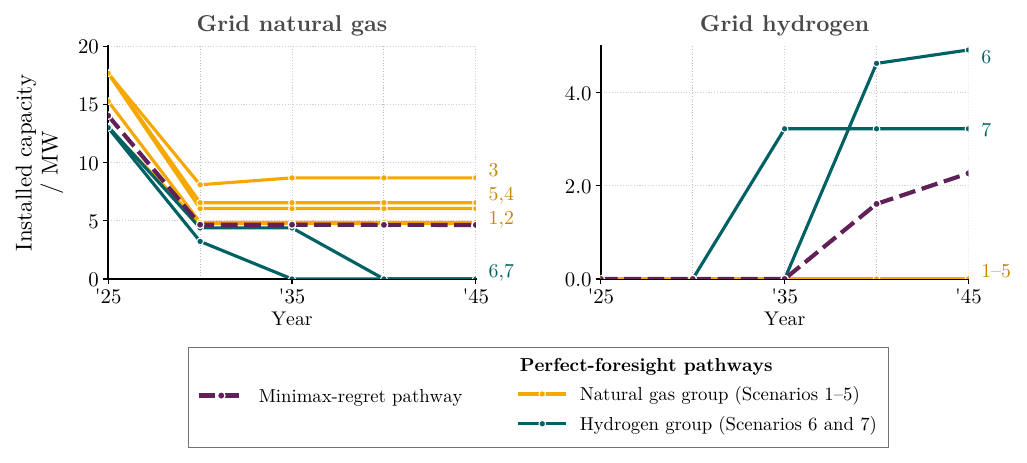}
\caption{Natural gas and hydrogen grid connection capacities of the seven perfect-foresight pathways and of the minimax-regret pathway (dashed purple line). Each perfect-foresight pathway is colored by its energy carrier group and annotated with its scenario identifier at the final year.}
\label{fig:grid-commitment}
\end{figure}

All installed capacities of the seven perfect-foresight pathways are reported in Section~S2.2 of the Supplementary Information. Each perfect-foresight pathway builds for the energy carrier that is cheaper in its own scenario, either natural gas including the CO$_2$ certificate costs of its combustion or hydrogen. This choice divides the seven scenarios, and the perfect-foresight pathways optimized for them, into a natural gas group (Scenarios 1--5) and a hydrogen group (Scenarios 6 and~7), as Figure~\ref{fig:grid-commitment} shows. In the natural gas group, natural gas including certificate costs stays cheaper than hydrogen from 2035 onwards, when hydrogen becomes available, so the pathways install natural-gas-fired components and no hydrogen grid connection. In the hydrogen group, the certificate price makes hydrogen the cheaper energy carrier by 2045, so the pathways install hydrogen-ready components and phase out the natural gas grid connection. No perfect-foresight pathway can therefore switch energy carriers. The regret matrix reflects this split. A pathway's regret averages 3.2\% under the scenarios of its own group and 10.4\% under those of the other group. The extreme case is the pathway optimized for Scenario~3, the scenario with the cheapest natural gas including certificate costs. This pathway keeps the largest natural-gas-fired capacity and reaches a regret of 30.9\% under Scenario~7, the scenario with the most expensive natural gas including certificate costs. Only the minimax-regret pathway installs both grid connections and can switch energy carriers.

The minimax-regret pathway (purple line in Figure~\ref{fig:range-clean}) contracts the full PPA capacity in all years and installs most of the maximum installable on-site renewable capacity in 2030. The natural-gas-fired capacity falls to a small margin from 2035 onwards, because the existing units retire and are not replaced. Instead, the pathway installs hydrogen-ready boilers in 2030, which are fired with natural gas until the hydrogen grid connection is contracted in 2040. The electrified heating and cooling capacity stays almost constant, because the retiring compression chillers are replaced in 2030 by a smaller chiller, an electric steam boiler, and an air-source heat pump. Thermal and battery storage are added in 2030, and a small hydrogen storage in 2045. On-site hydrogen production is not installed. The natural gas grid connection shrinks in 2030 but is kept until 2045. The following sections examine how much transition flexibility is available around this pathway, both in the choice of technologies and in the timing of their installation.

\subsection{Exploring low-regret pathway spaces}
\label{sec:results-space}

The low-regret pathway spaces are explored around the minimax-regret pathway at the four regret tolerances. In this section, we first verify the convergence of the pathway space explorations (Section~\ref{subsubsec:results-convergence}). We then report the individual capacity ranges (Section~\ref{subsubsec:results-ranges}) and separate the surplus capacity within them (Section~\ref{subsubsec:results-surplus}). Section~\ref{subsubsec:results-summary} summarizes the answer to \ref{rq:pathway-space}.

\subsubsection{Convergence of the exploration}
\label{subsubsec:results-convergence}

\begin{figure}[htbp]
\centering
\includegraphics[width=\textwidth]{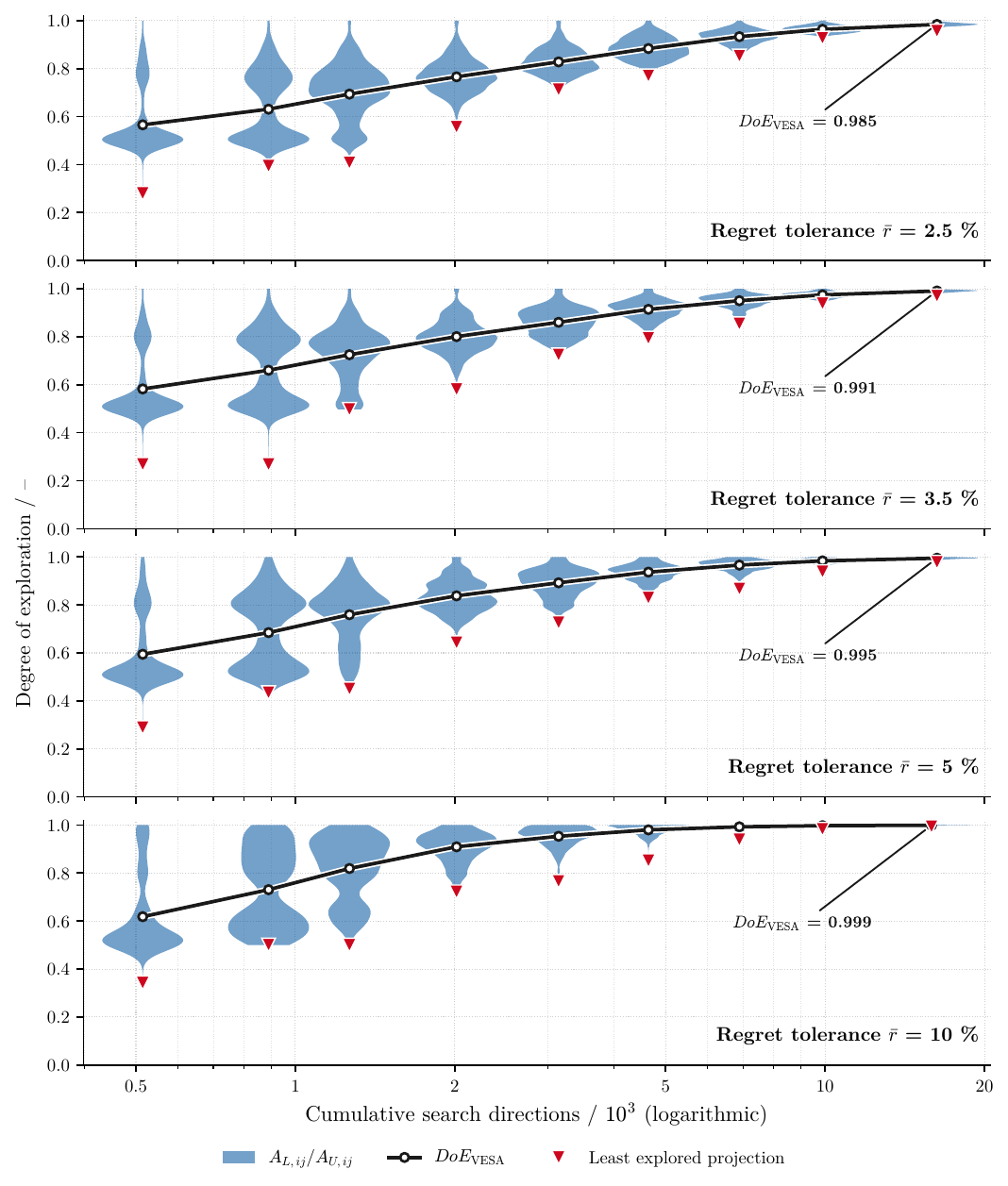}
\caption{Convergence of the exploration against the
cumulative number of explored search directions, at four regret tolerances $\bar r$. The compute budget is exhausted after about 16{,}000 directional optimizations per pathway space (Section~S2.1 of the Supplementary Information). Each violin is the distribution of the ratios
$A_{L,ij}/A_{U,ij}$ over the 1653 coordinate projections of the 58 explored dimensions,
where $A_{L,ij}$ and $A_{U,ij}$ are the areas of the inner approximation $\mathcal{L}$
and the outer approximation $\mathcal{U}$ projected onto the coordinate pair $(i,j)$.
The reported degree of exploration $\mathit{DoE}_{\mathrm{VESA}}$ follows
Eq.~\eqref{eq:doe-vesa}. Red triangles mark the least explored projection.}
\label{fig:doe-convergence}
\end{figure}

Figure~\ref{fig:doe-convergence} shows the convergence of the exploration against the cumulative number of explored search directions at the four regret tolerances. The exploration of the low-regret pathway spaces around the minimax-regret pathway converges at every regret tolerance. When the compute budget is exhausted, the degree of exploration $\mathit{DoE}_{\mathrm{VESA}}$ of Eq.~\eqref{eq:doe-vesa} reaches 0.985 to 0.9994. In this case study, higher regret tolerances converge earlier and leave a smaller remaining error.

\subsubsection{Capacity ranges}
\label{subsubsec:results-ranges}

\begin{figure}[htbp]
\centering
\includegraphics[width=\textwidth]{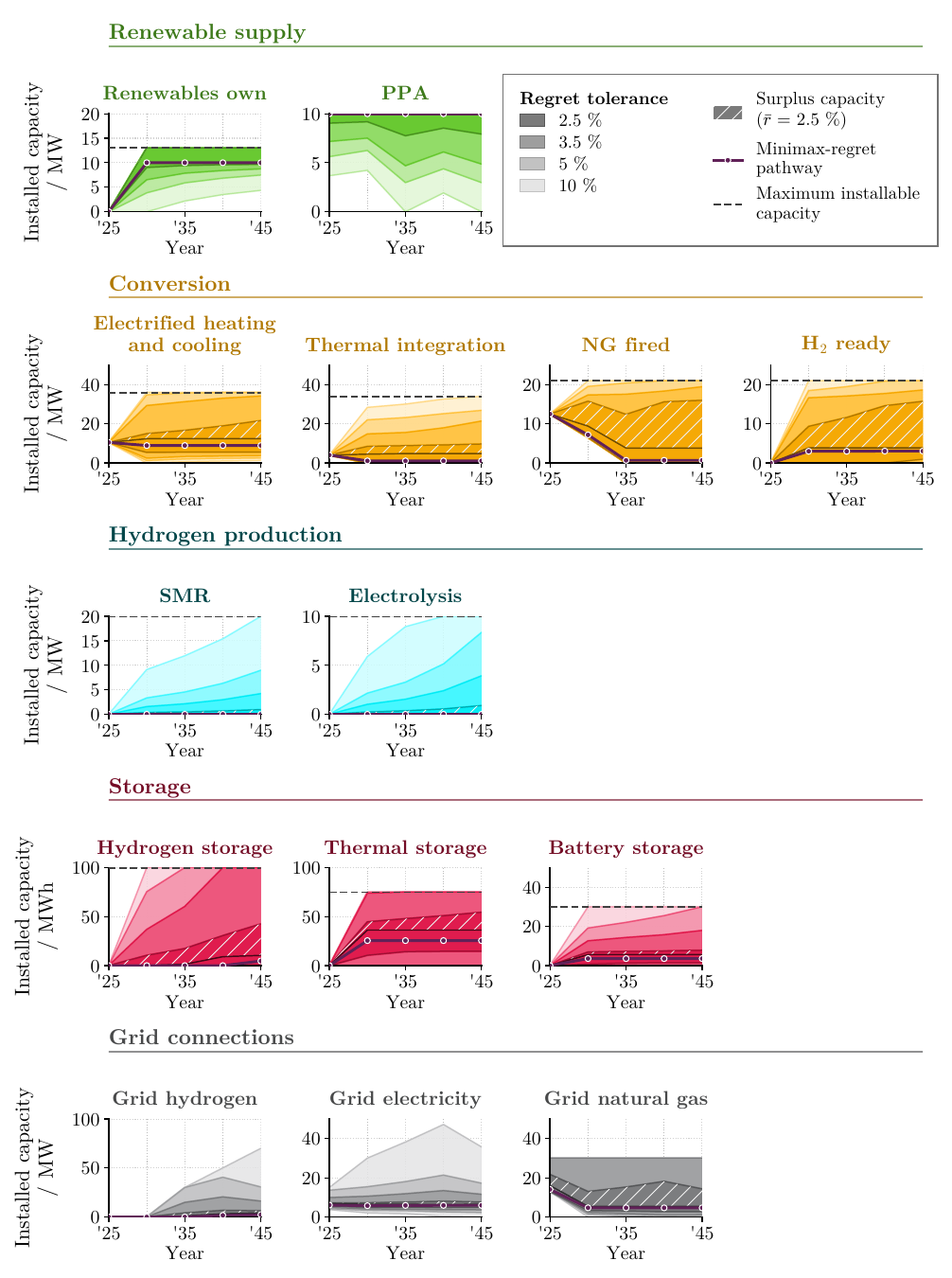}
\caption{Capacity ranges of 14 technology
families. Nested bands show the range of the installed and contracted capacities in every investment year at the four explored regret
tolerances $\bar r$. The hatched part of the darkest band
marks the surplus capacity at $\bar r=2.5\%$, i.e., the capacity that no
capacity-minimizing search reaches (Section~\ref{subsec:decision-support}). The purple line is the minimax-regret
pathway. The dashed grey line marks the maximum installable capacity of a technology family. Capacities change only in the investment years. The band edges and the pathway line are connected between investment years for readability and do not imply a continuous capacity expansion.}
\label{fig:range-clean}
\end{figure}

Figure~\ref{fig:range-clean} shows the low-regret capacity ranges of the 14 technology families at the four regret tolerances. From an investment perspective, the lower edge of a capacity range is the \emph{must-have} capacity, i.e., the minimum that every low-regret pathway installs. The upper edge is the \emph{can-have} capacity, i.e., the maximum that can be installed within the regret tolerance. In between, the lower part of the capacity range can be installed instead of other capacity. The hatched upper part is the surplus capacity, which is installed only in addition to other capacity.

The lower capacity range edges show what every low-regret pathway has in common: at 2.5\% regret tolerance, every pathway contracts a large share of the available PPA capacity already in 2025. From 2030 onwards, every pathway installs three further technology families: on-site renewables, thermal storage, and electrified heating and cooling. Every pathway also keeps an electricity and a natural gas grid connection over the whole planning horizon. All other decisions are optional apart from marginal additions. On-site hydrogen production is never strictly required, and a full exit from exclusively natural-gas-fired capacity is possible from 2035 onwards. At higher tolerances, the capacity ranges widen in both directions. The 2030 lower capacity range edges of on-site renewables and of electrified heating and cooling fall to zero at 10\%, i.e., they are not strictly required in 2030. The minimum cumulative added capacity per technology family and investment year is reported in Section~S2.3 of the Supplementary Information.

\subsubsection{Surplus capacity}
\label{subsubsec:results-surplus}

The individual capacity ranges are wide, especially for higher regret tolerances. The exploration reaches the upper edges of the capacity ranges mainly by search directions with positive coefficients, which reward additional capacity as long as the pathway stays within the regret tolerance, irrespective of whether the component is ever operated. Part of the width of the capacity ranges can therefore be capacity that is built without being operated, and including it would overstate the usable transition flexibility. We therefore report the utilization of the installed capacities and mark the surplus capacity of Section~\ref{subsec:decision-support}, i.e., the part of the capacity ranges that no capacity-minimizing search reaches.

\begin{figure}[htbp]
\centering
\includegraphics[width=0.99\textwidth]{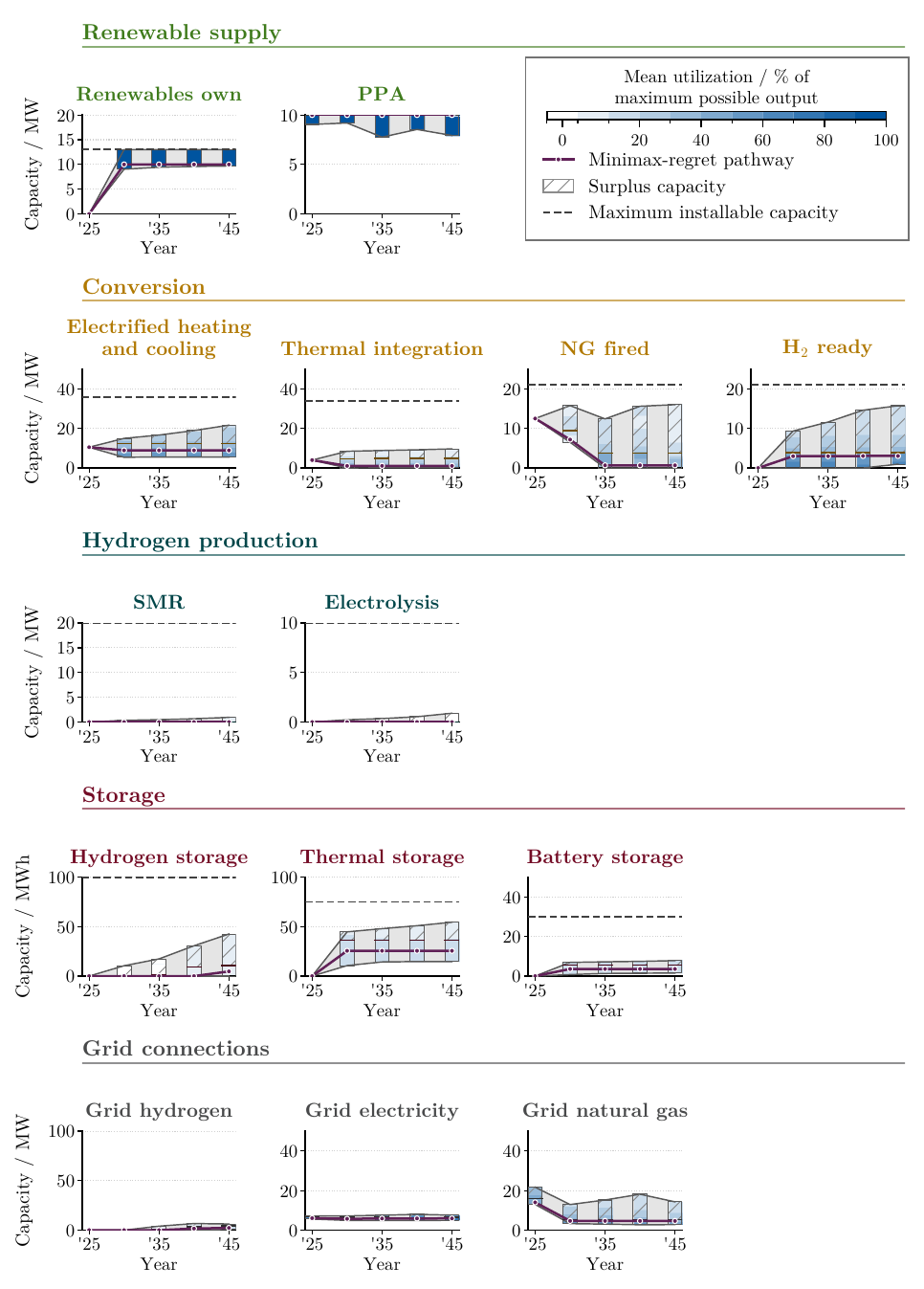}
\caption{Mean utilization of the installed capacities within the low-regret pathway space at the regret tolerance $\bar r=2.5\%$. The capacity range of each technology family (Figure~\ref{fig:range-clean}) is colored by the mean annual output of the explored pathways that install approximately that capacity, averaged across scenarios and divided by the maximum possible output of the installed capacity.}
\label{fig:utilization-band}
\end{figure}

Figure~\ref{fig:utilization-band} colors the capacity ranges at 2.5\% regret tolerance by the mean utilization of the installed capacity, i.e., the annual output relative to the maximum possible output (utilization of the individual components in Section~S2.4 of the Supplementary Information). In most technology families, the utilization decreases towards the upper edge of the capacity range. SMR and electrolysis capacities are barely operated. The hydrogen-ready capacities, in contrast, are highly utilized in the lower part of their capacity range. On-site renewable generation and the PPA capacity are highly utilized over their entire capacity ranges. For the PPA capacity, this follows from the pay-as-produced contract structure, because the site has to take the full generation of the contracted capacity irrespective of its own demand.

Figures~\ref{fig:range-clean} and~\ref{fig:utilization-band} hatch the surplus capacity. The surplus capacity coincides with the capacity values of low utilization. The surplus capacity is, however, an upper estimate of the capacity that is built without being operated. The capacity-minimizing searches also leave out capacity whose operation saves less than its investment costs. Additionally, a pathway that installs surplus capacity still operates it in the hours in which its operation is cheaper than the alternatives, e.g., an oversized electric boiler in hours with low electricity prices. The utilization within the hatched part therefore does not fall to zero in the conversion technology families. Hydrogen production capacities are almost entirely surplus. On-site renewables and PPA capacity hold no surplus capacity, because the capacity-minimizing searches install them up to the maximum installable capacity. The surplus capacities of the 3.5\%, 5\%, and 10\% spaces are identified in the same way and are shown in Section~S2.3 of the Supplementary Information.

\subsubsection{Summary}
\label{subsubsec:results-summary}

The results of Section~\ref{sec:results-space} answer \ref{rq:pathway-space}. The low-regret pathway spaces leave substantial transition flexibility in the sizing of the individual components. A large share of this transition flexibility is surplus capacity. Surplus capacity can be installed within the regret tolerance, e.g., as reserve capacity, but it is largely not operated under the considered scenarios. How much transition flexibility is gained with each increase of the regret tolerance is the subject of the next section.

\subsection{Regret--flexibility trade-off}
\label{sec:results-cost-of-flexibility}

\begin{figure}[htbp]
\centering
\includegraphics[width=0.55\textwidth]{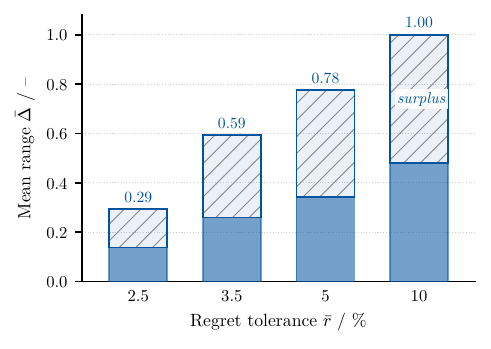}
\caption{The regret--flexibility trade-off: mean range
$\bar{\Delta}$ of the four explored low-regret pathway spaces. $\bar{\Delta}$ is normalized per variable to
the range of the $\bar r=10\%$ space. The full bar is the mean range $\bar{\Delta}$ over all searched directions. The solid lower part is the mean range $\bar{\Delta}$ over the capacity-minimizing searches
only. The hatched upper part is therefore the contribution of surplus capacity to the mean range, that is, of capacity that no capacity-minimizing search reaches.}
\label{fig:flexibility-tradeoff}
\end{figure}

Figure~\ref{fig:flexibility-tradeoff} shows the mean range $\bar{\Delta}$ of Section~\ref{subsec:decision-support} of the four explored pathway spaces against their regret tolerance. Increasing the regret tolerance widens the capacity ranges. The increase in transition flexibility is largest at low regret tolerances $\bar r$. At $\bar r = 2.5\%$, the pathway space has a mean range of about 30\% of the ranges at $ \bar r = 10\%$, and at $\bar r = 3.5\%$ already about 60\%. When taking capacity-minimizing searches only, the 2.5\% space likewise holds a share of 0.29 (solid part of the bars in Figure~\ref{fig:flexibility-tradeoff}). A higher regret tolerance thus widens the full capacity range as well as the surplus capacity range in a similar proportion.

Figure~\ref{fig:range-clean} shows how the normalized ranges separate over technologies. The share of the transition flexibility at 10\% regret tolerance that the 2.5\% space already holds differs strongly between technology families. For SMR and electrolysis, a substantial increase in transition flexibility is observed with higher regret tolerance, whereas for the natural-gas-fired and hydrogen-ready capacities the increase is much smaller.

These results answer \ref{rq:cost-of-flexibility} for the individual capacity ranges. Substantial transition flexibility is gained already at low regret tolerances. The capacity ranges are however one-dimensional projections of a pathway space, so not every combination of individually achievable capacity values is a low-regret pathway. The next section quantifies this combinability.

\subsection{Combinability of pathway decisions}
\label{sec:results-conditional}

The individual pathway decisions constrain one another. We quantify this dependence with the conditional range $R_i$ of Eq.~\eqref{eq:conditional-range}. The conditional range measures which share of their capacity ranges the other pathway decisions keep on average once a single pathway decision is fixed. We first examine the pathway decisions along the capacity ranges at a regret tolerance of 2.5\%. We then examine the dependence of the conditional range on the regret tolerance.

\begin{figure}[htbp]
\centering
\includegraphics[width=\textwidth]{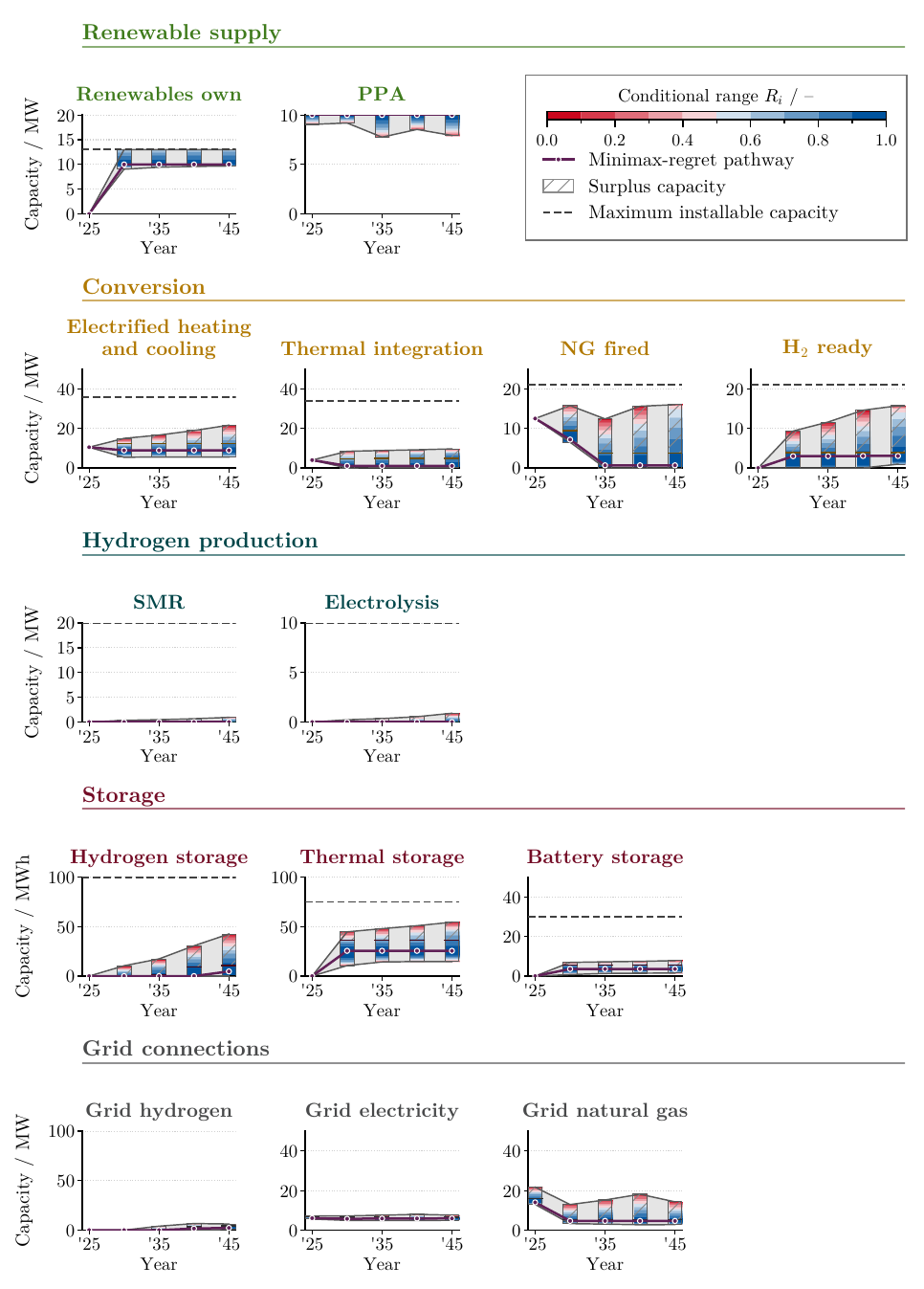}
\caption{Capacity ranges of the 14 technology families at the regret
tolerance $\bar r=2.5\%$, colored by the conditional range $R_i$. Values are shown as colored columns at the five
investment years. Hatching marks the
surplus capacity that no capacity-minimizing search reaches (Figure~\ref{fig:range-clean}). The purple line is the
minimax-regret pathway and the dashed grey line the maximum installable
capacity.}
\label{fig:flexibility-band}
\end{figure}

A single fixed pathway decision narrows the capacity ranges of the other pathway decisions. Figure~\ref{fig:flexibility-band} colors the capacity ranges by the conditional range. Blue regions mark fixed pathway decisions that leave the other capacity ranges almost unchanged, and red regions mark fixed pathway decisions that narrow the other capacity ranges strongly. The conditional range tends to be highest around the minimax-regret pathway (purple line) and to decrease when moving away from it towards the capacity range edges.

\begin{figure}[htbp]
\centering
\includegraphics[width=0.75\textwidth]{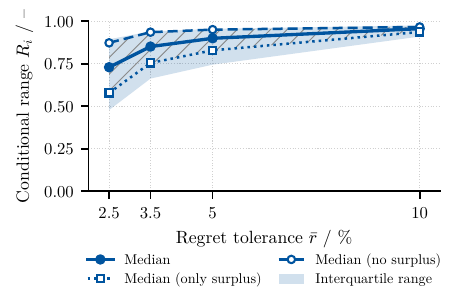}
\caption{Combinability of the pathway decisions across the explored regret
tolerances: median conditional range $R_i$. The solid curve takes
the median over all fixed pathway decisions of Figure~\ref{fig:flexibility-band}. The dotted curve takes the median over the
fixed pathway decisions that install surplus capacity, i.e., the hatched regions of
Figure~\ref{fig:flexibility-band}, and the dashed curve over the remaining pathway decisions. The shaded band is the interquartile range over all fixed pathway decisions. The hatching marks the gap between the two medians.}
\label{fig:combination-evolution}
\end{figure}

Figure~\ref{fig:combination-evolution} shows how the combinability of the pathway decisions changes across the four regret tolerances. Pathway decisions that install surplus capacity are less combinable. At 2.5\% regret tolerance, their median conditional range (dotted curve) is 0.58, compared with 0.87 for the remaining pathway decisions (dashed curve). The median conditional range (solid curve) rises from 0.73 at the lowest to 0.96 at the highest explored regret tolerance (the combinability plots at the higher regret tolerances are shown in Section~S2.3 of the Supplementary Information). This rise is carried by the pathway decisions that install surplus capacity, whose median rises from 0.58 to 0.94, whereas the median of the remaining pathway decisions rises from 0.87 to 0.97. At 10\% regret tolerance, the gap between the two medians has almost closed.

These results answer \ref{rq:residual-flexibility} and complete the answer to \ref{rq:cost-of-flexibility}. A single fixed pathway decision leaves the other pathway decisions about three quarters of their capacity ranges at 2.5\% regret tolerance. Pathway decisions that install surplus capacity are less combinable. A higher regret tolerance raises the median conditional range to 0.96 and thus yields, besides wider individual capacity ranges, a higher combinability of pathway decisions.
\clearpage
\section{Discussion}
\label{sec:discussion}

This section compares the principal findings with existing studies (Section~\ref{subsec:discussion-findings}) and discusses the limitations of the method and the case study together with directions for future work (Section~\ref{subsec:discussion-limitations}).

\subsection{Principal findings and comparison with existing studies}
\label{subsec:discussion-findings}

The minimax-regret pathway installs both grid connections and therefore keeps the choice of the imported energy carrier (natural gas or hydrogen) open, so it can switch the purchased energy carrier according to the price scenario. This hedging by operational flexibility is typical for minimax-regret energy system designs, e.g., as shown by \citet{Mavromatidis2018}.

The low-regret pathway space is defined through scenario-wise regret constraints around a stochastic minimax-regret pathway. Previous studies intersect near-optimal solution spaces that are explored separately per scenario \citep{Grochowicz2023,vanGreevenbroek2025b,vanGreevenbroek2025a}. Each of our exploration runs considers all scenarios simultaneously and therefore requires more computational effort than a single-scenario run. The intersection of separately explored near-optimal solution spaces, however, requires the facets of the convex hulls, which are intractable to enumerate in high dimensions \citep{Barber1996}. As a result, the previous studies intersect near-optimal solution spaces in at most eight dimensions. Our exploration likewise does not cover all dimensions but only their two-dimensional projections (Section~\ref{subsubsec:directions-termination}). Two-dimensional projections of separately explored near-optimal solution spaces would be cheap to intersect. This intersection, however, overestimates the projected low-regret pathway space. A point of the projection plane lies in the intersection if every scenario has a low-regret pathway that projects onto this point, but these pathways may differ between the scenarios in the remaining dimensions. The projected low-regret pathway space of Eq.~\eqref{eq:projected-space} instead contains a point only if a single pathway satisfies the regret constraints of all scenarios at once, as every boundary point of Eq.~\eqref{eq:directional-lp} does. The intersection of the projections is therefore only an outer bound of the projected low-regret pathway space. Our simultaneous exploration of all scenarios is preferable for high-dimensional spaces and few scenarios, whereas intersecting separately explored near-optimal solution spaces is preferable for lower-dimensional spaces and many scenarios. A further advantage of exploring high-dimensional pathway spaces is that it preserves the timing of the pathway decisions, which the decisions cumulated over the planning horizon in \citet{Price2017} and \citet{Sinha2024} do not.

Moving to the results of the case study, we observe that a higher regret tolerance yields a higher transition flexibility of two types: wider capacity ranges for individual pathway decisions and a higher combinability of pathway decisions. The capacity ranges are already wide at low regret tolerances. This agrees with the flatness of near-optimal spaces that, e.g., \citet{Pedersen2021} report. At higher regret tolerances the increase in the capacity ranges decreases, which follows from linear programming theory: the regret tolerance enters the linear program only on the right-hand side of the regret constraints (Eq.~\eqref{eq:low-regret-decision-space}). The optimal objective value of a linear program is a concave function of the right-hand side for maximization and a convex function for minimization \citep{Boyd2004}. Each capacity range, as the difference between the maximizing and the minimizing search (Section~\ref{subsec:decision-support}), is therefore a concave function of the regret tolerance. Intuitively, the first deviations from the minimax-regret pathway shift capacity between technologies with nearly equal costs, e.g., electric boilers replace natural-gas-fired boilers in the hours with the lowest electricity prices, whereas further deviations require increasingly expensive substitutes or reach the maximum installable capacity.

The conditional range is highest around the minimax-regret pathway and decreases towards the capacity range edges (Section~\ref{sec:results-conditional}). The fixed value $\tilde{x}_i$ is again a right-hand-side parameter of the searches, so the conditional range $R_i$ is a concave function of it and decreases in both directions from its maximum. The regret slack left for the other pathway decisions is largest at the minimax-regret pathway, which is consistent with the location of the maximum. The median conditional range rises from 0.73 to 0.96 between 2.5~\% and 10~\% regret tolerance, and the rise per percentage point of regret tolerance diminishes as for the capacity ranges. Unlike for the capacity ranges, this diminishing rise is not guaranteed by linear programming theory, since the conditional range relates the ranges after fixing to the full ranges, which widen as well. The rise is mainly driven by a higher combinability of pathway decisions that install surplus capacity. To our knowledge, the dependence of the combinability of pathway decisions itself on the cost slack or on the tolerated regret has not been quantified before. The conditional range averages the ratios of the ranges after fixing to the full ranges, as the indicator of \citet{vanGreevenbroek2025b} does.
Since the degree of exploration reaches at least 0.985 at every regret tolerance, these effects of the regret tolerance on the capacity ranges and on the conditional range are not biased by an incomplete exploration (Section~\ref{sec:results-space}). Our inner and outer approximations bound the low-regret pathway space from both sides. Monitoring the inner approximation alone would be weaker, since an inner approximation that no longer grows may still be far from converged \citep{Turan2026}.

\subsection{Limitations and outlook}
\label{subsec:discussion-limitations}

The considered energy system in this study is modeled as a two-stage stochastic linear program. Convexity is a precondition of our method, so modeling details of distributed energy systems that require integer variables, e.g., minimum installable capacities, economies of scale, and part-load behavior, are omitted. The two-stage formulation fixes one non-anticipative investment pathway and lets only operations adapt. Later investment decisions, however, can in practice still be adapted once more information about the uncertainty is revealed. As a consequence, the two-stage approach might be overconservative \citep{Huang2009,McDonald2024}. An extension of our method to a multi-stage formulation would be an interesting direction for future work.

We reduce the temporal representation of the energy system model in the operation year and in the planning horizon. One year of operation is aggregated to typical days. Methods for selecting the temporal aggregation of near-optimality analyses show that near-optimal solutions are less sensitive to the aggregation than the cost optimum \citep{Simarik2025}. The sensitivity of the low-regret pathway space to the aggregation should nevertheless be investigated in more detail. The finite optimization horizon is a known limitation of transition pathway formulations \citep{RuizManuel2026}. We extend the last investment period so that every investment within the horizon is charged its full annuities, but the operation after 2045 is represented by the input parameters of 2045. A residual end effect therefore remains in the last investment year.

The low-regret pathway space is valid only for the seven considered scenarios. Every additional scenario adds a regret constraint (Eq.~\eqref{eq:low-regret-decision-space}) and can therefore only narrow the low-regret pathway space at a given regret tolerance. Conversely, not all considered scenarios shape the low-regret pathway space (Section~S2.5 of the Supplementary Information). Since every scenario enters every directional optimization, dropping the scenarios that never bind would reduce the size of the linear programs without changing the explored space. Such a scenario reduction is a promising direction to reduce the computational effort of the method.

Capacity built without being operated is an undesirable feature of near-optimal alternatives \citep{Lombardi2025} that likewise enters the low-regret pathway space. Surplus capacity, however, may serve as a margin for higher demand \citep{Juergens2026} or extreme weather years \citep{Lombardi2020}. Our model does not consider the security of supply explicitly, i.e., capacity that an industrial site holds in reserve against outages appears as surplus capacity. The surplus capacity thus marks the capacity that the price scenarios alone do not justify, not capacity without any use. For these reasons, our results mark the surplus capacity instead of excluding it. Future work could include the security of supply explicitly in the model, so that reserve capacity no longer appears as surplus capacity.

The value of the regret tolerance is difficult to set a priori. Realized transitions of the European electricity systems deviated from the cost-optimal solution by more than 40\,\% over 30 years, far above typical slack values in MGA \citep{FossasTenas2026,Lombardi2025}. These deviations, however, reflect national policy and public acceptance, and no comparable evidence exists for distributed energy systems, e.g., for industrial sites. The regret tolerance therefore remains a case-specific choice. This choice does not have to be made before the exploration. Exploring the low-regret pathway space at several regret tolerances shows the transition flexibility that each regret tolerance yields (Section~\ref{sec:results-cost-of-flexibility}). Planners can then select the regret tolerance by weighing the tolerated regret against the transition flexibility that it yields. Planners can also refine the near-optimal space interactively, as in human-in-the-loop MGA \citep{Lombardi2025b}. Surrogate models of the energy system \citep{Neumann2023,Juergens2026} could accelerate such interactive refinement.
\section{Conclusion}
\label{sec:conclusion}

This paper explores low-regret transition pathway spaces around a stochastic minimax-regret pathway and quantifies the transition flexibility. The low-regret pathway space is defined through scenario-wise regret constraints, so that a transition pathway is low-regret only if it is near-optimal in every scenario. Instead of intersecting separately explored near-optimal solution spaces, the exploration method searches the intersection directly, directs each search toward the largest gaps between the inner and outer approximation, and bounds the remaining exploration error. The method is demonstrated on the transition planning of an industrial energy system under uncertain energy and CO$_2$ certificate prices. In this case study, at 2.5\% regret tolerance, every low-regret pathway contracts a large share of the available PPA capacity already in 2025. From 2030 onwards, every low-regret pathway installs on-site renewables, thermal storage, and electrified heating and cooling. The further pathway decisions leave substantial transition flexibility in the sizing of the individual components. A large share of this transition flexibility is surplus capacity. Surplus capacity can be installed within the regret tolerance, e.g., as reserve capacity, but it is largely not operated under the considered scenarios. The increase in transition flexibility is largest at low regret tolerances. A single fixed pathway decision narrows the capacity ranges of the other pathway decisions, in the median to about three quarters at 2.5\% regret tolerance. At 10\% regret tolerance, the other pathway decisions keep almost their full capacity ranges.

The relevant planning question for energy system transitions is therefore not only which pathway is optimal today, but which decisions keep transition flexibility open. Exploring the low-regret pathway space answers this question with the decisions that every low-regret pathway shares and the transition flexibility that remains around them.
\section*{Nomenclature}
\noindent
\begingroup
\small
\setlength{\tabcolsep}{0pt}
\renewcommand{\arraystretch}{1.12}
\begin{longtable}{@{}>{\raggedright\arraybackslash}p{3.4cm}>{\raggedright\arraybackslash}p{\dimexpr\textwidth-3.4cm\relax}@{}}
\toprule
\textbf{Symbol} & \textbf{Explanation} \\
\midrule
\endhead
\bottomrule
\endlastfoot

\multicolumn{2}{@{}l}{\textit{Sets}} \\[2pt]
$c \in \mathcal{C}$      & Set of energy system components \\
$y \in \mathcal{Y}$      & Set of investment years \\
$\omega \in \Omega$      & Set of scenarios \\[6pt]

$\mathcal{X}$             & Feasible set of the full stochastic energy system model \\
$\mathcal{X}_{\omega}$    & Scenario-specific feasible set for scenario $\omega$ \\
$\mathcal{X}_{\bar r}$    & Low-regret pathway space under regret tolerance $\bar r$ \\
$\mathcal{P}_{\bar r}$    & Projected low-regret pathway space \\
$\mathcal{L}$             & Inner approximation of the projected low-regret pathway space \\
$\mathcal{U}$             & Outer approximation of the projected low-regret pathway space \\
$\mathcal{E}$             & Set of explored pathways \\
$\mathcal{I}$             & Set of exploration variables with a nonzero range at $\bar r_{\max}$ \\
$\mathcal{J}_i$           & Set of exploration variables $j\neq i$ with a nonzero range \\
$\mathcal{E}_{ij}$        & Explored pathway set projected onto the plane of variables $i$ and $j$ \\[6pt]

\multicolumn{2}{@{}l}{\textit{Decision variables}} \\[2pt]
$\boldsymbol{x}$              & Full decision vector of the energy system model \\
$\boldsymbol{x}_{\omega}$     & Perfect-foresight pathway of scenario $\omega$ \\
$\boldsymbol{d}$              & First-stage investment decision vector \\
$d_{c,y}$                 & Investment decision for component $c$ in investment year $y$ \\
$\boldsymbol{o}_{\omega}$     & Second-stage operational decision vector in scenario $\omega$ \\
$R$                       & Worst-case regret across scenarios \\
$\hat{\boldsymbol{x}} \in \mathbb{R}^{m}$ & Projected pathway vector \\
$\hat{x}_{i}$                             & Exploration variable $i$, i.e., $i$-th component of $\hat{\boldsymbol{x}}$ \\[6pt]

\multicolumn{2}{@{}l}{\textit{Derived variables}} \\[2pt]
$\mathit{TAC}_{\omega}(\boldsymbol{x})$      & Scenario-dependent total annualized cost of transition pathway $\boldsymbol{x}$ \\
$\mathit{COST}_{y,\omega}(\boldsymbol{x})$   & Annual cost in investment year $y$ and scenario $\omega$ \\
$r_{\omega}(\boldsymbol{x})$                 & Scenario-wise relative regret of transition pathway $\boldsymbol{x}$ \\
$r^{*}$                                  & Minimum achievable worst-case regret \\
$\beta(\boldsymbol{n})$                      & Support value of $\mathcal{P}_{\bar r}$ in direction $\boldsymbol{n}$ \\
$A_{L},\;A_{U}$                          & Aggregated projected areas of the inner and outer approximation \\
$\mathit{DoE}_{\mathrm{VESA}}$           & VESA-based projected-area degree of exploration \\
$\beta_{\mathcal{L}}(\boldsymbol{n}_e),\;\beta_{\mathcal{U}}(\boldsymbol{n}_e)$ & Support values of the inner and outer approximation along the normal $\boldsymbol{n}_e$ of projection edge $e$ \\
$\gamma_e$                               & Support gap of projection edge $e$ \\
$\hat{x}_i^{\max},\; \hat{x}_i^{\min}$   & Maximum and minimum of exploration variable $i$ across the explored pathway set $\mathcal{E}$ \\
$\Delta_i$                               & Range of exploration variable $i$ across the explored pathway set $\mathcal{E}$ \\
$\bar{\Delta}(\bar r)$                   & Mean range at regret tolerance $\bar r$, normalized to $\bar r_{\max}$ \\
$\tilde{x}_i$                            & Fixed value of exploration variable $i$ \\
$\Delta_j(\tilde{x}_i)$                  & Range of exploration variable $j$ after fixing variable $i$ at $\tilde{x}_i$ \\
$R_i(\tilde{x}_i)$                       & Conditional range \\[6pt]

\multicolumn{2}{@{}l}{\textit{Parameters}} \\[2pt]
$\mathbf{S}$               & Projection matrix \\
$\mathit{TAC}_{\omega}^{*}$ & Perfect-foresight optimal total annualized cost in scenario $\omega$ \\
$\alpha_{y}$               & Discount factor for investment year $y$ \\
$i$                        & Interest rate for discounting \\
$y_{\mathrm{ref}}$         & Reference year for discounting \\
$m$                        & Number of exploration variables \\
$\boldsymbol{n}$           & Search direction in the projected pathway space \\
$\boldsymbol{n}_h$         & Searched direction $h$, i.e., a solved search direction \\
$H$                        & Number of searched directions \\
$\bar r$                   & Regret tolerance \\
$\bar r_{\max}$            & Highest explored regret tolerance \\
\end{longtable}
\endgroup
\section*{Data availability}

The source code of the method will be made openly available under an MIT license upon publication of this article. The case study was solved with Gurobi~12.0.0, which requires a separate license. The method can also be used with suitable open-source solvers. A supporting web tool for exploring low-regret pathways is available at \url{https://apps.ltt.rwth-aachen.de/transition-pathways}.

\section*{CRediT authorship contribution statement}

\textbf{Hendrik Schricker:} Conceptualization, Data curation, Formal analysis, Funding acquisition, Investigation, Methodology, Supervision, Project administration, Software, Validation, Visualization, Writing -- original draft, Writing -- review \& editing. \textbf{Maria Louise Schmitt:} Methodology, Software, Data curation, Visualization. \textbf{Philipp Kleusberg:} Methodology, Software, Data curation, Visualization. \textbf{Valentin Dinges:} Methodology, Software. \textbf{Benedikt Schuler:} Methodology, Software, Data curation. \textbf{Frederike Kuperjans:} Conceptualization, Funding acquisition, Methodology, Writing -- review \& editing. \textbf{Boyung Jürgens:} Conceptualization, Funding acquisition, Methodology, Supervision, Writing -- review \& editing. \textbf{Niklas von der Assen:} Conceptualization, Funding acquisition, Resources, Supervision, Writing -- review \& editing.

\section*{Declaration of competing interest}

The authors declare that they have no known competing financial interests or personal relationships that could have appeared to influence the work reported in this paper.

\section*{Funding}

We gratefully acknowledge funding from the German Federal Ministry of Economic Affairs and Energy (BMWE) for the project ``FlexiT'' (03EI1113A).

\section*{Acknowledgments}

The authors thank Conrad Luka and Jonas Padberg for contributing to the model framework. Simulations were performed with computing resources granted by RWTH Aachen University under project ID rwth1542. We thank Karan Anand and Yannick Becker for many fruitful discussions.

\section*{Declaration of generative AI and AI-assisted technologies in the manuscript preparation process}

During the preparation of this work, the authors used Claude (Anthropic) in order to revise the structure, wording, and clarity of text drafted by the authors, and to assist with programming tasks, including the code for data analysis and figure generation. After using this tool, the authors reviewed and edited the content as needed and take full responsibility for the content of the published article.

\begingroup
\sloppy
\bibliographystyle{elsarticle-harv-nourl}
\bibliography{references}
\endgroup

\end{document}


\begin{center}
{\large Supplementary Information}\\[10pt]
{\LARGE Exploring Low-Regret Energy System Transition Pathways with Modeling to Generate Alternatives}\\[14pt]
{\large Hendrik Schricker, Maria Louise Schmitt, Philipp Kleusberg, Valentin Dinges, Benedikt Schuler, Frederike Kuperjans, Boyung Jürgens, Niklas von der Assen}\\[8pt]
\textit{\small Institute of Technical Thermodynamics, RWTH Aachen University, Schinkelstr.~8, 52062 Aachen, Germany}
\end{center}

\noindent\rule{\linewidth}{0.4pt}

\begin{singlespace}
\setcounter{tocdepth}{2}
\tableofcontents
\end{singlespace}
\newpage

\section{Supplementary Case Study}
\label{sec:si-casestudy}

This section documents all techno-economic parameters and assumptions of the industrial energy system case study. Section~\ref{sec:si-demand} documents the energy demands, Sections~\ref{sec:si-conversion} and~\ref{sec:si-storage-tech} the assumptions and parameters for the conversion and storage components, and Section~\ref{sec:si-procurement} for energy purchases and the price scenarios.

\subsection{Energy demands}
\label{sec:si-demand}

\begin{figure}[H]
\centering
\includegraphics[width=\textwidth]{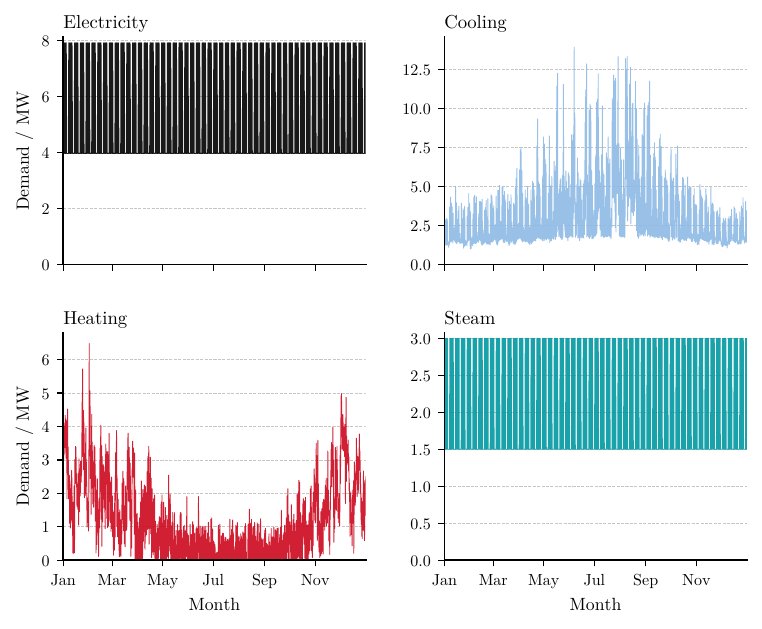}
\caption{Hourly demand time series per energy carrier. All demand time series are invariant across scenarios and investment years.}
\label{fig:si-demand-profiles}
\end{figure}

The energy system supplies electricity, heating, cooling, and steam to an industrial site. Demand time series are based on \citet{Baumgaertner2019}. Table~\ref{tab:si-demand} reports the annual average demands of the four energy carriers. The annual average demands are assumed to be constant over the investment years and the scenarios. Of the energy supplied as steam, the industrial production system is assumed to return 40\% as recoverable excess heat. Figure~\ref{fig:si-demand-profiles} shows the demand time series of the four energy carriers for one year.

\begin{table}[htbp]
\begin{singlespace}
\centering
\caption{Annual average energy demands taken from \cite{Baumgaertner2019}.}
\label{tab:si-demand}
\begin{tabular}{L{4.5cm}R{4.5cm}}
\toprule
Energy carrier & Annual average demand (MW) \\
\midrule
Electricity     & 5.45 \\
Cooling  & 3.12 \\
Heating  & 1.22 \\
Steam & 1.99 \\
\bottomrule
\end{tabular}
\end{singlespace}
\end{table}

\newpage
\subsection{Conversion components}
\label{sec:si-conversion}

This section documents the investment costs, lifetimes, and maintenance costs (Section~\ref{sec:si-conv-economics}), the conversion coefficients (Section~\ref{sec:si-conv-coefficients}), the availability factors (Section~\ref{sec:si-availability-factors}), and the existing and installable capacities (Section~\ref{sec:si-conv-capacities}) of the conversion components. In the first investment year 2025, only the existing capacities are available. All further investment options are available from 2030 onwards. The tables therefore report the investment years 2030 to 2045.

\subsubsection{Investment costs, lifetimes, and maintenance costs}
\label{sec:si-conv-economics}

Table~\ref{tab:si-invest-conv} reports the specific investment costs for all installable conversion components. Investment costs are annualized over an assumed economic period of 10 years. All costs are discounted with interest rate $i = 5\%$ \citep{Reinert2023} to the reference year $y_{\mathrm{ref}} = 2025$. Each investment year represents the years until the next investment year, and the annual costs are aggregated over this period with the corresponding annuity present value factor. The investment years 2025 to 2040 represent five years each. The last investment year 2045 represents ten years, so that every investment made within the optimization horizon is charged its full annuities. Annuities are charged only for the years of the annuity period that fall within an investment period. Maintenance costs and the capacity charges of the grid connections are charged for every investment period in which a component is available. Table~\ref{tab:si-lifetime-conv} reports the lifetimes of the conversion components. Table~\ref{tab:si-maintenance-conv} reports the yearly maintenance costs of the conversion components.
\newpage
\begingroup
\singlespacing\sitablefont
\setlength{\tabcolsep}{4pt}
\begin{xltabular}{\textwidth}{>{\raggedright\arraybackslash}Xrrrrc}
\caption{Specific investment costs (in kEUR/MW) for conversion components by investment year. If only one value is reported, this value applies to all investment years. CCS: carbon capture and storage, CHP: combined heat and power, ELEC: electric, H$_2$R: hydrogen-ready, NG: natural gas, PEM: proton exchange membrane, SMR: steam methane reformer.}
\label{tab:si-invest-conv}\\
\toprule
Component & 2030 & 2035 & 2040 & 2045 & Reference \\
\midrule
\endfirsthead
\multicolumn{6}{l}{\itshape Table~\ref{tab:si-invest-conv} (continued)}\\
\toprule
Component & 2030 & 2035 & 2040 & 2045 & Reference \\
\midrule
\endhead
\midrule
\multicolumn{6}{r}{\itshape continued on next page}\\
\endfoot
\bottomrule
\multicolumn{6}{@{}p{\textwidth}@{}}{\footnotesize\raggedright
\tmark{a}~\cite{ChillerRef2022},
\tmark{b}~\cite{DEA_electricity},
\tmark{c}~\cite{Gorbach2022} (hydrogen-ready components: natural gas component costs increased by 20\%),
\tmark{d}~own assumption,
\tmark{e}~\cite{Tartiere2017},
\tmark{f}~\cite{DEA_fuels},
\tmark{g}~\cite{ACS_SMR},
\tmark{h}~\cite{Valancius2015},
\tmark{i}~\cite{DEA_industrial},
\tmark{*}~derived from the availability factor and the collector efficiency, see Section~\ref{sec:si-availability-factors}.}\\
\endlastfoot
Chiller (absorption) & \multicolumn{4}{c}{241} & \tmark{a} \\
Chiller (compression) & \multicolumn{4}{c}{170} & \tmark{a} \\
CHP (H$_2$R) & 1148 & 1132 & 1117 & 1101 & \tmark{b,c} \\
CHP (NG) & 957 & 944 & 930 & 917 & \tmark{b} \\
Heat exchanger & \multicolumn{4}{c}{50} & \tmark{d} \\
Heat pump (air source) & \multicolumn{4}{c}{910} & \tmark{b} \\
Heat pump (excess heat) & \multicolumn{4}{c}{810} & \tmark{b} \\
Organic Rankine cycle & \multicolumn{4}{c}{1500} & \tmark{e} \\
PEM electrolyzer & 1624 & 1400 & 1177 & 965 & \tmark{f} \\
Photovoltaics & 570 & 515 & 460 & 435 & \tmark{b} \\
SMR & 809 & 793 & 777 & 761 & \tmark{g} \\
SMR with CCS & 1256 & 1237 & 1219 & 1200 & \tmark{g} \\
Solar thermal & 568 & 563 & 558 & 553 & \tmark{b,h,*} \\
Steam boiler (ELEC) & \multicolumn{4}{c}{80} & \tmark{i} \\
Steam boiler (H$_2$R) & \multicolumn{4}{c}{54} & \tmark{i,c} \\
Steam boiler (NG) & \multicolumn{4}{c}{45} & \tmark{i} \\
Water boiler (ELEC) & \multicolumn{4}{c}{70} & \tmark{i} \\
Water boiler (H$_2$R) & \multicolumn{4}{c}{43.2} & \tmark{i,c} \\
Water boiler (NG) & \multicolumn{4}{c}{36} & \tmark{i} \\
Wind & 1150 & 1130 & 1110 & 1100 & \tmark{b} \\
\end{xltabular}
\endgroup

\begin{table}[H]
\begin{singlespace}
\centering
\sitablefont
\setlength{\tabcolsep}{4pt}
\caption{Lifetimes (in years) for conversion components. CCS: carbon capture and storage, CHP: combined heat and power, ELEC: electric, H$_2$R: hydrogen-ready, NG: natural gas, PEM: proton exchange membrane, SMR: steam methane reformer.}
\label{tab:si-lifetime-conv}
\sbox0{\begin{tabular}{lrc}
\toprule
Component & Lifetime & Reference \\
\midrule
Chiller (absorption) & 25 & \tmark{a} \\
Chiller (compression) & 25 & \tmark{a} \\
CHP (H$_2$R) & 25 & \tmark{b} \\
CHP (NG) & 25 & \tmark{b} \\
Heat exchanger & 20 & \tmark{c} \\
Heat pump (air source) & 25 & \tmark{b} \\
Heat pump (excess heat) & 25 & \tmark{b} \\
Organic Rankine cycle & 30 & \tmark{d} \\
PEM electrolyzer & 25 & \tmark{e} \\
Photovoltaics & 40 & \tmark{b} \\
SMR & 25 & \tmark{f} \\
SMR with CCS & 25 & \tmark{f} \\
Solar thermal & 30 & \tmark{b} \\
Steam boiler (ELEC) & 25 & \tmark{g} \\
Steam boiler (H$_2$R) & 25 & \tmark{g} \\
Steam boiler (NG) & 25 & \tmark{g} \\
Water boiler (ELEC) & 25 & \tmark{g} \\
Water boiler (H$_2$R) & 25 & \tmark{g} \\
Water boiler (NG) & 25 & \tmark{g} \\
Wind & 30 & \tmark{b} \\
\bottomrule
\end{tabular}}
\begin{minipage}{\wd0}
\usebox0
\par\smallskip
{\footnotesize\raggedright \tmark{a}~\cite{Grosse2017}, \tmark{b}~\cite{DEA_electricity}, \tmark{c}~own assumption, \tmark{d}~\cite{Reinert2023}, \tmark{e}~\cite{DEA_fuels}, \tmark{f}~\cite{ACS_SMR}, \tmark{g}~\cite{DEA_industrial}.\par}
\end{minipage}
\end{singlespace}
\end{table}

\clearpage

\begingroup
\singlespacing\sitablefont
\setlength{\tabcolsep}{4pt}
\begin{xltabular}{\textwidth}{Yrrrrc}
\caption{Yearly maintenance costs (in kEUR/MW per year, or as a percentage of the specific investment cost) for conversion components by investment year. If only one value is reported, this value applies to all investment years. CCS: carbon capture and storage, CHP: combined heat and power, ELEC: electric, H$_2$R: hydrogen-ready, NG: natural gas, PEM: proton exchange membrane, SMR: steam methane reformer.}
\label{tab:si-maintenance-conv}\\
\toprule
Component & 2030 & 2035 & 2040 & 2045 & Reference \\
\midrule
\endfirsthead
\multicolumn{6}{l}{\itshape Table~\ref{tab:si-maintenance-conv} (continued)}\\
\toprule
Component & 2030 & 2035 & 2040 & 2045 & Reference \\
\midrule
\endhead
\midrule
\multicolumn{6}{r}{\itshape continued on next page}\\
\endfoot
\bottomrule
\endlastfoot
Chiller (absorption) & \multicolumn{4}{c}{2.5\% of invest} & \tmark{a} \\
Chiller (compression) & \multicolumn{4}{c}{1.0\% of invest} & \tmark{a} \\
CHP (H$_2$R) & 11.87 & 11.61 & 11.36 & 11.10 & \tmark{b,c} \\
CHP (NG) & 9.89 & 9.68 & 9.46 & 9.25 & \tmark{b} \\
Heat exchanger & \multicolumn{4}{c}{1.0\% of invest} & \tmark{a} \\
Heat pump (air source) & \multicolumn{4}{c}{2.13} & \tmark{b} \\
Heat pump (excess heat) & \multicolumn{4}{c}{2.13} & \tmark{b} \\
Organic Rankine cycle & \multicolumn{4}{c}{2.0\% of invest} & \tmark{d} \\
PEM electrolyzer & \multicolumn{4}{c}{2.0\% of invest} & \tmark{e} \\
Photovoltaics & 8.10 & 7.70 & 7.30 & 7.05 & \tmark{b} \\
SMR & \multicolumn{4}{c}{4.7\% of invest} & \tmark{f} \\
SMR with CCS & \multicolumn{4}{c}{3.0\% of invest} & \tmark{f} \\
Solar thermal & \multicolumn{4}{c}{0.285} & \tmark{b} \\
Steam boiler (ELEC) & 1.06 & 1.03 & 1.00 & 0.98 & \tmark{g} \\
Steam boiler (H$_2$R) & 2.28 & 2.22 & 2.16 & 2.10 & \tmark{g,c} \\
Steam boiler (NG) & 1.90 & 1.85 & 1.80 & 1.75 & \tmark{g} \\
Water boiler (ELEC) & 0.92 & 0.90 & 0.87 & 0.85 & \tmark{g} \\
Water boiler (H$_2$R) & 2.28 & 2.22 & 2.16 & 2.10 & \tmark{g,c} \\
Water boiler (NG) & 1.90 & 1.85 & 1.80 & 1.75 & \tmark{g} \\
Wind & 16.66 & 16.31 & 15.97 & 15.78 & \tmark{b} \\
\end{xltabular}
\vspace{-6pt}
{\footnotesize\raggedright
\tmark{a}~own assumption,
\tmark{b}~\cite{DEA_electricity},
\tmark{c}~\cite{Gorbach2022},
\tmark{d}~\cite{Reinert2023},
\tmark{e}~\cite{DEA_fuels},
\tmark{f}~\cite{ACS_SMR},
\tmark{g}~\cite{DEA_industrial}.\par}
\endgroup
\newpage
\subsubsection{Conversion coefficients}
\label{sec:si-conv-coefficients}

The conversion coefficients give the amount of each energy carrier a component consumes or produces. Table~\ref{tab:si-efficiency} reports the conversion coefficients normalized to 1~MW output of each component's reference energy carrier. All coefficients are constant within an investment year, except for the electricity input coefficient of the air-source heat pump, which is scaled by an hourly multiplier to consider the influence of varying ambient temperatures on the coefficient of performance (COP).

\begingroup
\singlespacing\footnotesize
\setlength{\tabcolsep}{2pt}
\begin{xltabular}{\textwidth}{>{\hsize=1.35\hsize}Y>{\hsize=0.85\hsize}Y>{\hsize=0.80\hsize}Yrrrr}
\caption{Conversion coefficients by investment year. All coefficients are normalized to 1~MW of the component's reference energy carrier and give the MW of each energy carrier consumed (negative) or produced (positive) per MW of this reference output. The hydrogen-ready components consume the energy carrier fuel, which is either natural gas or hydrogen. CCS: carbon capture and storage, CHP: combined heat and power, ELEC: electric, H$_2$R: hydrogen-ready, NG: natural gas, PEM: proton exchange membrane, SMR: steam methane reformer.}
\label{tab:si-efficiency}\\
\toprule
Component & Reference energy carrier & Energy carrier & 2030 & 2035 & 2040 & 2045 \\
\midrule
\endfirsthead
\multicolumn{7}{l}{\itshape Table~\ref{tab:si-efficiency} (continued)}\\
\toprule
Component & Reference energy carrier & Energy carrier & 2030 & 2035 & 2040 & 2045 \\
\midrule
\endhead
\midrule
\multicolumn{7}{r}{\itshape continued on next page}\\
\endfoot
\bottomrule
\multicolumn{7}{@{}p{\dimexpr\textwidth-2\tabcolsep\relax}@{}}{\footnotesize\raggedright
\textit{References.}~\tmark{a}~\cite{ChillerRef2022},
\tmark{b}~\cite{DEA_electricity},
\tmark{c}~own assumption,
\tmark{d}~\cite{Reinert2023},
\tmark{e}~\cite{DEA_fuels},
\tmark{f}~\cite{Antonini2020},
\tmark{g}~\cite{DEA_industrial},
\tmark{*}~annual-average electricity input coefficient, scaled by an hourly multiplier to consider the effect of varying ambient temperatures on the coefficient of performance,
\tmark{**}~receives a CO$_2$ certificate credit for its captured emissions, documented in the text following the table.}\\
\endlastfoot
Chiller (absorption)\tmark{a} & Cold & Heat & \multicolumn{4}{c}{$-0.714$} \\
Chiller (compression)\tmark{a} & Cold & Electricity & \multicolumn{4}{c}{$-0.141$} \\
\multirow{2}{=}{CHP (H$_2$R)\tmark{b}} & \multirow{2}{=}{Electricity} & Heat & 1.010 & 0.998 & 0.986 & 0.974 \\*
 & & Fuel & $-2.128$ & $-2.116$ & $-2.105$ & $-2.094$ \\
\multirow{2}{=}{CHP (NG)\tmark{b}} & \multirow{2}{=}{Electricity} & Heat & 1.010 & 0.998 & 0.986 & 0.974 \\*
 & & Natural gas & $-2.128$ & $-2.116$ & $-2.105$ & $-2.094$ \\
Heat exchanger\tmark{c} & Heat & Steam & \multicolumn{4}{c}{$-1.020$} \\
Heat pump (air source)\tmark{b,*} & Heat & Electricity & $-0.313$ & $-0.308$ & $-0.303$ & $-0.299$ \\
\multirow{2}{=}{Heat pump (excess heat)\tmark{b}} & \multirow{2}{=}{Heat} & Electricity & $-0.208$ & $-0.206$ & $-0.204$ & $-0.202$ \\*
 & & Excess heat & $-0.833$ & $-0.836$ & $-0.838$ & $-0.840$ \\
Organic Rankine cycle\tmark{d} & Electricity & Heat & \multicolumn{4}{c}{$-8.43$} \\
\multirow{2}{=}{PEM electrolyzer\tmark{e}} & \multirow{2}{=}{Hydrogen} & Electricity & $-1.710$ & $-1.665$ & $-1.620$ & $-1.565$ \\*
 & & Excess heat & $+0.453$ & $+0.410$ & $+0.368$ & $+0.311$ \\
\multirow{2}{=}{SMR\tmark{f}} & \multirow{2}{=}{Hydrogen} & Natural gas & \multicolumn{4}{c}{$-1.299$} \\*
 & & Electricity & \multicolumn{4}{c}{$+0.041$} \\
\multirow{2}{=}{SMR with CCS\tmark{f,**}} & \multirow{2}{=}{Hydrogen} & Natural gas & \multicolumn{4}{c}{$-1.299$} \\*
 & & Electricity & \multicolumn{4}{c}{$+0.016$} \\
Steam boiler (ELEC)\tmark{g} & Steam & Electricity & \multicolumn{4}{c}{$-1.010$} \\
Steam boiler (H$_2$R)\tmark{g} & Steam & Fuel & $-1.075$ & $-1.075$ & $-1.075$ & $-1.070$ \\
Steam boiler (NG)\tmark{g} & Steam & Natural gas & $-1.075$ & $-1.075$ & $-1.075$ & $-1.070$ \\
Water boiler (ELEC)\tmark{g} & Heat & Electricity & \multicolumn{4}{c}{$-1.010$} \\
Water boiler (H$_2$R)\tmark{g} & Heat & Fuel & $-1.064$ & $-1.064$ & $-1.064$ & $-1.058$ \\
Water boiler (NG)\tmark{g} & Heat & Natural gas & $-1.064$ & $-1.064$ & $-1.064$ & $-1.058$ \\
\end{xltabular}
\endgroup

\begin{figure}[H]
\centering
\includegraphics[width=\textwidth]{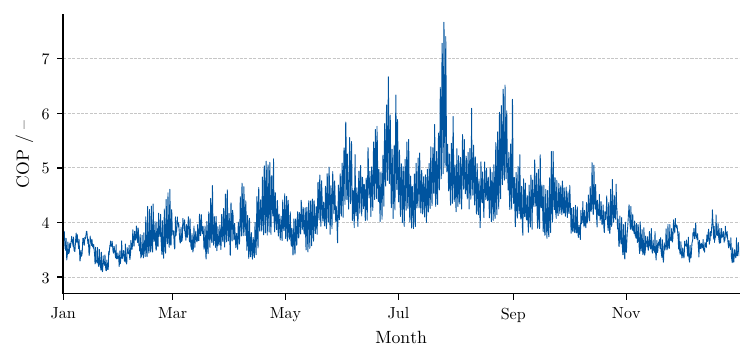}
\caption{Assumed hourly COP profile of the air-source heat pump for 2045.}
\label{fig:si-heat-pump-cop}
\end{figure}

The electricity input of the air-source heat pump in Table~\ref{tab:si-efficiency} is $-1/\mathrm{COP}$ with the nominal COP \citep{DEA_electricity}, which increases from 3.20 in 2030 to 3.35 in 2045 and refers to an ambient temperature of $-1$\,$^\circ$C. An hourly multiplier scales this coefficient with the ambient temperature time series obtained from Renewables.ninja \citep{Staffell2016, Pfenninger2016}. The multiplier assumes a constant supply temperature of 70\,$^\circ$C and a constant ratio of the actual COP to the Carnot COP, which is the supply temperature divided by the temperature lift between supply and ambient temperature. Because the temperature lift grows as the ambient temperature falls, the electricity demand per unit of heat exceeds the nominal value in hours colder than $-1$\,$^\circ$C and falls below it in warmer hours. Figure~\ref{fig:si-heat-pump-cop} shows the resulting hourly COP profile for 2045.

SMR with CCS receives a CO$_2$ certificate credit of 140.2~kg~CO$_2$ per MWh of hydrogen. Its total natural gas input of 1.299~MWh per MWh of hydrogen (Table~\ref{tab:si-efficiency}) would be subject to a certificate charge of 260.8~kg~CO$_2$ at the emission factor of Section~\ref{sec:si-price-scenarios}, of which 120.6~kg~CO$_2$ remain as direct emissions. The credit is the difference between these two values. The direct emissions reflect pre-combustion capture of 90\% of the CO$_2$ in the syngas, corresponding to an overall capture rate of 55\% for the SMR because the furnace flue gas is not captured \citep{Antonini2020}.

\subsubsection{Availability factors}
\label{sec:si-availability-factors}

\begin{figure}[H]
\centering
\includegraphics[width=\textwidth]{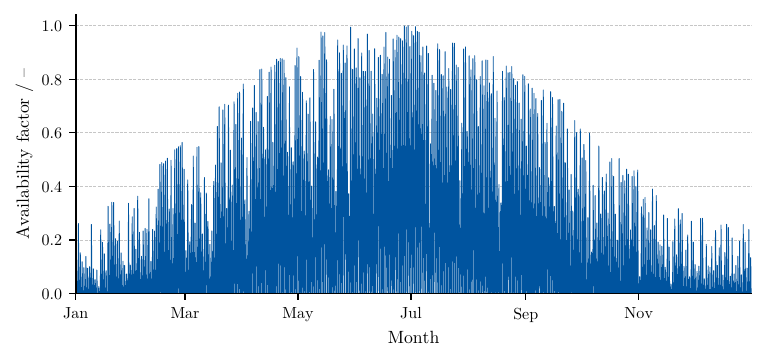}
\caption{Assumed hourly availability factor of the solar thermal collector.}
\label{fig:si-solar-thermal-availability}
\end{figure}

The availability factors give the hourly limit on the output of photovoltaics, wind, and solar thermal as a fraction of their installed capacity. The availability factors of the on-site photovoltaic, wind, and of the PPA units are taken from \citet{Schricker2026b}. The availability factor of the solar thermal collector is the irradiance time series obtained from Renewables.ninja \citep{Staffell2016, Pfenninger2016}, normalized by its annual maximum. Figure~\ref{fig:si-solar-thermal-availability} shows the resulting hourly availability factor.

The investment and maintenance costs of the solar thermal collector in Tables~\ref{tab:si-invest-conv} and \ref{tab:si-maintenance-conv} follow from this profile. The collector efficiencies increasing from 0.48 to 0.495 for 2030--2045 \citep{DEA_electricity} give the collector area required for one MW of thermal peak capacity. The investment cost is this area multiplied by the size-dependent installation cost $740.07\,A^{-0.14}$~EUR/m$^2$ of \citet{Valancius2015}, evaluated at the area $A$ of one MW of thermal peak capacity in m$^2$. The maintenance cost is the specific maintenance cost of 0.2~EUR/MWh \citep{DEA_electricity} applied to the maximum possible annual yield of the profile.

\subsubsection{Existing and installable capacities}
\label{sec:si-conv-capacities}

Table~\ref{tab:si-existing} reports the existing capacities. The existing technology set is that of the case study of \cite{Schricker2026b}. The capacities of these units result from a cost optimization for the investment year 2020, rounded up to 0.5~MW. The existing capacities phase out within the optimization horizon and can be replaced by newly installable capacities. Table~\ref{tab:si-capacity-conv} reports the assumed maximum installable capacity of each conversion component.

\begin{table}[H]
\begin{singlespace}
\centering
\sitablefont
\setlength{\tabcolsep}{4pt}
\caption{Assumed existing capacities for energy conversion. All existing capacities are available for operation from the start of the optimization horizon and phase out at the stated end of life. CHP: combined heat and power, NG: natural gas.}
\label{tab:si-existing}
\begin{tabular}{lrc}
\toprule
Existing unit & Capacity (MW) & End of life \\
\midrule
Chiller (absorption) & 4.0 & 2029 \\
Chiller (compression) & 10.5 & 2029 \\
CHP (NG) & 6.0 & 2029 \\
Steam boiler (NG) & 3.0 & 2034 \\
Water boiler (NG) & 3.5 & 2034 \\
\bottomrule
\end{tabular}
\end{singlespace}
\end{table}

\begin{table}[H]
\begin{singlespace}
\centering
\sitablefont
\setlength{\tabcolsep}{4pt}
\caption{Maximum installable capacity (in MW) of the conversion components. CCS: carbon capture and storage, CHP: combined heat and power, ELEC: electric, H$_2$R: hydrogen-ready, NG: natural gas, PEM: proton exchange membrane, SMR: steam methane reformer.}
\label{tab:si-capacity-conv}
\begin{tabular}{lr}
\toprule
Component & Maximum installable capacity \\
\midrule
Chiller (absorption) & 16 \\
Chiller (compression) & 16 \\
CHP (H$_2$R) & 6 \\
CHP (NG) & 6 \\
Heat exchanger & 5 \\
Heat pump (air source) & 5 \\
Heat pump (excess heat) & 5 \\
Organic Rankine cycle & 8 \\
PEM electrolyzer & 10 \\
Photovoltaics & 5 \\
SMR & 10 \\
SMR with CCS & 10 \\
Solar thermal & 3 \\
Steam boiler (ELEC) & 5 \\
Steam boiler (H$_2$R) & 5 \\
Steam boiler (NG) & 5 \\
Water boiler (ELEC) & 10 \\
Water boiler (H$_2$R) & 10 \\
Water boiler (NG) & 10 \\
Wind & 5 \\
\bottomrule
\end{tabular}
\end{singlespace}
\end{table}

\subsection{Storage components}
\label{sec:si-storage-tech}

This section documents the investment costs, lifetimes, and maintenance costs (Section~\ref{sec:si-stor-economics}), the technical parameters (Section~\ref{sec:si-stor-parameters}), and the installable capacities (Section~\ref{sec:si-stor-capacities}) of the storage components. In the first investment year 2025, only the existing capacities are available. All further investment options are available from 2030 onwards. The tables therefore report the investment years 2030 to 2045.

\subsubsection{Investment costs, lifetimes, and maintenance costs}
\label{sec:si-stor-economics}

Table~\ref{tab:si-invest-stor} reports the specific investment costs for all installable storage components. Investment costs are annualized and discounted as for the conversion components (Section~\ref{sec:si-conv-economics}). Table~\ref{tab:si-lifetime-stor} reports the lifetimes of the storage components. Table~\ref{tab:si-maintenance-stor} reports the yearly maintenance costs of the storage components.

\begin{table}[H]
\begin{singlespace}
\centering
\sitablefont
\setlength{\tabcolsep}{4pt}
\caption{Specific investment costs (in kEUR/MWh) for storage components by investment year. If only one value is reported, this value applies to all investment years.}
\label{tab:si-invest-stor}
\begin{tabularx}{\textwidth}{Yrrrrc}
\toprule
Component & 2030 & 2035 & 2040 & 2045 & Reference \\
\midrule
Cold water storage & \multicolumn{4}{c}{55} & \tmark{a} \\
Hot water storage & \multicolumn{4}{c}{27.5} & \tmark{b} \\
Hydrogen tank & 48 & 38.5 & 29 & 25.5 & \tmark{c} \\
Lithium-ion battery & 280 & 280 & 300 & 310 & \tmark{c} \\
Vanadium redox flow battery & 550 & 460 & 390 & 355 & \tmark{c} \\
\bottomrule
\end{tabularx}
\par\smallskip
{\footnotesize\raggedright \tmark{a}~own assumption, \tmark{b}~\cite{Reinert2023}, \tmark{c}~\cite{DEA_storage}.\par}
\end{singlespace}
\end{table}

\begin{table}[H]
\begin{singlespace}
\centering
\sitablefont
\setlength{\tabcolsep}{4pt}
\caption{Lifetimes (in years) for storage components by investment year. If only one value is reported, this value applies to all investment years.}
\label{tab:si-lifetime-stor}
\begin{tabularx}{\textwidth}{Yrrrrc}
\toprule
Component & 2030 & 2035 & 2040 & 2045 & Reference \\
\midrule
Cold water storage & \multicolumn{4}{c}{20} & \tmark{a} \\
Hot water storage & \multicolumn{4}{c}{20} & \tmark{a} \\
Hydrogen tank & \multicolumn{4}{c}{30} & \tmark{b} \\
Lithium-ion battery & 20 & 25 & 30 & 30 & \tmark{b} \\
Vanadium redox flow battery & 30 & 30 & 35 & 40 & \tmark{b} \\
\bottomrule
\end{tabularx}
\par\smallskip
{\footnotesize\raggedright \tmark{a}~own assumption, \tmark{b}~\cite{DEA_storage}.\par}
\end{singlespace}
\end{table}

\clearpage

\begin{table}[H]
\begin{singlespace}
\centering
\sitablefont
\setlength{\tabcolsep}{4pt}
\caption{Yearly maintenance costs (in kEUR/MWh per year, or as a percentage of the specific investment cost) for storage components by investment year. If only one value is reported, this value applies to all investment years.}
\label{tab:si-maintenance-stor}
\begin{tabularx}{\textwidth}{Yrrrrc}
\toprule
Component & 2030 & 2035 & 2040 & 2045 & Reference \\
\midrule
Cold water storage & \multicolumn{4}{c}{1.5\% of invest} & \tmark{a} \\
Hot water storage & \multicolumn{4}{c}{1.0\% of invest} & \tmark{b} \\
Hydrogen tank & \multicolumn{4}{c}{2.0\% of invest} & \tmark{a} \\
Lithium-ion battery & 8.09 & 8.14 & 8.68 & 8.95 & \tmark{c} \\
Vanadium redox flow battery & 3.39 & 2.88 & 2.43 & 2.21 & \tmark{c} \\
\bottomrule
\end{tabularx}
\par\smallskip
{\footnotesize\raggedright \tmark{a}~own assumption, \tmark{b}~\cite{Reinert2023}, \tmark{c}~\cite{DEA_storage}.\par}
\end{singlespace}
\end{table}

\subsubsection{Technical parameters}
\label{sec:si-stor-parameters}

Table~\ref{tab:si-storage} reports the technical storage parameters per investment year: charge and discharge efficiencies, self-discharge, and the power-to-energy (P/E) ratio. Self-discharge is implemented as an hourly loss factor on the state of charge. The table reports the equivalent loss per day.

\begingroup
\singlespacing\sitablefont
\setlength{\tabcolsep}{2.5pt}
\begin{longtable}{L{2.3cm}L{4.3cm}rrrr}
\caption{Technical storage parameters by investment year. Self-discharge is the equivalent daily loss implied by the hourly loss factor of the model. If only one value is reported, this value applies to all investment years. P/E: power-to-energy.}
\label{tab:si-storage}\\
\toprule
Component & Parameter & 2030 & 2035 & 2040 & 2045 \\
\midrule
\endfirsthead
\multicolumn{6}{l}{\itshape Table~\ref{tab:si-storage} (continued)}\\
\toprule
Component & Parameter & 2030 & 2035 & 2040 & 2045 \\
\midrule
\endhead
\midrule
\multicolumn{6}{r}{\itshape continued on next page}\\
\endfoot
\bottomrule
\endlastfoot
\multirow{4}{*}{\parbox{2.3cm}{Cold water\\storage}} & Charge efficiency\tmark{a} & \multicolumn{4}{c}{0.90} \\*
 & Discharge efficiency\tmark{a} & \multicolumn{4}{c}{0.95} \\*
 & Self-discharge (\%/day)\tmark{a} & \multicolumn{4}{c}{1.20} \\*
 & P/E ratio (MW/MWh)\tmark{a} & \multicolumn{4}{c}{0.5} \\
\addlinespace
\multirow{4}{*}{\parbox{2.3cm}{Hot water\\storage}} & Charge efficiency\tmark{b} & \multicolumn{4}{c}{0.95} \\*
 & Discharge efficiency\tmark{b} & \multicolumn{4}{c}{0.95} \\*
 & Self-discharge (\%/day)\tmark{a} & \multicolumn{4}{c}{0.48} \\*
 & P/E ratio (MW/MWh)\tmark{a} & \multicolumn{4}{c}{0.5} \\
\addlinespace
\multirow{4}{*}{\parbox{2.3cm}{Hydrogen\\tank}} & Charge efficiency\tmark{c} & 0.89 & 0.895 & 0.90 & 0.90 \\*
 & Discharge efficiency\tmark{c} & \multicolumn{4}{c}{1.00} \\*
 & Self-discharge (\%/day)\tmark{c} & \multicolumn{4}{c}{0.14} \\*
 & P/E ratio (MW/MWh)\tmark{a} & \multicolumn{4}{c}{0.2} \\
\addlinespace
\multirow{4}{*}{\parbox{2.3cm}{Lithium-ion\\battery}} & Charge efficiency\tmark{c} & 0.98 & 0.985 & 0.985 & 0.985 \\*
 & Discharge efficiency\tmark{c} & 0.97 & 0.975 & 0.975 & 0.975 \\*
 & Self-discharge (\%/day)\tmark{c} & \multicolumn{4}{c}{0.10} \\*
 & P/E ratio (MW/MWh)\tmark{c} & \multicolumn{4}{c}{1.0} \\
\addlinespace
\multirow{4}{*}{\parbox{2.3cm}{Vanadium\\redox flow\\battery}} & Charge efficiency\tmark{c} & 0.894 & 0.894 & 0.911 & 0.917 \\*
 & Discharge efficiency\tmark{c} & 0.894 & 0.894 & 0.911 & 0.917 \\*
 & Self-discharge (\%/day)\tmark{c} & 0.20 & 0.20 & 0.15 & 0.13 \\*
 & P/E ratio (MW/MWh)\tmark{c} & \multicolumn{4}{c}{0.83} \\
\end{longtable}
\vspace{-6pt}
{\footnotesize\raggedright \tmark{a}~own assumption, \tmark{b}~\cite{Reinert2023}, \tmark{c}~\cite{DEA_storage}.\par}
\endgroup

\subsubsection{Installable capacities}
\label{sec:si-stor-capacities}

No storage capacity exists at the start of the optimization horizon. Table~\ref{tab:si-capacity-stor} reports the assumed maximum installable capacity of each storage component.

\begin{table}[H]
\begin{singlespace}
\centering
\sitablefont
\setlength{\tabcolsep}{4pt}
\caption{Assumed maximum installable capacity (in MWh) of the storage components.}
\label{tab:si-capacity-stor}
\begin{tabular}{lR{3.2cm}}
\toprule
Component & Maximum installable capacity \\
\midrule
Cold water storage & 50 \\
Hot water storage & 25 \\
Hydrogen tank & 100 \\
Lithium-ion battery & 20 \\
Vanadium redox flow battery & 10 \\
\bottomrule
\end{tabular}
\end{singlespace}
\end{table}

\subsection{Energy purchases and price scenarios}
\label{sec:si-procurement}

This section documents the grid connections (Section~\ref{sec:si-grid-connections}), the power purchase agreements (Section~\ref{sec:si-ppa}), and the price scenarios (Section~\ref{sec:si-price-scenarios}).

\subsubsection{Grid connections}
\label{sec:si-grid-connections}

The grid connections carry the energy purchased at the spot market and, for electricity, also the energy contracted under the PPAs of Section~\ref{sec:si-ppa}. The grid connection capacity is contracted per investment year for the investment period and bounds the hourly grid usage. Yearly capacity charges apply to the contracted grid connection capacity. The capacity charges are 113~kEUR/MW for electricity \cite{HirthEicke2023} and 13.28~kEUR/MW for natural gas \cite{Westnetz2024}, both taken from \cite{Schricker2026b}, and 44.34~kEUR/MW for hydrogen. The hydrogen capacity charge is derived from the natural gas capacity charge under the assumption that pipeline capacity charges are proportional to the volume flow, and is therefore scaled by the ratio of the volumetric lower heating values of natural gas and hydrogen. The capacity charges are constant across investment years. Table~\ref{tab:si-grid-limits} reports the assumed grid connection limits up to which grid connection capacity can be installed.

\begin{table}[H]
\begin{singlespace}
\centering
\sitablefont
\setlength{\tabcolsep}{4pt}
\caption{Assumed grid connection limits (in MW) per investment year. The hydrogen grid connection is assumed to become available in 2035.}
\label{tab:si-grid-limits}
\begin{tabularx}{\textwidth}{Yrrrrr}
\toprule
Grid connection & 2025 & 2030 & 2035 & 2040 & 2045 \\
\midrule
Electricity grid & 15 & 30 & 45 & 60 & 75 \\
Natural gas grid & 30 & 30 & 30 & 30 & 30 \\
Hydrogen grid & 0 & 0 & 30 & 50 & 70 \\
\bottomrule
\end{tabularx}
\end{singlespace}
\end{table}

\subsubsection{Power purchase agreements (PPAs)}
\label{sec:si-ppa}

PPAs are modeled as in \cite{Schricker2026b} as ten-year pay-as-produced contracts for photovoltaics, onshore wind, and offshore wind. Contracted electricity that cannot be used at the site is resold at the spot market with the selling price reduction assumptions of \cite{Schricker2026b} and is not delivered physically. Electricity generated or stored on site cannot be sold. Table~\ref{tab:si-ppa} reports the maximum contract sizes and the contract prices. The 2025 prices are the market prices of \cite{SPGlobal2024} used in \cite{Schricker2026b}. The 2035 and 2045 prices are the mean levelized cost of electricity that \cite{FraunhoferISE2024} projects for the respective technology, multiplied by a constant factor. The factor is chosen per technology such that the same calculation reproduces the 2025 price. The 2030 and 2040 prices are linearly interpolated. The contract prices are assumed to be identical across scenarios, because the scenarios differ in fossil fuel and CO$_2$ certificate prices, which drive the spot market price but not the cost of the renewable generation on which the contract prices are based.

\begin{table}[H]
\begin{singlespace}
\centering
\sitablefont
\setlength{\tabcolsep}{4pt}
\caption{PPA contract prices (in EUR/MWh) by investment year and maximum contracted capacity.}
\label{tab:si-ppa}
\begin{tabularx}{\textwidth}{Ycrrrrrc}
\toprule
PPA type & Max.\ cap.\ (MW) & 2025 & 2030 & 2035 & 2040 & 2045 & Reference \\
\midrule
Photovoltaics & 2.0 & 91.0 & 81.3 & 71.6 & 69.0 & 66.4 & \tmark{a} \\
Wind onshore & 3.0 & 89.0 & 85.7 & 82.3 & 79.8 & 77.3 & \tmark{a} \\
Wind offshore & 5.0 & 97.7 & 95.7 & 93.7 & 91.7 & 89.7 & \tmark{a} \\
\bottomrule
\end{tabularx}
\par\smallskip
{\footnotesize\raggedright \tmark{a}~\cite{SPGlobal2024, FraunhoferISE2024}.\par}
\end{singlespace}
\end{table}

\subsubsection{Price scenarios}
\label{sec:si-price-scenarios}

The seven price scenarios and their 2025 electricity, natural gas, and CO$_2$ certificate prices are taken from \cite{Schricker2026b}. The present study adds hydrogen prices and extends the annual average prices of each scenario until 2045. Table~\ref{tab:si-energy-prices} reports the annual average CO$_2$ certificate and energy carrier prices of all seven scenarios per investment year. The sources of Scenarios~1--4 report no hydrogen prices, so these scenarios share the hydrogen price trajectory of \cite{Agora2021}. Scenarios~5--7 use the hydrogen prices of \cite{Gorbach2022}. The industrial purchase prices include the supplier margins and fees assumed in \cite{Schricker2026b}, but not the capacity charges of Section~\ref{sec:si-grid-connections}.

CO$_2$ certificate costs are charged at the natural gas grid connection with the constant emission factor of 200.8~kg~CO$_2$/MWh of the German national emission reporting \citep{UBA2022}. SMR with CCS receives a certificate credit for the captured share of these emissions (Section~\ref{sec:si-conv-coefficients}).

\begingroup
\singlespacing\sitablefont
\setlength{\tabcolsep}{4pt}
\begin{longtable}{lcrrrrrc}
\caption{Annual average CO$_2$ certificate prices (in EUR/tCO$_2$) and energy carrier prices (in EUR/MWh) by scenario and investment year. Hydrogen is purchasable from 2035, when the hydrogen grid connection becomes available (Table~\ref{tab:si-grid-limits}).}
\label{tab:si-energy-prices}\\
\toprule
Commodity & Scenario & 2025 & 2030 & 2035 & 2040 & 2045 & Reference \\
\midrule
\endfirsthead
\multicolumn{8}{l}{\itshape Table~\ref{tab:si-energy-prices} (continued)}\\
\toprule
Commodity & Scenario & 2025 & 2030 & 2035 & 2040 & 2045 & Reference \\
\midrule
\endhead
\midrule
\multicolumn{8}{r}{\itshape continued on next page}\\
\endfoot
\bottomrule
\multicolumn{8}{@{}p{\textwidth}@{}}{\footnotesize\raggedright
\tmark{a}~\cite{Prognos2024},
\tmark{b}~\cite{TUM_thesis2020},
\tmark{c}~\cite{Gorbach2022},
\tmark{d}~\cite{Bundesnetzagentur2024},
\tmark{e}~\cite{Agora2021}.}\\
\endlastfoot
\multirow{7}{*}{CO$_2$ certificate}
 & 1 & 76 & 137 & 218.5 & 300 & 300 & \tmark{a} \\*
 & 2 & 72 & 118 & 141 & 164 & 187 & \tmark{a} \\*
 & 3 & 62 & 64 & 83.5 & 103 & 122.5 & \tmark{a} \\*
 & 4 & 31 & 41.8 & 52.7 & 63.5 & 74.4 & \tmark{b} \\*
 & 5 & 45 & 72.5 & 86.3 & 100 & 113.8 & \tmark{c} \\*
 & 6 & 201 & 222 & 243 & 264 & 285 & \tmark{c} \\*
 & 7 & 100 & 400 & 500 & 600 & 700 & \tmark{c} \\
\addlinespace
\multirow{7}{*}{Electricity}
 & 1 & 358.0 & 207.7 & 202.5 & 225.6 & 189.8 & \tmark{a,d} \\*
 & 2 & 187.6 & 170.4 & 168.7 & 172.7 & 155.7 & \tmark{a,d} \\*
 & 3 & 172.5 & 133.8 & 135.4 & 140.0 & 131.1 & \tmark{a,d} \\*
 & 4 & 142.2 & 147.2 & 153.4 & 157.8 & 158.1 & \tmark{b} \\*
 & 5 & 191.0 & 191.0 & 191.0 & 191.0 & 191.0 & \tmark{c} \\*
 & 6 & 241.0 & 241.0 & 241.0 & 241.0 & 241.0 & \tmark{c} \\*
 & 7 & 191.0 & 191.0 & 191.0 & 191.0 & 191.0 & \tmark{c} \\
\addlinespace
\multirow{7}{*}{Hydrogen}
 & 1 & -- & -- & 138.0 & 131.0 & 124.0 & \tmark{e} \\*
 & 2 & -- & -- & 138.0 & 131.0 & 124.0 & \tmark{e} \\*
 & 3 & -- & -- & 138.0 & 131.0 & 124.0 & \tmark{e} \\*
 & 4 & -- & -- & 138.0 & 131.0 & 124.0 & \tmark{e} \\*
 & 5 & -- & -- & 188.5 & 164.0 & 139.5 & \tmark{c} \\*
 & 6 & -- & -- & 155.0 & 120.0 &  85.0 & \tmark{c} \\*
 & 7 & -- & -- & 177.3 & 149.3 & 121.3 & \tmark{c} \\
\addlinespace
\multirow{7}{*}{Natural gas}
 & 1 & 240.1 &  66.1 &  50.4 &  51.3 &  52.2 & \tmark{a,d} \\*
 & 2 &  66.1 &  49.1 &  48.1 &  47.2 &  46.2 & \tmark{a,d} \\*
 & 3 &  58.3 &  38.7 &  39.0 &  39.3 &  39.6 & \tmark{a,d} \\*
 & 4 &  53.3 &  54.9 &  56.8 &  57.0 &  57.0 & \tmark{b} \\*
 & 5 &  48.6 &  48.6 &  48.6 &  48.6 &  48.6 & \tmark{c} \\*
 & 6 & 100.0 & 100.0 & 100.0 & 100.0 & 100.0 & \tmark{c} \\*
 & 7 &  97.2 &  97.2 &  97.2 &  97.2 &  97.2 & \tmark{c} \\
\end{longtable}
\endgroup

\begin{figure}[H]
\centering
\includegraphics[width=\textwidth]{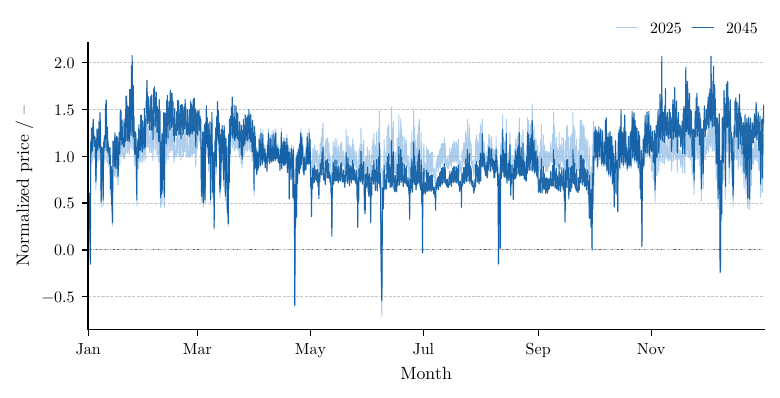}
\caption{Normalized hourly electricity price profile in 2025 and 2045.}
\label{fig:si-electricity-price}
\end{figure}

Electricity and natural gas have intra-annual price profiles. Hydrogen and CO$_2$ certificate prices are assumed to be constant within each investment year. The hourly spot market price is the annual average price of the scenario multiplied by a normalized hourly profile with an annual mean of one, taken from the price time series of \cite{Reinert2023}. The natural gas profile is applied identically in every investment year. The electricity profile of each investment year keeps the hourly pattern within each month and rescales the monthly means to the monthly seasonality that \cite{Prognos2024} projects for that investment year. The electricity profile therefore becomes more volatile toward the later investment years. Figure~\ref{fig:si-electricity-price} shows the normalized hourly electricity price profile in 2025 and 2045.

\newpage
\section{Computational Setup and Supplementary Results}
\label{sec:si-results}

This section documents the computational setup and presents supplementary results of our study. Section~\ref{sec:si-convergence} describes the computational setup of the exploration, Section~\ref{sec:si-range-det} documents the installed capacities of the perfect-foresight pathways, Section~\ref{sec:si-layers} the low-regret pathway space at higher regret tolerances, Section~\ref{sec:si-tech-ranges} the capacity ranges and utilization of the conversion and storage components across the low-regret pathway space, and Section~\ref{sec:si-regret-scenarios} the scenario-wise regret of the explored pathways.

\subsection{Computational setup}
\label{sec:si-convergence}

The implementation uses the open-source model SecMOD \citep{Reinert2022} and Gurobi~12.0.0. Following \citet{Schricker2026b}, the full hourly time series of every investment year and scenario is aggregated into four representative typical days by hierarchical clustering with the duration representation which conserves the duration curves of the original time series \citep{Kotzur2018,Hoffmann2022}. Storage components are modeled with a maximum storage duration of 24 hours. A feasibility scenario containing the time series extrema is appended to the scenario set to enable operational feasibility during events with peak demand and low availability as in \citet{Schricker2026b}.

The exploration allocates the optimization runs in batches. Each regret tolerance is explored with 96 parallel workers in batches of 384 search directions, at a mean runtime of about 33 minutes per optimization. Each exploration terminates at the 96-hour wall-clock limit, after 16{,}268 directional optimizations at the regret tolerances of 2.5\%, 3.5\%, and 5\% and after 15{,}884 at 10\%. All optimization problems of Eq.~(8) of the main text are linear programs solved at the default Gurobi feasibility and optimality tolerances of $10^{-6}$.

\begin{figure}[H]
\centering
\includegraphics[width=\textwidth]{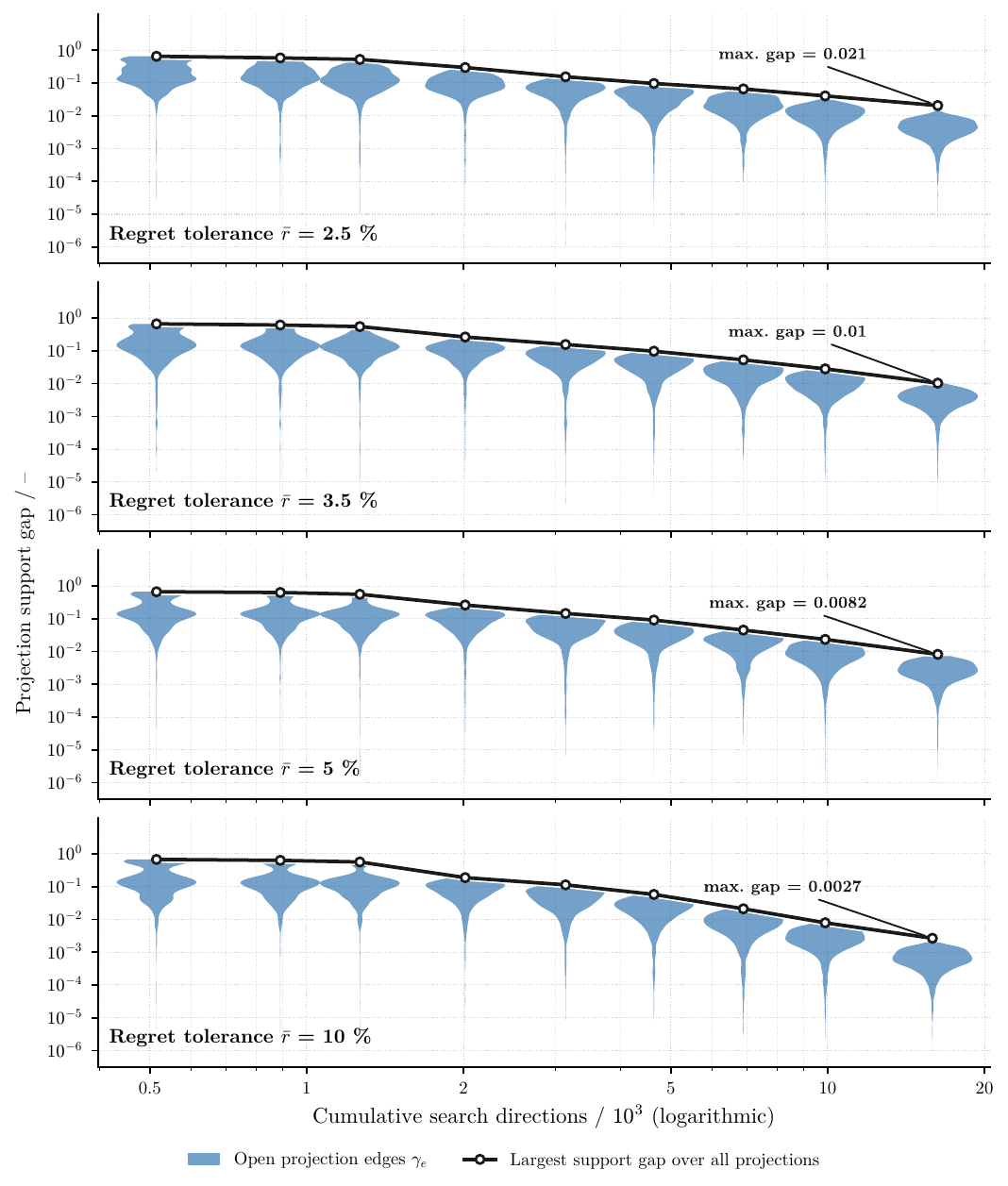}
\caption{Support gaps $\gamma_e$ of the projection edges against the cumulative number of explored search directions, at four regret tolerances $\bar r$. The horizontal axis is logarithmic. Each violin is the distribution of the support gaps over a random sample of 2{,}500 projection edges. The black line is the largest support gap over all projection edges.}
\label{fig:si-gap}
\end{figure}

Figure~\ref{fig:si-gap} shows the support gaps of the projection edges (Section~2.3 of the main text) against the cumulative number of explored search directions. Since every exploration variable is rescaled to $[0,1]$, a support gap is a fraction of a capacity range. The largest support gap over all projections falls to 0.021 at $\bar r=2.5\%$, 0.010 at 3.5\%, 0.008 at 5\%, and 0.003 at 10\%. The boundary of every two-dimensional projection of the low-regret pathway space therefore lies at most about 2\% of a capacity range outside the explored pathway set.

\subsection{Perfect-foresight pathways}
\label{sec:si-range-det}

\begin{figure}[H]
\centering
\includegraphics[width=\textwidth]{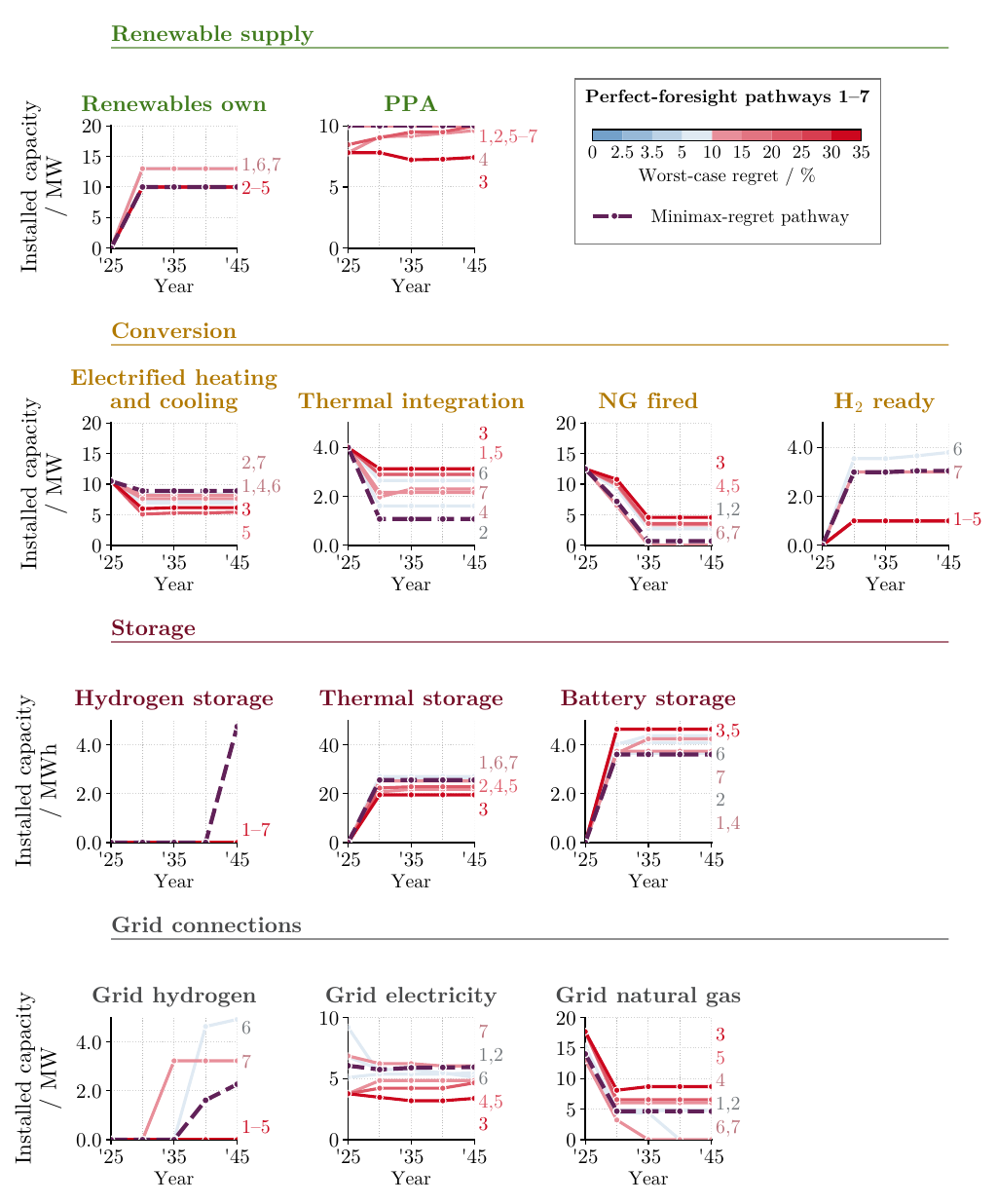}
\caption{Installed capacities of the seven perfect-foresight pathways,
together with the minimax-regret pathway (dashed purple line). Each pathway is
colored by its worst-case regret and annotated with its
scenario identifier at the final year. Capacities change only in the investment years. The pathway lines are connected between investment years for readability and do not imply a continuous capacity expansion. Electrolysis and SMR are omitted because no pathway installs them.}
\label{fig:si-range-det}
\end{figure}

Figure~\ref{fig:si-range-det} shows the installed capacities of the seven perfect-foresight pathways together with the minimax-regret pathway. In the pathways optimized for Scenarios~1--5, natural gas fired steam boilers reach their maximum installable capacity (Table~\ref{tab:si-capacity-conv}), so additional steam boiler capacity is provided by hydrogen-ready boilers even though they are never operated with hydrogen.

\newpage
\subsection{Low-regret pathway space at higher regret tolerances}
\label{sec:si-layers}

\begin{figure}[H]
\centering
\includegraphics[width=\textwidth]{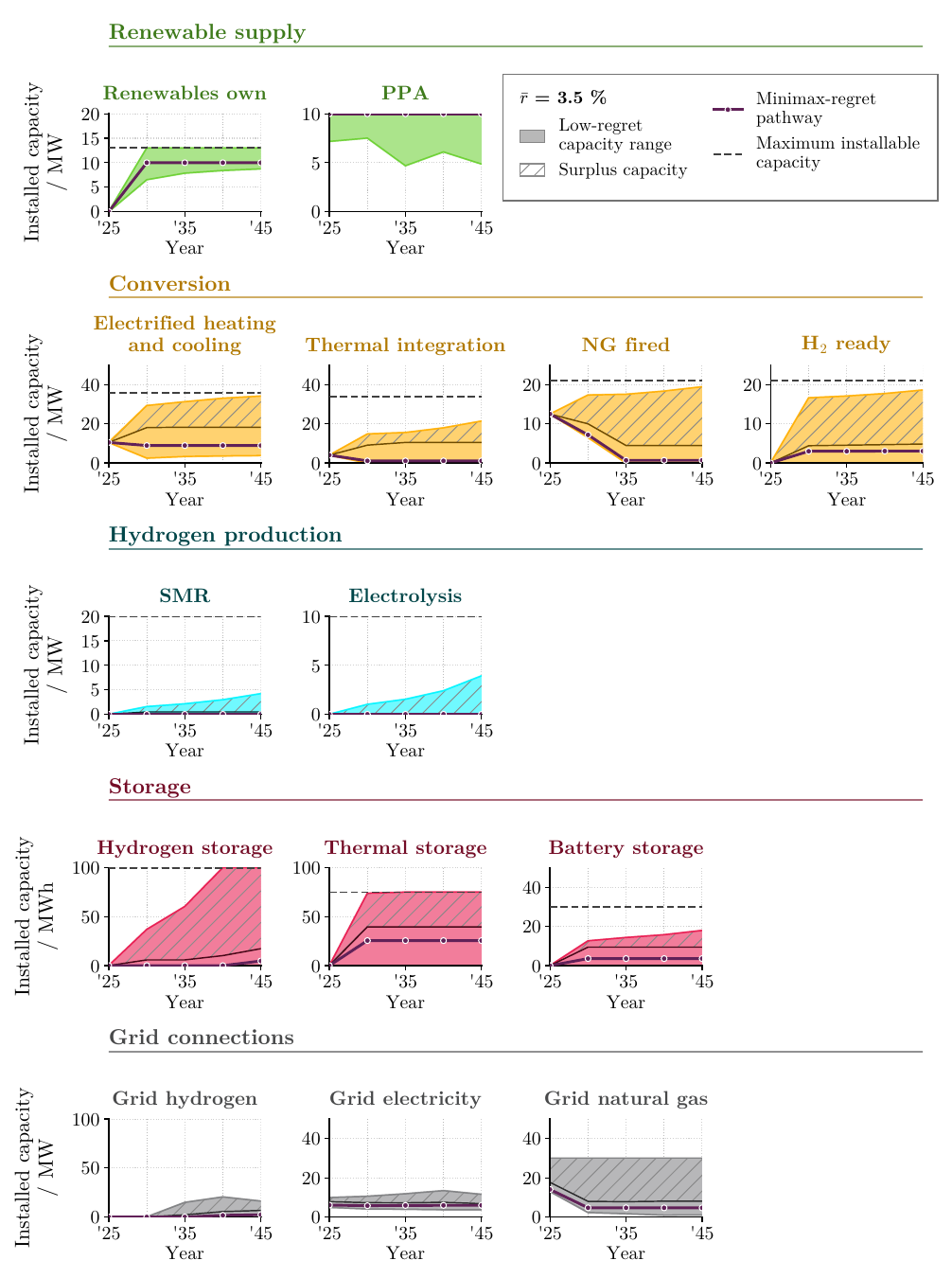}
\caption{Low-regret capacity ranges at $\bar r=3.5\%$. The hatched part of each band is the surplus capacity. The purple line is the minimax-regret pathway, the dashed grey line the maximum installable capacity of a technology family.}
\label{fig:si-layers-3-5}
\end{figure}

\begin{figure}[H]
\centering
\includegraphics[width=\textwidth]{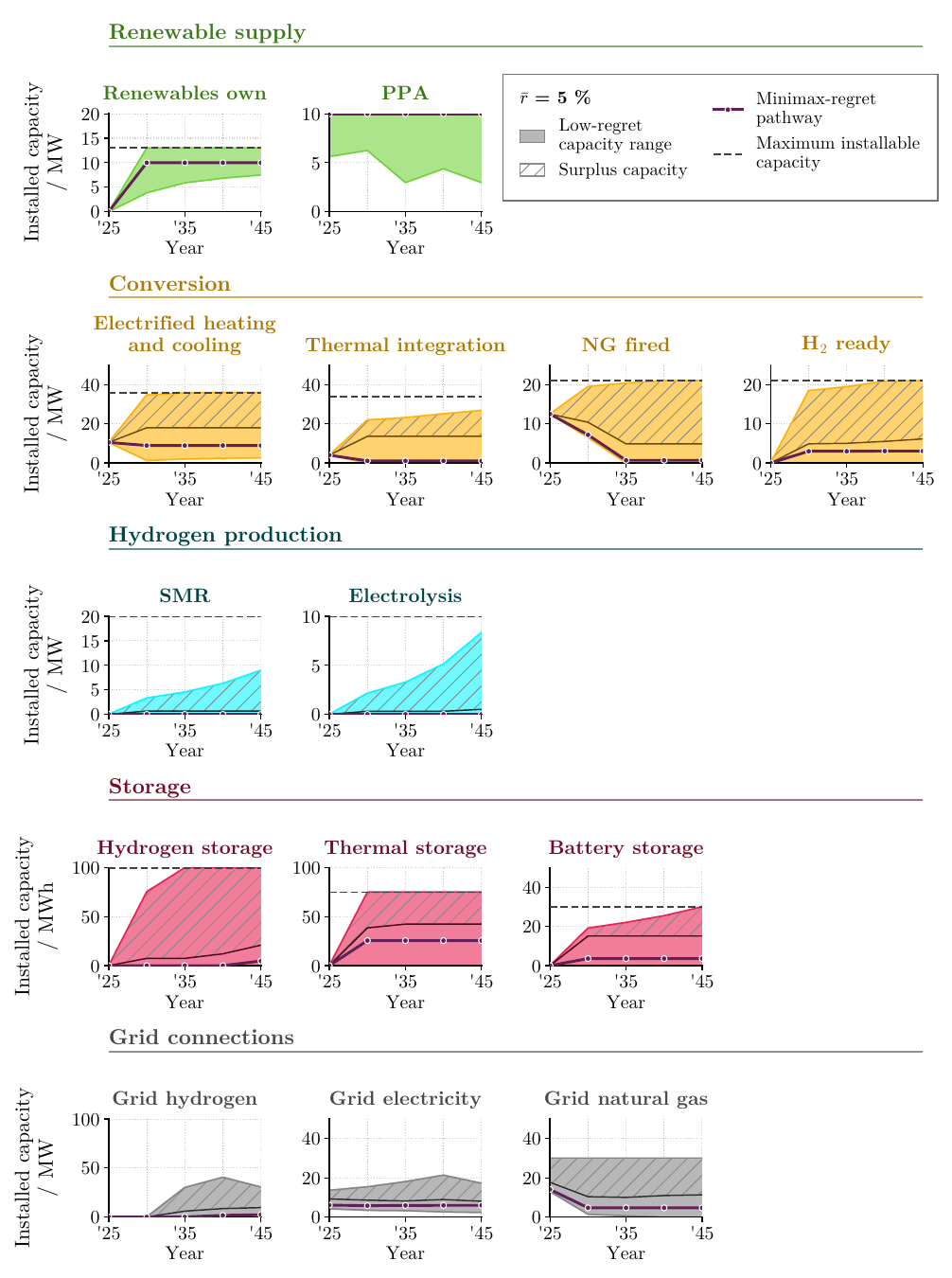}
\caption{Low-regret capacity ranges at $\bar r=5\%$. The hatched part of each band is the surplus capacity. The purple line is the minimax-regret pathway, the dashed grey line the maximum installable capacity of a technology family.}
\label{fig:si-layers-5}
\end{figure}

\begin{figure}[H]
\centering
\includegraphics[width=\textwidth]{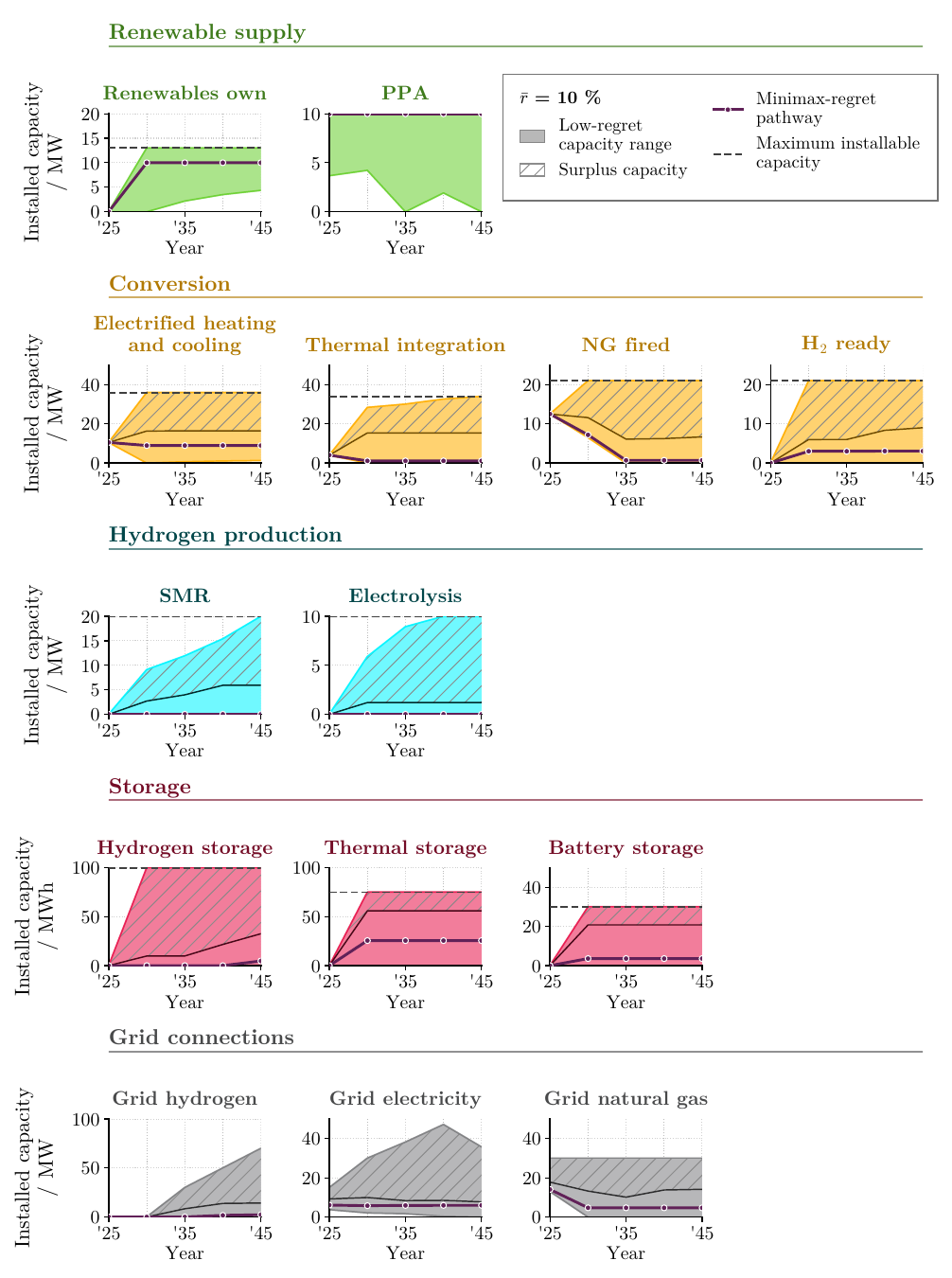}
\caption{Low-regret capacity ranges at $\bar r=10\%$. The hatched part of each band is the surplus capacity. The purple line is the minimax-regret pathway, the dashed grey line the maximum installable capacity of a technology family.}
\label{fig:si-layers-10}
\end{figure}

\begin{figure}[H]
\centering
\includegraphics[width=\textwidth]{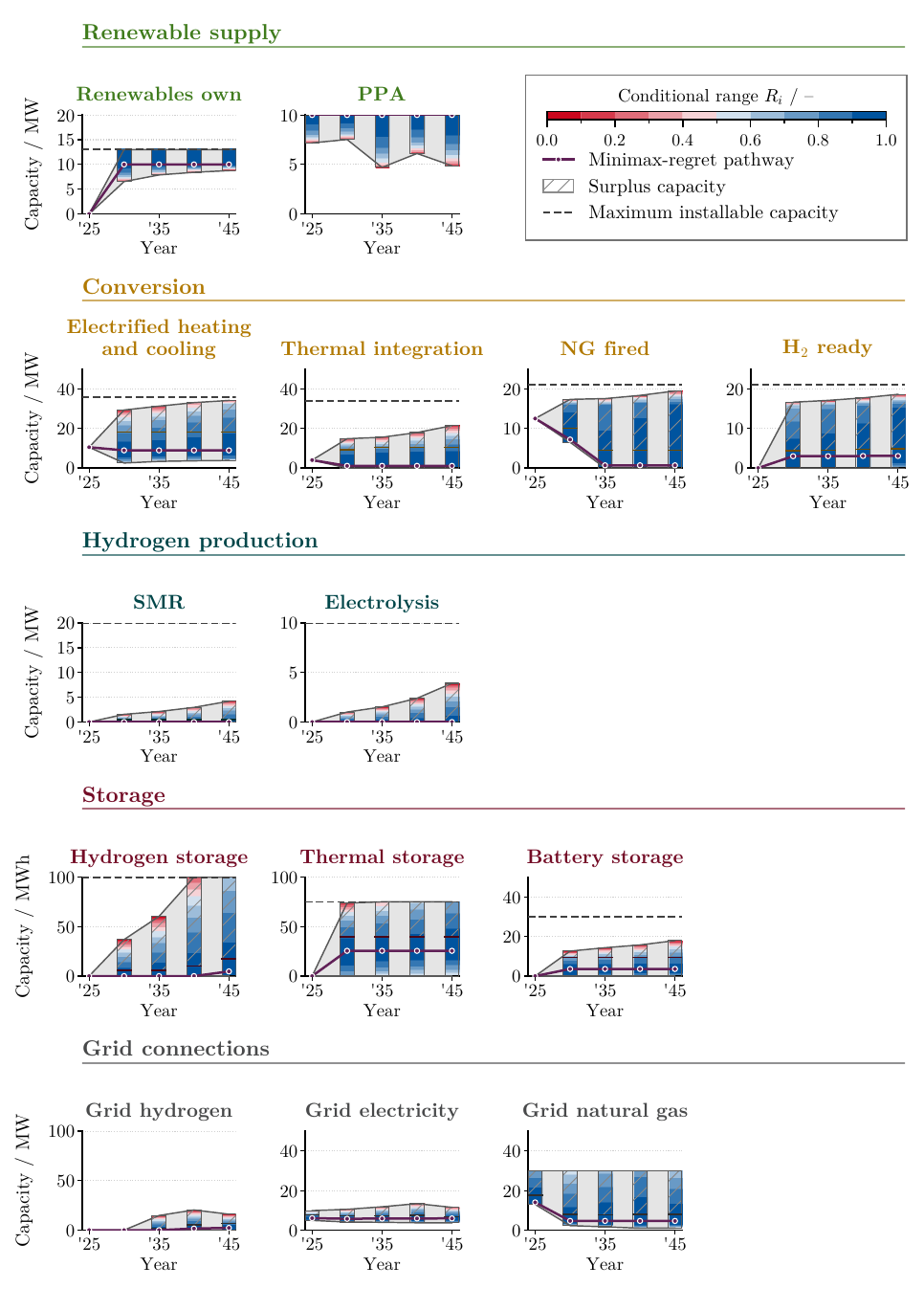}
\caption{Low-regret capacity ranges at $\bar r=3.5\%$, colored by the conditional range $R_i$ at the five investment years. The hatched part of each band is the surplus capacity. The purple line is the minimax-regret pathway, the dashed grey line the maximum installable capacity of a technology family.}
\label{fig:si-logdensity-3-5}
\end{figure}

\begin{figure}[H]
\centering
\includegraphics[width=\textwidth]{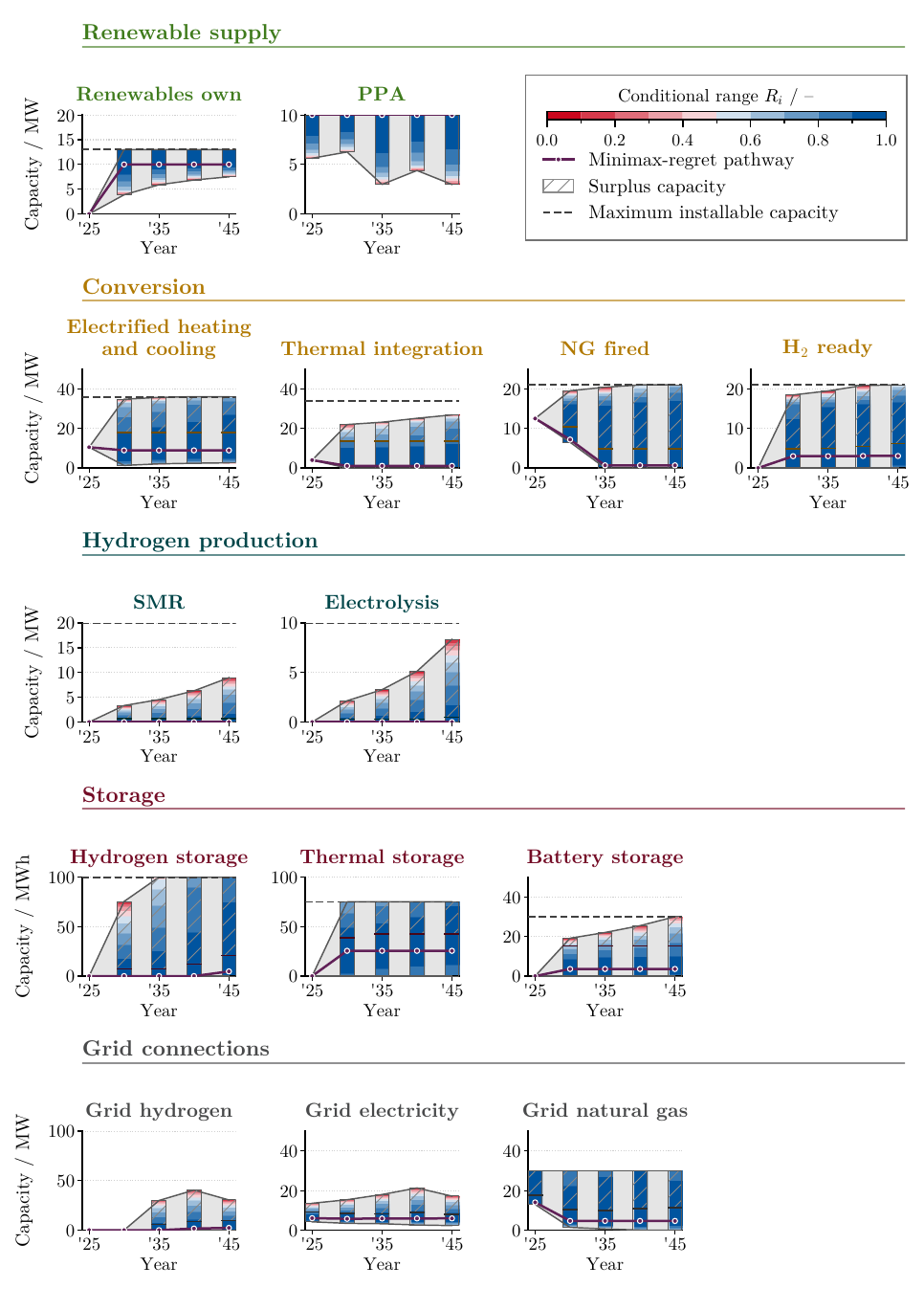}
\caption{Low-regret capacity ranges at $\bar r=5\%$, colored by the conditional range $R_i$ at the five investment years. The hatched part of each band is the surplus capacity. The purple line is the minimax-regret pathway, the dashed grey line the maximum installable capacity of a technology family.}
\label{fig:si-logdensity-5}
\end{figure}

\begin{figure}[H]
\centering
\includegraphics[width=\textwidth]{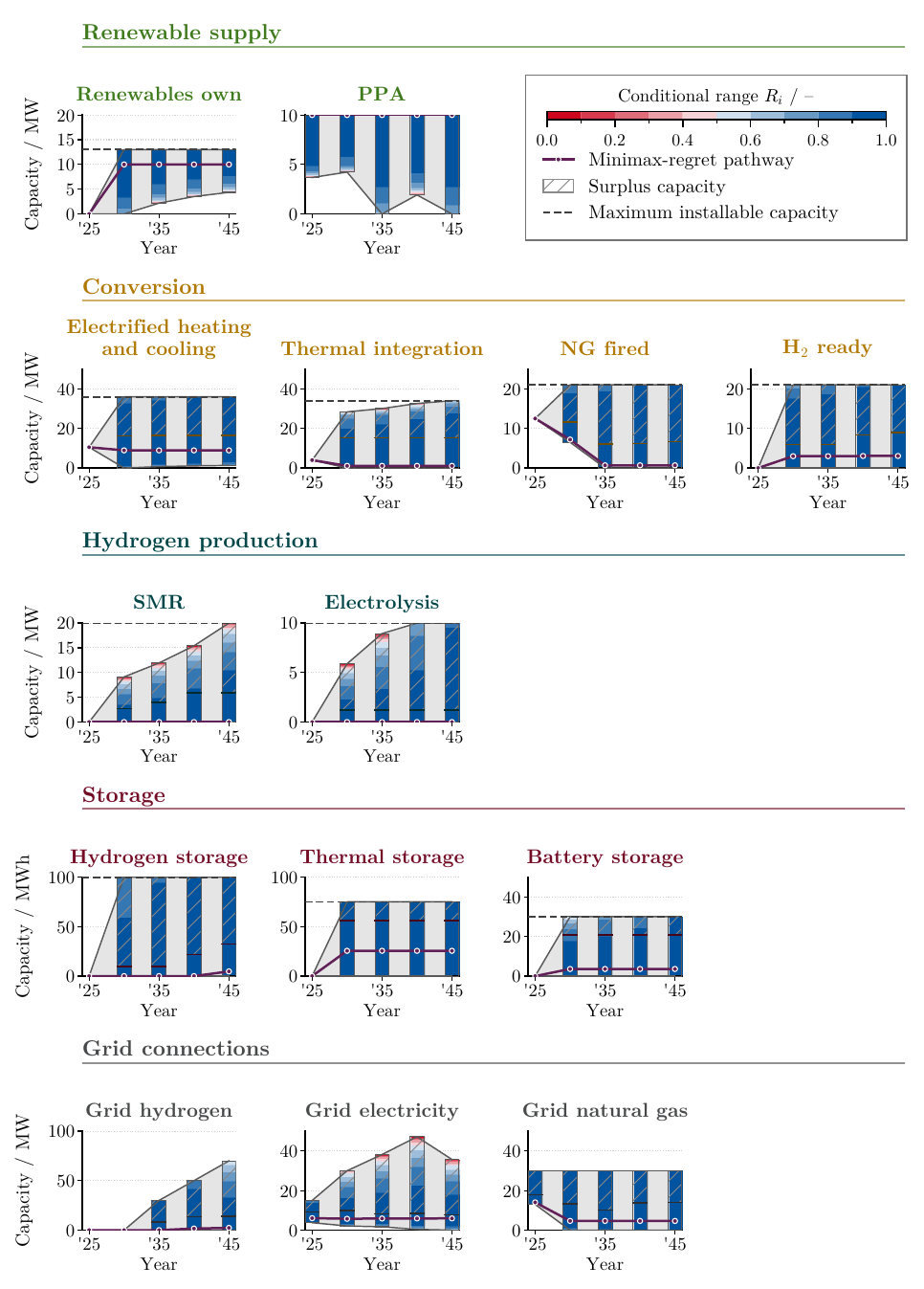}
\caption{Low-regret capacity ranges at $\bar r=10\%$, colored by the conditional range $R_i$ at the five investment years. The hatched part of each band is the surplus capacity. The purple line is the minimax-regret pathway, the dashed grey line the maximum installable capacity of a technology family.}
\label{fig:si-logdensity-10}
\end{figure}

\begin{figure}[H]
\centering
\includegraphics[width=\textwidth]{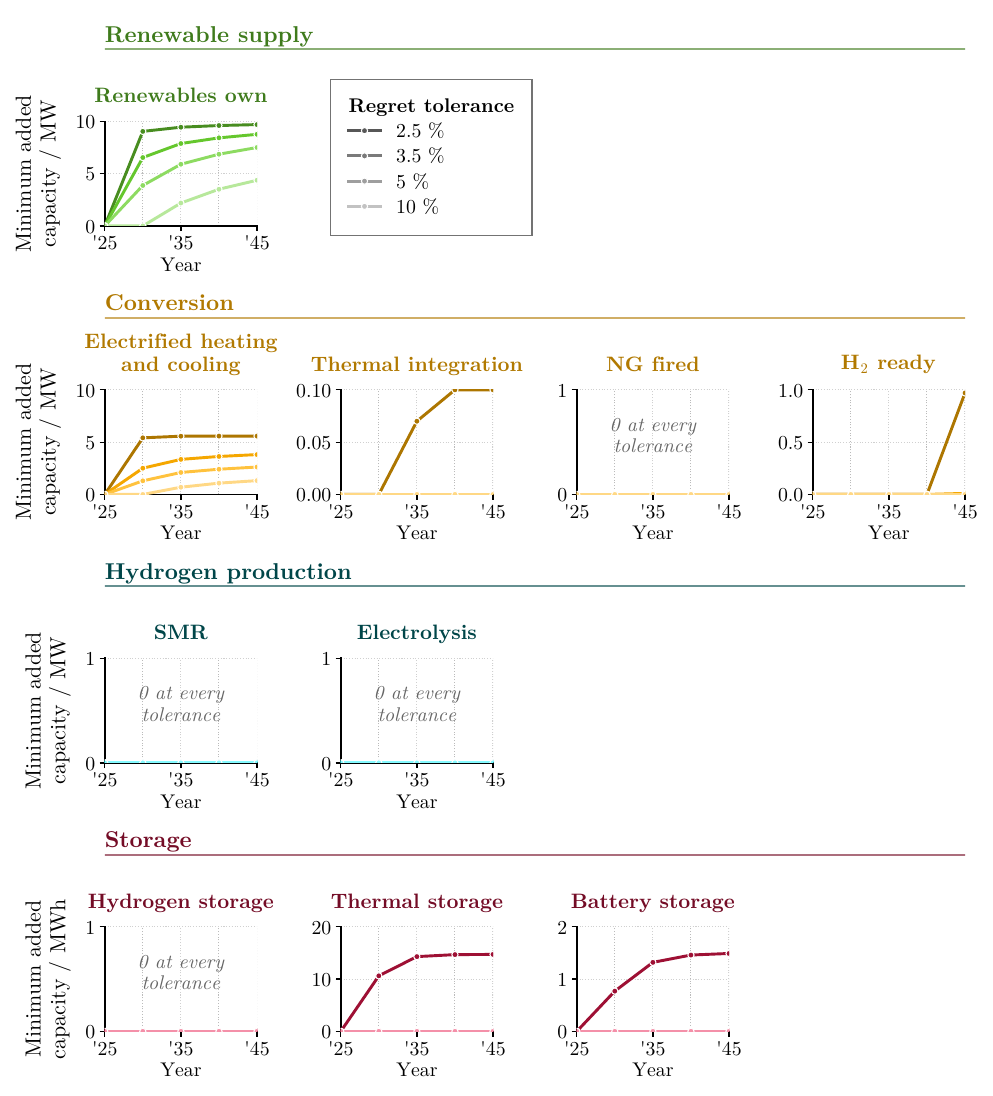}
\caption{Minimum cumulative added capacity per technology family and
investment year over all explored pathways, at the four regret tolerances
$\bar r$.}
\label{fig:si-min-additions}
\end{figure}

Figures~\ref{fig:si-layers-3-5}, \ref{fig:si-layers-5}, and~\ref{fig:si-layers-10} show the 3.5\%, 5\%, and 10\% low-regret pathway spaces. The hatched part of each band is the surplus capacity at that tolerance, i.e., the capacity that no capacity-minimizing search reaches (Section~2.4 of the main text). The capacity-minimizing searches are the searches with non-positive direction coefficients. Searches whose weighted variables all end at their lower bound, i.e., at an installed capacity of zero or at the capacity of the existing components, are excluded from this set, because the search direction then no longer determines the other variables of the solution. Figures~\ref{fig:si-logdensity-3-5}, \ref{fig:si-logdensity-5}, and~\ref{fig:si-logdensity-10} show the low-regret capacity ranges at the same tolerances, colored by the conditional range $R_i$.
Figure~\ref{fig:si-min-additions} shows the minimum cumulative added capacity per technology family and investment year over all explored pathways, at the four regret tolerances $\bar r$. Cumulative added capacity cannot decrease along a pathway. Wherever the minimum is zero through an investment year, the explored pathway that reaches this minimum has consequently not installed capacities in any earlier year either.

\newpage

\subsection{Component capacity ranges and utilization}
\label{sec:si-tech-ranges}

\begin{figure}[H]
\centering
\includegraphics[width=\textwidth]{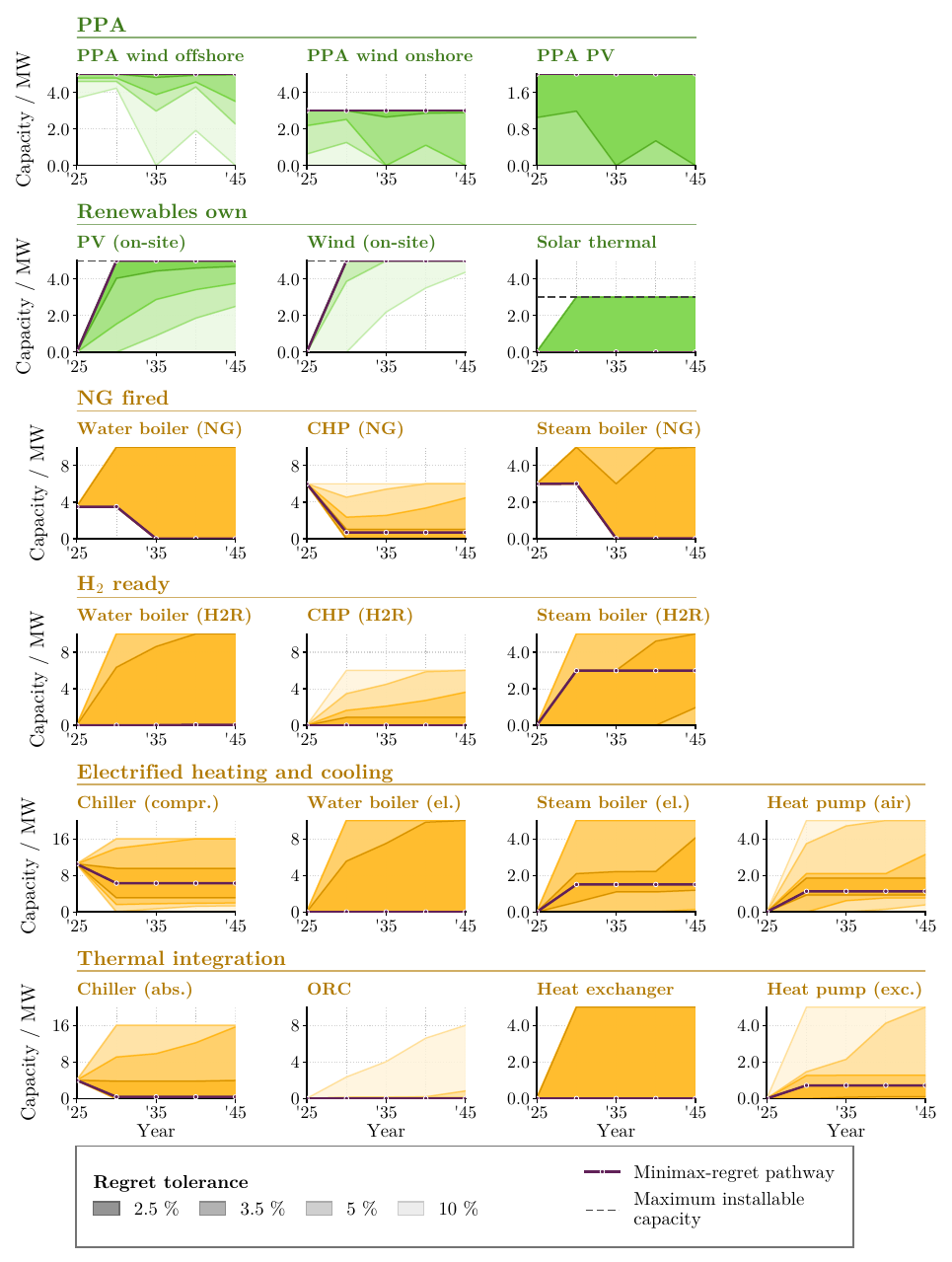}
\caption{Low-regret capacity ranges of the
conversion components. The dashed grey line marks the maximum installable capacity
of a component (Tables~\ref{tab:si-ppa}, \ref{tab:si-capacity-conv},
and~\ref{tab:si-capacity-stor}).}
\label{fig:si-tech-ranges-a}
\end{figure}

\begin{figure}[H]
\centering
\includegraphics[width=\textwidth]{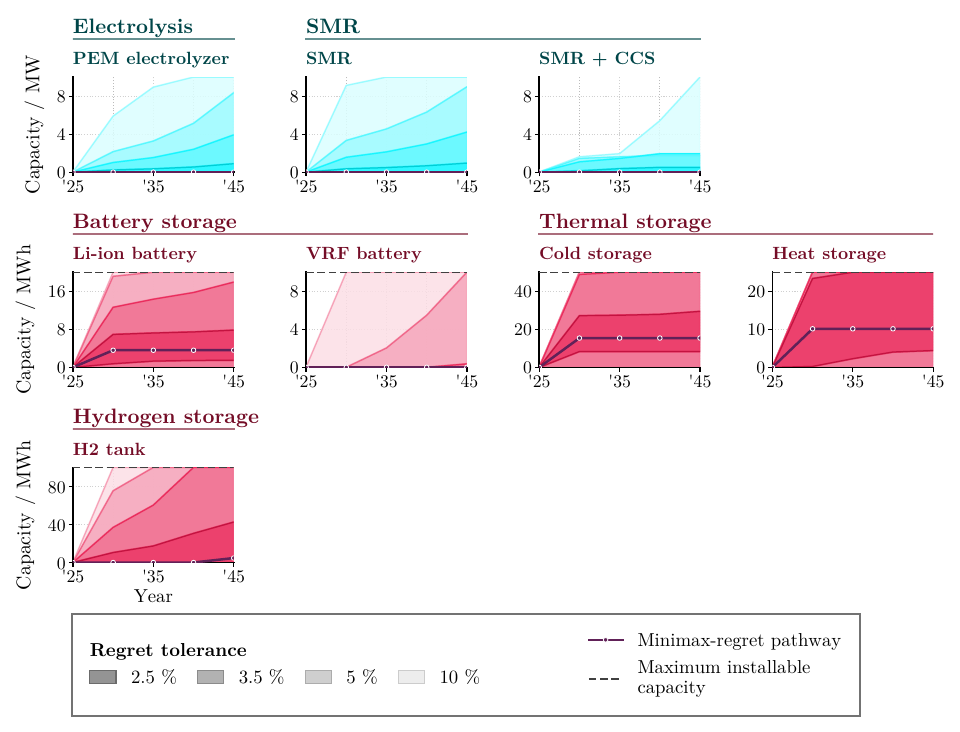}
\caption{Low-regret capacity ranges of the hydrogen production and
storage components. The dashed grey line marks the maximum installable capacity.}
\label{fig:si-tech-ranges-b}
\end{figure}

\begin{figure}[H]
\centering
\includegraphics[width=\textwidth]{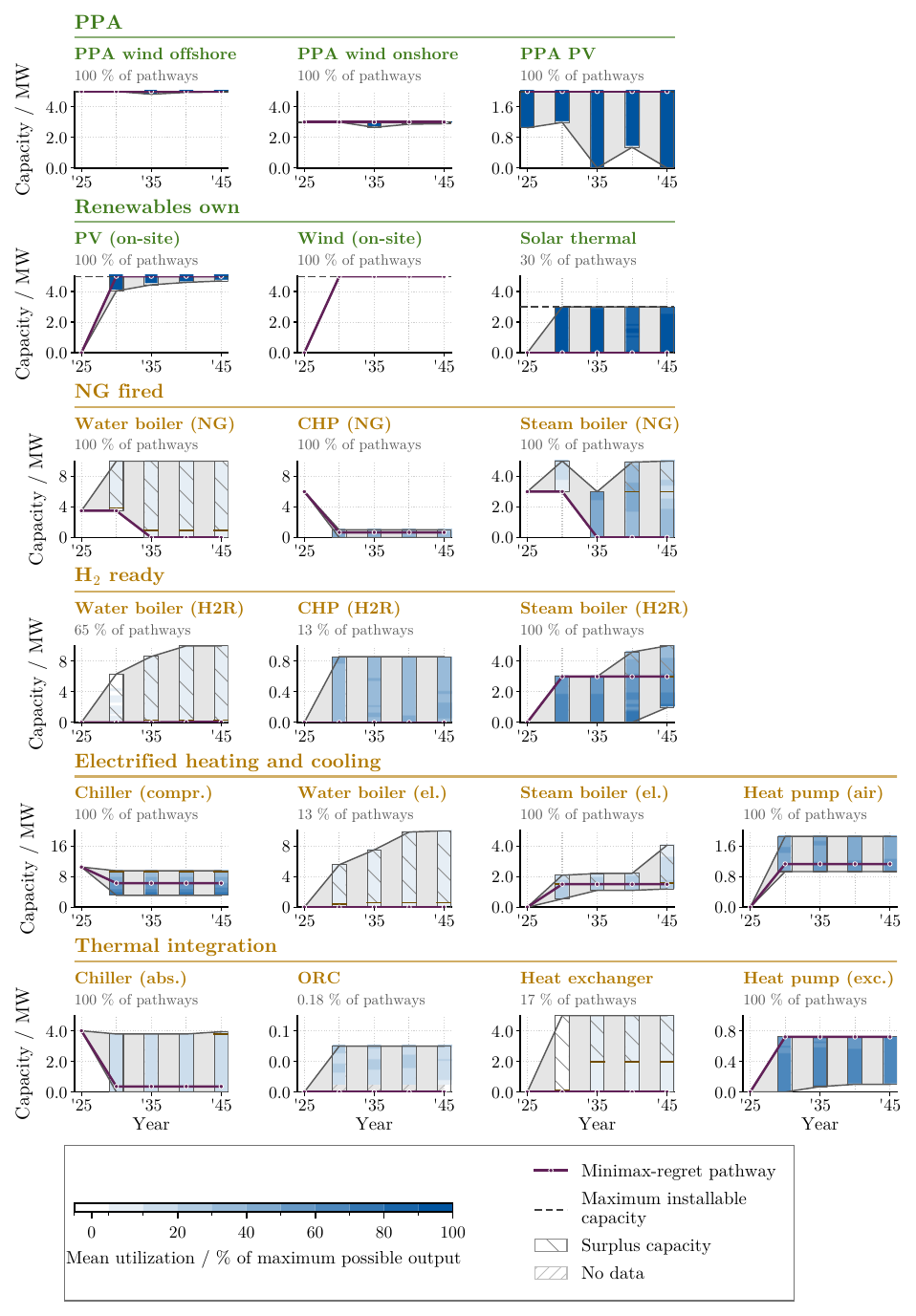}
\caption{Mean utilization of the conversion
components within the low-regret pathway space at the regret tolerance
$\bar r=2.5\%$. Utilization values are shown as colored columns at the five investment years, and the grey band connects the capacity range edges linearly. The dashed grey line marks the maximum installable capacity.}
\label{fig:si-tech-flh-a}
\end{figure}

\begin{figure}[H]
\centering
\includegraphics[width=\textwidth]{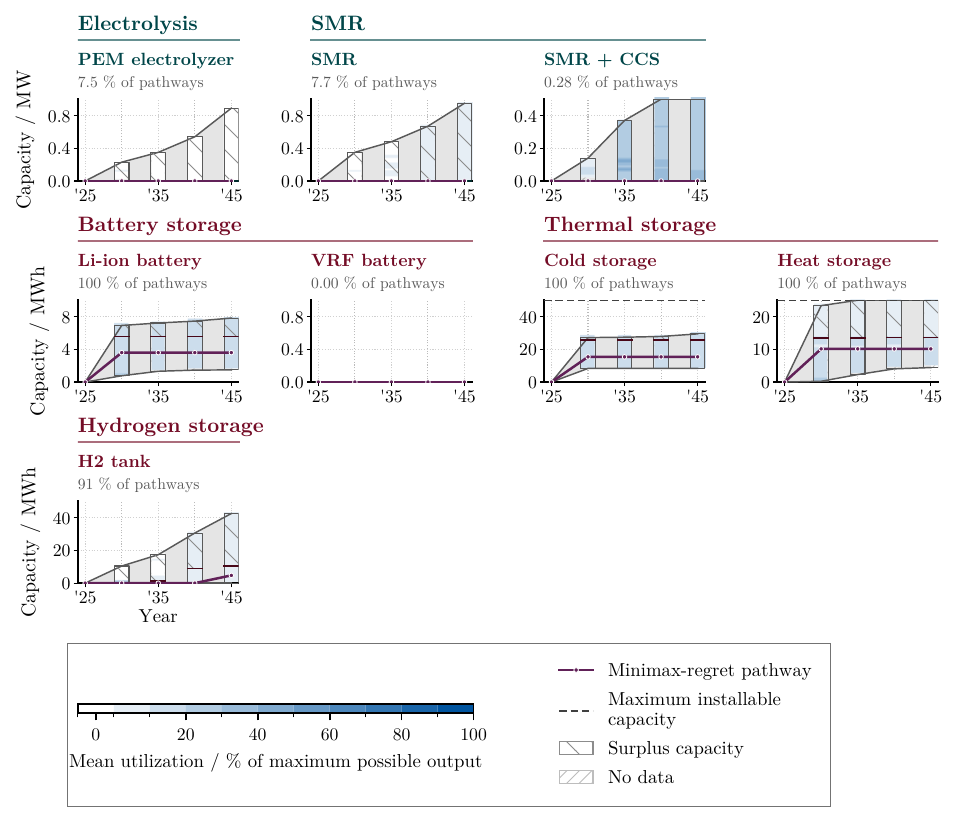}
\caption{Mean utilization of the hydrogen-production and storage
components within the low-regret pathway space at the regret tolerance
$\bar r=2.5\%$. Utilization values are shown as colored columns at the five investment years, and the grey band connects the capacity range edges linearly. The dashed grey line marks the maximum installable capacity.}
\label{fig:si-tech-flh-b}
\end{figure}

Figures~\ref{fig:si-tech-ranges-a} and~\ref{fig:si-tech-ranges-b} show the capacity ranges for the conversion and storage components. Figures~\ref{fig:si-tech-flh-a} and~\ref{fig:si-tech-flh-b} show the observed operation of the components, coloring each technology's capacity band at a regret tolerance of 2.5\% by its mean utilization relative to the maximum possible output. At a given capacity value, the color is the mean over the explored pathways that install approximately that capacity of the component's annual output, averaged across scenarios and divided by the output its installed capacity could deliver at full availability. That reference output is 8{,}760 full-load hours for conversion components with full availability. The maximum possible full-load hours for the renewables is about 1{,}200~h for photovoltaics, 2{,}400~h for onshore wind, 5{,}500~h for offshore wind, and 1{,}400~h for solar thermal. For storage components, we refer utilization to the maximum annual throughput permitted by the charge and discharge power limit for storage (Section~\ref{sec:si-stor-parameters}). The hatched part of each band is the surplus capacity of the component. The percentage above each panel is the share of explored pathways that install the component.

\newpage

\subsection{Scenario-wise regret of the explored pathways}
\label{sec:si-regret-scenarios}

\begin{figure}[H]
\centering
\includegraphics[width=\textwidth]{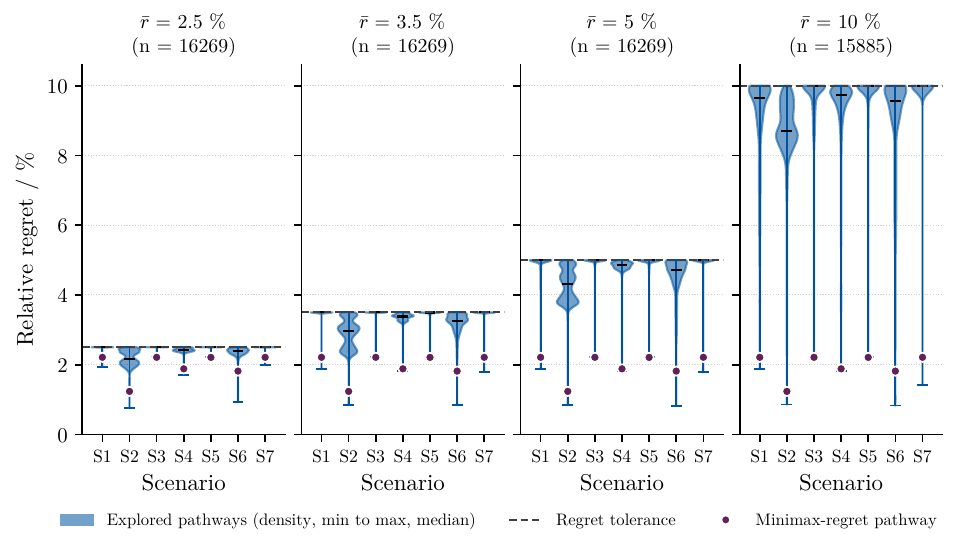}
\caption{Relative regret of the explored pathways per price scenario at the four regret tolerances $\bar r$. Each violin shows the distribution over all explored pathways of one low-regret pathway space from the minimum to the maximum, the black line marks the median. The dashed line is the regret tolerance, the purple markers are the minimax-regret pathway.}
\label{fig:si-regret-scenarios}
\end{figure}

Figure~\ref{fig:si-regret-scenarios} shows the relative regret of the explored pathways per price scenario. Each panel collects the explored pathways of one low-regret pathway space, and each violin shows the distribution of their regret in one scenario. The dashed line marks the regret tolerance, which no explored pathway exceeds in any scenario, and the markers show the minimax-regret pathway.

\clearpage
\begingroup
\sloppy
\bibliographystyle{elsarticle-harv-nourl}
\bibliography{references}
\endgroup